\documentclass[12pt]{article}
\usepackage{full page, 
amssymb, amscd, graphicx, 
}
    \title{{\bf Logarithmic
tensor category theory for generalized modules for a conformal
vertex algebra, I: Introduction and strongly graded
algebras and their generalized modules}}
    \author{Yi-Zhi Huang, James Lepowsky and Lin Zhang}
    \date{}

\begin{document}
    \bibliographystyle{alpha}
    \maketitle

    \newtheorem{rema}{Remark}[section]
    \newtheorem{propo}[rema]{Proposition}
    \newtheorem{theo}[rema]{Theorem}
   \newtheorem{defi}[rema]{Definition}
    \newtheorem{lemma}[rema]{Lemma}
    \newtheorem{corol}[rema]{Corollary}
     \newtheorem{exam}[rema]{Example}
\newtheorem{assum}[rema]{Assumption}
     \newtheorem{nota}[rema]{Notation}
        \newcommand{\ba}{\begin{array}}
        \newcommand{\ea}{\end{array}}
        \newcommand{\be}{\begin{equation}}
        \newcommand{\ee}{\end{equation}}
        \newcommand{\bea}{\begin{eqnarray}}
        \newcommand{\eea}{\end{eqnarray}}
        \newcommand{\nno}{\nonumber}
        \newcommand{\nn}{\nonumber\\}
        \newcommand{\lbar}{\bigg\vert}
        \newcommand{\p}{\partial}
        \newcommand{\dps}{\displaystyle}
        \newcommand{\bra}{\langle}
        \newcommand{\ket}{\rangle}
 \newcommand{\res}{\mbox{\rm Res}}
\newcommand{\wt}{\mbox{\rm wt}\;}
\newcommand{\swt}{\mbox{\scriptsize\rm wt}\;}
 \newcommand{\pf}{{\it Proof}\hspace{2ex}}
 \newcommand{\epf}{\hspace{2em}$\square$}
 \newcommand{\epfv}{\hspace{1em}$\square$\vspace{1em}}
        \newcommand{\ob}{{\rm ob}\,}
        \renewcommand{\hom}{{\rm Hom}}
\newcommand{\C}{\mathbb{C}}
\newcommand{\R}{\mathbb{R}}
\newcommand{\Z}{\mathbb{Z}}
\newcommand{\N}{\mathbb{N}}
\newcommand{\A}{\mathcal{A}}
\newcommand{\Y}{\mathcal{Y}}
\newcommand{\Arg}{\mbox{\rm Arg}\;}
\newcommand{\comp}{\mathrm{COMP}}
\newcommand{\lgr}{\mathrm{LGR}}

\newcommand{\dlt}[3]{#1 ^{-1}\delta \bigg( \frac{#2 #3 }{#1 }\bigg) }

\newcommand{\dlti}[3]{#1 \delta \bigg( \frac{#2 #3 }{#1 ^{-1}}\bigg) }

 \makeatletter
\newlength{\@pxlwd} \newlength{\@rulewd} \newlength{\@pxlht}
\catcode`.=\active \catcode`B=\active \catcode`:=\active \catcode`|=\active
\def\sprite#1(#2,#3)[#4,#5]{
   \edef\@sprbox{\expandafter\@cdr\string#1\@nil @box}
   \expandafter\newsavebox\csname\@sprbox\endcsname
   \edef#1{\expandafter\usebox\csname\@sprbox\endcsname}
   \expandafter\setbox\csname\@sprbox\endcsname =\hbox\bgroup
   \vbox\bgroup
  \catcode`.=\active\catcode`B=\active\catcode`:=\active\catcode`|=\active
      \@pxlwd=#4 \divide\@pxlwd by #3 \@rulewd=\@pxlwd
      \@pxlht=#5 \divide\@pxlht by #2
      \def .{\hskip \@pxlwd \ignorespaces}
      \def B{\@ifnextchar B{\advance\@rulewd by \@pxlwd}{\vrule
         height \@pxlht width \@rulewd depth 0 pt \@rulewd=\@pxlwd}}
      \def :{\hbox\bgroup\vrule height \@pxlht width 0pt depth
0pt\ignorespaces}
      \def |{\vrule height \@pxlht width 0pt depth 0pt\egroup
         \prevdepth= -1000 pt}
   }
\def\endsprite{\egroup\egroup}
\catcode`.=12 \catcode`B=11 \catcode`:=12 \catcode`|=12\relax
\makeatother

\def\hboxtr{\FormOfHboxtr} 
\sprite{\FormOfHboxtr}(25,25)[0.5 em, 1.2 ex] 

:BBBBBBBBBBBBBBBBBBBBBBBBB |
:BB......................B |
:B.B.....................B |
:B..B....................B |
:B...B...................B |
:B....B..................B |
:B.....B.................B |
:B......B................B |
:B.......B...............B |
:B........B..............B |
:B.........B.............B |
:B..........B............B |
:B...........B...........B |
:B............B..........B |
:B.............B.........B |
:B..............B........B |
:B...............B.......B |
:B................B......B |
:B.................B.....B |
:B..................B....B |
:B...................B...B |
:B....................B..B |
:B.....................B.B |
:B......................BB |
:BBBBBBBBBBBBBBBBBBBBBBBBB |

\endsprite

\def\shboxtr{\FormOfShboxtr} 
\sprite{\FormOfShboxtr}(25,25)[0.3 em, 0.72 ex] 

:BBBBBBBBBBBBBBBBBBBBBBBBB |
:BB......................B |
:B.B.....................B |
:B..B....................B |
:B...B...................B |
:B....B..................B |
:B.....B.................B |
:B......B................B |
:B.......B...............B |
:B........B..............B |
:B.........B.............B |
:B..........B............B |
:B...........B...........B |
:B............B..........B |
:B.............B.........B |
:B..............B........B |
:B...............B.......B |
:B................B......B |
:B.................B.....B |
:B..................B....B |
:B...................B...B |
:B....................B..B |
:B.....................B.B |
:B......................BB |
:BBBBBBBBBBBBBBBBBBBBBBBBB |

\endsprite


\begin{abstract}
This is the first part in a series of papers in which we introduce and
develop a natural, general tensor category theory for suitable module
categories for a vertex (operator) algebra.  This theory generalizes
the tensor category theory for modules for a vertex operator algebra
previously developed in a series of papers by the first two authors to
suitable module categories for a ``conformal vertex algebra'' or even
more generally, for a ``M\"obius vertex algebra.''  We do not require
the module categories to be semisimple, and we accommodate modules
with generalized weight spaces.  As in the earlier series of papers,
our tensor product functors depend on a complex variable, but in the
present generality, the logarithm of the complex variable is required;
the general representation theory of vertex operator algebras requires
logarithmic structure.  This work includes the complete proofs in the
present generality and can be read independently of the earlier series
of papers.  Since this is a new theory, we present it in detail,
including the necessary new foundational material.  In addition, with
a view toward anticipated applications, we develop and present the
various stages of the theory in the natural, general settings in which
the proofs hold, settings that are sometimes more general than what we
need for the main conclusions. In this paper (Part I), we give a
detailed overview of our theory, state our main results and introduce
the basic objects that we shall study in this work.  We include a
brief discussion of some of the recent applications of this theory,
and also a discussion of some recent literature.
\end{abstract}


\tableofcontents
\vspace{2em}

In this paper, Part I of a series of eight papers, we give a detailed
overview of logarithmic tensor category theory, state our main results
and introduce the basic objects that we shall study in this work.  We
include a brief discussion of some of the recent applications of this
theory, and also a discussion of some recent literature.  The
sections, equations, theorems and so on are numbered globally in the
series of papers rather than within each paper, so that for example
equation (a.b) is the b-th labeled equation in Section a, which is
contained in the paper indicated as follows: The present paper, Part
I, contains Sections 1 and 2.  In Part II \cite{HLZ2}, which contains
Section 3, we develop logarithmic formal calculus and study
logarithmic intertwining operators.  In Part III \cite{HLZ3}, which
contains Section 4, we introduce and study intertwining maps and
tensor product bifunctors.  In Part IV \cite{HLZ4}, which contains
Sections 5 and 6, we give constructions of the $P(z)$- and
$Q(z)$-tensor product bifunctors using what we call ``compatibility
conditions'' and certain other conditions.  In Part V \cite{HLZ5},
which contains Sections 7 and 8, we study products and iterates of
intertwining maps and of logarithmic intertwining operators and we
begin the development of our analytic approach.  In Part VI
\cite{HLZ6}, which contains Sections 9 and 10, we construct the
appropriate natural associativity isomorphisms between triple tensor
product functors.  In Part VII \cite{HLZ7}, which contains Section 11,
we give sufficient conditions for the existence of the associativity
isomorphisms.  In Part VIII \cite{HLZ8}, which contains Section 12, we
construct braided tensor category structure.

\renewcommand{\theequation}{\thesection.\arabic{equation}}
\renewcommand{\therema}{\thesection.\arabic{rema}}
\setcounter{equation}{0}
\setcounter{rema}{0}

\section{Introduction}

\subsection*{A brief description of the present work}

In the representation theory of many important algebraic structures,
such as Lie algebras, groups (or group algebras), commutative
associative algebras, Hopf algebras or quantum groups, tensor product
operations among modules play a central role.  They not only give new
modules from known ones, but they of course also provide a powerful
tool for studying modules.  More significantly, suitable categories of
modules for such algebras, equipped with tensor product operations and
appropriate natural isomorphisms, and so on, become symmetric or
braided tensor categories, and this tensor category structure is
always used, even when it is not explicitly discussed.

Vertex operator algebras, and more generally, vertex algebras, are a
fundamental class of algebraic structures whose extensive theory has
been developed and used in recent years to provide the means to
illuminate and to solve many problems in a wide variety of areas of
mathematics and theoretical physics.  In particular, the
representation theory of vertex (operator) algebras plays deep roles
in the construction and study of infinite-dimensional Lie algebra
representations, of structures linking sporadic finite simple groups
to string theory and to the theory of modular functions, and of knot
invariants and $3$-manifold invariants, in mathematics; and of
conformal field theory and string theory, in physics.

The present work is devoted to introducing and developing a natural,
general tensor category theory for suitable module categories for a
vertex (operator) algebra.  This tensor category theory, and
consequently, the representation theory, of vertex (operator)
algebras, is much, much more elaborate and more difficult than that of
Lie algebras, commutative associative algebras, Hopf algebras or
quantum groups.  In fact, the vertex-operator-algebraic analogues of
even the most elementary parts of the tensor product theory of an
algebra such as one of those are highly nontrivial, and the theory
needs to be developed with completely new ideas and strategies (and
with great care!).  The present theory was what we needed to carry out
in order to obtain the appropriate vertex-operator-algebraic analogue
of the following routine triviality in the representation theory of
(for example) Lie algebras: ``Given a Lie algebra $\mathfrak{g}$,
consider the symmetric tensor category of $\mathfrak{g}$-modules.''  A
vertex operator algebra ``wants to be'' the space of primitive
elements of a Hopf algebra (as is a Lie algebra, for example; this
immediately yields the tensor category of modules), but a vertex
operator algebra is {\it not} the space of primitive elements of any
Hopf algebra, and this is the beginning of why the problem of
constructing a tensor product theory and a tensor category theory of
modules for a vertex operator algebra was (and is) hard.  Yet it is at
least as important to have a theory of tensor products and tensor
categories of modules for a vertex operator algebra as it is in
classical theories such as Lie algebra theory (where such tensor
products and tensor categories of modules exist ``automatically'').

In Lie algebra theory (among other theories), many important module
categories are semisimple, that is, every module is completely
reducible, while on the other hand, many important module categories
are not.  Earlier, the first two authors developed a theory of braided
tensor categories for the module category of a what we call a
``finitely reductive'' vertex operator algebra satisfying certain
additional conditions; finite reductivity means that the module
category is semisimple and that certain finiteness conditions hold.
But it is just as natural and important to develop a theory for
non-semisimple module categories in vertex operator algebra theory as
it is in Lie algebra theory.  Also, in any one of the classical
theories such as Lie algebra theory, observing that there is a tensor
category of modules is just as easy for not-necessarily-semisimple
modules as it is for semisimple modules.  For these and many other
reasons, we considered it a natural problem to generalize the tensor
category theory for vertex operator algebras from the finitely
reductive case to the general case.

The present work accomplishes this goal, culminating in the
construction of a braided tensor category structure on a suitable
module category, not assumed semisimple, for a vertex (operator)
algebra.  It turns out the non-semisimplicity of modules is intimately
linked to the presence of logarithms in the basic ingredients of the
theory, beginning with intertwining operators among modules, and this
is why we call the present theory ``logarithmic tensor category
theory.''  We must in fact consider ``generalized modules''---
structures for which a certain basic operator has generalized
eigenvectors in addition to ordinary eigenvectors.  This basic
operator is contained in a natural copy of the three-dimensional
simple Lie algebra, which plays the role of the Lie algebra of the
group of M\"obius symmetries; this Lie algebra is in turn a subalgebra
of a natural copy of the Virasoro algebra, a central extension of a
Lie algebra of conformal symmetries.  In this work, we carry out our
theory for suitable categories of generalized modules for a
``conformal vertex algebra,'' which includes a copy of the Virasoro
algebra, and even more generally, for a ``M\"obius vertex algebra,''
which has the M\"obius symmetries but not all of the conformal
symmetries.  The present theory explicitly includes the earlier
finitely reductive theory as a special case; however, the present
theory is (necessarily) much more elaborate and subtle than the
finitely reductive theory.

In both the finitely reductive and the logarithmic generality, even
the construction of the tensor product (generalized) modules is
nontrivial; the correct tensor product module of two modules (when it
exists) is not at all based on the tensor product vector space of the
two underlying vector spaces.  Moreover, the construction of the
necessary natural associativity isomorphisms among triples of modules
is highly nontrivial.  While in classical tensor product theories the
natural associativity isomorphisms among triples of modules are given
by the usual trivial maps, in the tensor product theory of modules for
a vertex (operator) algebra, the corresponding statement is not at all
true, and indeed, there are not even any candidates for easy
associativity isomorphisms.  These and many related issues require the
present tensor product and tensor category theory to be elaborate.

A crucial discovery in the work of the first two authors in the
finitely reductive case was the existence of natural tensor products
of two or more elements in the algebraic completions of tensor product
modules.  All of the categorical structures and properties are
formulated, constructed and/or proved using tensor products of
elements.  In the finitely reductive case, tensor products of elements
were defined using intertwining operators (without logarithms).  In
order to develop the tensor category theory in the general setting of
the present work, it is again crucial to establish the existence of
tensor products of elements and to prove the fundamental properties of
these tensor product elements, and to do this, we are inevitably led
to the development of the theory of logarithmic intertwining
operators.

The structures of tensor product module, natural associativity
isomorphisms, and resulting braided tensor category structure
incorporating these, constructed in the present work, are assumed to
exist in a number of research works in mathematics and physics.  The
results in the present work allow one to remove assumptions of this
type.  We provide a mathematical foundation for such results and for
ongoing and future research involving the representation theory of
vertex (operator) algebras.

In fact, what we actually construct in this work is structure much
stronger than braided tensor category structure: The natural
associativity isomorphisms are constructed by means of a ``logarithmic
operator product expansion'' theorem for logarithmic intertwining
operators.  This logarithmic operator product expansion is in fact the
starting point of ``logarithmic conformal field theory,'' which has
been studied extensively by physicists as well as mathematicians.
Here, this logarithmic operator product expansion is established as a
mathematical theorem.

Moreover, our constructions and proofs in this work actually give what
the first two authors have called ``vertex-tensor-categorical
structure,'' in which the tensor product bifunctors depend crucially
on complex variables.  This structure is necessary for producing the
desired braided tensor category struture, through the use of the
tensor product elements and logarithmic operator product expansion
mentioned above, and our construction of braided tensor category
structure involves a ``limiting process'' in which the
complex-analytic information is ``forgotten'' and only the
``topological'' information associated with braided tensor category
structure is retained.  When we perform this specialization to the
``limiting case'' of braided tensor category structure, tensor
products of three or more elements are no longer defined.

The word ``algebra'' appears in the phrases ``vertex operator
algebra'' and ``vertex algebra,'' but beginning at the stage of the
theory where one must compose intertwining operators, or rather,
intertwining maps, among (generalized) modules, one must use analysis
as well as algebra, starting even from the definition of composition
of intertwining maps.  The kind of algebra on which the theory is
largely based, and which is needed throughout, is called ``formal
calculus,'' which we must in fact extensively develop in the course of
the work.  We must also enhance formal calculus with a great deal of
analytic reasoning, and the synthesized theory is no longer ``pure
algebra.''

This work includes the complete proofs in the present generality and
can be read independently of the first two authors' earlier series of
papers carrying out the finitely reductive theory.  Since this is a
new theory, we present it in detail, including the necessary new
foundational material.  In addition, we develop and present the
various stages of the theory in the natural, general settings in which
the proofs hold, settings that are sometimes more general than
what we need for the main conclusions.  This will allow for the future
use of the intermediate results in a variety of directions.

Later in the Introduction, we mention some of
the recent applications of the present theory, and we 
include a discussion of some recent literature. We state the main results 
of the present work at the end of the Introduction.

The main results presented here have been announced in \cite{HLZ}.

\paragraph{Acknowledgments}
We are grateful to Bertram Kostant for asking whether the work of the
first two authors described in \cite{tensorK} could (be generalized
to) recover Kazhdan-Lusztig's braided tensor category structure
\cite{KL1}--\cite{KL5}; the present work accomplishes this, as
discussed below, and in fact, this question was one of our main
initial motivations for embarking on this work.
We would like to thank Sasha Kirillov Jr., Masahiko Miyamoto, Kiyokazu
Nagatomo and Bin Zhang for their early interest in our work on
logarithmic tensor product theory, and Bin Zhang in particular for
inviting L.~Z. to lecture on this ongoing work at Stony Brook.
Y.-Z.~H. would also like to thank Paul Fendley for raising his
interest in logarithmic conformal field theory.  We are grateful to
the participants in our courses and seminars for their insightful
comments on this work, especially Geoffrey Buhl, William Cook, Robert McRae,
Thomas Robinson, Brian Vancil, and most especially, Shashank Kanade.  We 
also thank Ingo Runkel for asking us whether our theory satisfies
various conditions that he and his collaborators invoke in their work.
(It does.)  The authors gratefully acknowledge partial support {}from
NSF grants DMS-0070800 and DMS-0401302. Y.-Z.~H. is also grateful for
partial support {}from NSF grant PHY-0901237 and for the hospitality
of Institut des Hautes \'{E}tudes Scientifiques in the fall of 2007.

\subsection*{Introduction}

In a series of papers (\cite{tensorAnnounce}, \cite{tensorK},
\cite{tensor1}, \cite{tensor2}, \cite{tensor3}, \cite{tensor4}), the
first two authors have developed a tensor product and tensor category
theory for modules for a vertex operator algebra under suitable
conditions.  A structure called ``vertex tensor category structure''
(see \cite{tensorK}), which is much richer than tensor category
structure, has thereby been established for many important categories
of modules for classes of vertex operator algebras, since the
conditions needed for invoking the general theory have been verified
for these categories.  The most important such families of examples of
this theory are listed in Section 1.1 below.  In the present work,
which has been announced in \cite{HLZ}, we generalize this tensor
category theory to a larger family of module categories, for a
``conformal vertex algebra,'' or even more generally, for a ``M\"obius
vertex algebra,'' under suitably relaxed conditions.  A conformal
vertex algebra is just a vertex algebra in the sense of Borcherds
\cite{B} equipped with a conformal vector satisfying the usual axioms;
a M\"obius vertex algebra is a variant of a ``quasi-vertex operator
algebra'' as in \cite{FHL}.  Central features of the present work are
that we do not require the modules in our categories to be completely
reducible and that we accommodate modules with generalized weight
spaces.

As in the earlier series of papers, our tensor product functors depend
on a complex variable, but in the present generality, the logarithm of
the complex variable is required.  The first part of this work is
devoted to the study of logarithmic intertwining operators and their
role in the construction of the tensor product functors.  The
remainder of this work is devoted to the construction of the
appropriate natural associativity isomorphisms between triple tensor
product functors, to the proof of their fundamental properties, and to
the construction of the resulting braided tensor category structure.
This leads to vertex tensor category structure for further important
families of examples, or, in the M\"obius case, to ``M\"obius vertex
tensor category'' structure.

We emphasize that we develop our representation theory (tensor
category theory) in a very general setting; the vertex (operator)
algebras that we consider are very general, and the ``modules'' that
we consider are very general.  We call them ``generalized modules'';
they are not assumed completely reducible.  Many extremely important
(and well-understood) vertex operator algebras have semisimple module
categories, but in fact, now that the theory of vertex operator
algebras and of their representations is as highly developed as it has
come to be, it is in fact possible, and very fruitful, to work in the
greater generality.  Focusing mainly on the representation theory of
those vertex operator algebras for which every module is completely
reducible would be just as restrictive as focusing, classically, on
the representation theory of semisimple Lie algebras as opposed to the
representation theory of Lie algebras in general.  In addition, once
we consider suitably general vertex (operator) algebras, it is
unnatural to focus on only those modules that are completely
reducible.  As we explain below, such a general representation theory
of vertex (operator) algebras requires logarithmic structure.

A general representation theory of vertex operator algebras is crucial
for a range of applications, and we expect that it will be a
foundation for future developments.  One example is that the original
formulation of the uniqueness conjecture \cite{FLM2} for the moonshine
module vertex operator algebra $V^{\natural}$ (again see \cite{FLM2})
requires (general) vertex operator algebras whose modules might not be
completely reducible.  Another example is that this general theory is
playing a deep role in the (mathematical) construction of conformal
field theories (cf. \cite{HPNAS}, \cite{Hconference},
\cite{HVerlindeconjecture}, \cite{rigidity}, \cite{LPNAS}), which in
turn correspond to the perturbative part of string theory.  Just as
the classical (general) representation theory of groups, or of Lie
groups, or of Lie algebras, is not about any particular group or Lie
group or Lie algebra (although one of its central goals is certainly
to understand the representations of particular structures), the
general representation theory of suitably general vertex operator
algebras is ``background independent,'' in the terminology of string
theory.  In addition, the general representation theory of vertex
(operator) algebras can be thought of as a ``symmetry'' theory, where
vertex (operator) algebras play a role analogous to that of groups or
of Lie algebras in classical theories; deep and well-known analogies
between the notion of vertex operator algebra and the classical notion
of, for example, Lie algebra are discussed in several places,
including \cite{FLM2}, \cite{FHL} and \cite{LL}.

As we mentioned above, the present work includes the complete proofs
in the present generality and can be read independently of the earlier
series of papers of the first two authors constructing tensor
categories.  Our treatment is based on the theory of vertex operator
algebras and their modules as developed in \cite{B}, \cite{FLM2},
\cite{FHL}, \cite{DL} and \cite{LL}.  Throughout the work, we must,
and do, develop new algebraic and analytic methods, including a
synthesis of the ``formal calculus'' of vertex operator algebra theory
with analysis.

\subsection{Tensor category theory for finitely reductive 
vertex operator algebras}

The main families for which the conditions needed for invoking the
first two authors' general tensor category theory have been verified,
thus yielding vertex tensor category structure \cite{tensorK} on these
module categories, include the module categories for the following
classes of vertex operator algebras (or, in the last case, vertex
operator superalgebras):

\begin{enumerate}
\item The vertex operator algebras $V_L$ associated with positive
definite even lattices $L$; see \cite{B}, \cite{FLM2} for these vertex
operator algebras and see \cite{D1}, \cite{DL} for the conditions
needed for invoking the general tensor category theory.

\item The vertex operator algebras $L(k,0)$ associated with affine Lie
algebras and positive integers $k$; see \cite{FZ} for these vertex
operator algebras and \cite{FZ}, \cite{HLaffine} for the conditions.

\item The ``minimal series'' of vertex operator algebras associated
with the Virasoro algebra; see \cite{FZ} for these vertex operator
algebras and \cite{W}, \cite{H3} for the conditions.

\item Frenkel, Lepowsky and Meurman's moonshine module $V^{\natural}$;
see \cite{FLM1}, \cite{B}, \cite{FLM2} for this vertex operator
algebra and \cite{D2} for the conditions.

\item The fixed point vertex operator subalgebra of $V^{\natural}$
under the standard involution; see \cite{FLM1}, \cite{FLM2} for this
vertex operator algebra and \cite{D2}, \cite{H4} for the conditions.

\item The ``minimal series'' of vertex operator superalgebras
(suitably generalized vertex operator algebras) associated with the
Neveu-Schwarz superalgebra and also the ``unitary series'' of vertex
operator superalgebras associated with the $N=2$ superconformal
algebra; see \cite{KW} and \cite{A2} for the corresponding $N=1$ and
$N=2$ vertex operator superalgebras, respectively, and \cite{A1},
\cite{A3}, \cite{HM1}, \cite{HM2} for the conditions.
\end{enumerate}

In addition, vertex tensor category structure has also been
established for the module categories for certain vertex operator
algebras built {}from the vertex operator algebras just mentioned, such
as tensor products of such algebras; this is carried out in certain of
the papers listed above.

For all of the six classes of vertex operator algebras (or
superalgebras) listed above, each of the algebras is ``rational'' in
the specific sense of Huang-Lepowsky's work on tensor category theory.
This particular ``rationality'' property is easily proved to be a
sufficient condition for insuring that the tensor product modules
exist; see for instance \cite{tensor1}.  Unfortunately, the phrase
``rational vertex operator algebra'' also has several other distinct
meanings in the literature.  Thus we find it convenient at this time
to assign a new term, ``finite reductivity,'' to our particular notion
of ``rationality'': We say that a vertex operator algebra (or
superalgebra) $V$ is {\it finitely reductive} if:
\begin{enumerate}
\item Every $V$-module is completely reducible.
\item There are only finitely many irreducible $V$-modules (up to
equivalence).
\item All the fusion rules (the dimensions of the spaces of
intertwining operators among triples of modules) for $V$ are finite.
\end{enumerate}
We choose the term ``finitely reductive'' because we think of the term
``reductive'' as describing the complete reducibility---the first of
the conditions (that is, the algebra ``(completely) reduces'' every
module); the other two conditions are finiteness conditions.

The vertex-algebraic study of tensor category structure on module
categories for certain vertex algebras was stimulated by the work of
Moore and Seiberg \cite{MS1} \cite{MS}, in which, in the study of what they
termed ``rational'' conformal field theory, they obtained a set of polynomial 
equations  based on the assumption of the existence of a
suitable operator product expansion for ``chiral vertex operators''
(which correspond to intertwining operators in vertex algebra theory)
and observed an analogy between the theory of this set of polynomial equations
and the theory of tensor
categories.
Earlier, in \cite{BPZ}, Belavin, Polyakov, and Zamolodchikov had
already formalized the relation between the (nonmeromorphic) operator
product expansion, chiral correlation functions and representation
theory, for the Virasoro algebra in particular, and Knizhnik and
Zamolodchikov \cite{KZ} had established fundamental relations between
conformal field theory and the representation theory of affine Lie
algebras.  As we have discussed in the introductory material in
\cite{tensorK}, \cite{tensor1} and \cite{HLaffine}, such study of
conformal field theory is deeply connected with the vertex-algebraic
construction and study of tensor categories, and also with other
mathematical approaches to the construction of tensor categories in
the spirit of conformal field theory.  Concerning the latter
approaches, we would like to mention that the method used by Kazhdan
and Lusztig, especially in their construction of the associativity
isomorphisms, in their breakthrough work in \cite{KL1}--\cite{KL5}, is
related to the algebro-geometric formulation and study of
conformal-field-theoretic structures in the influential works of
Tsuchiya-Ueno-Yamada \cite{TUY}, Drinfeld \cite{Dr} and
Beilinson-Feigin-Mazur \cite{BFM}.  See also the important work of
Deligne \cite{De}, Finkelberg (\cite{F1} \cite{F2}), Bakalov-Kirillov
\cite{BK} and Nagatomo-Tsuchiya \cite{NT} on the construction of
tensor categories in the spirit of this approach to conformal field
theory, and also the discussions in Remark \ref{hist-btc} and in Section 
\ref{literature} below.

\subsection{Logarithmic tensor category theory}

The semisimplicity of the module categories mentioned in the examples
above is related to another property of these modules, namely, that
each module is a direct sum of its ``weight spaces,'' which are the
eigenspaces of a special operator $L(0)$ coming {}from the Virasoro
algebra action on the module.  But there are important situations in
which module categories are not semisimple and in which modules are
not direct sums of their weight spaces.  Notably, for the vertex
operator algebras $L(k,0)$ associated with affine Lie algebras, when
the sum of $k$ and the dual Coxeter number of the corresponding Lie
algebra is not a nonnegative rational number, the vertex operator
algebra $L(k,0)$ is not finitely reductive, and, working with Lie
algebra theory rather than with vertex operator algebra theory,
Kazhdan and Lusztig constructed a natural braided tensor category
structure on a certain category of modules of level $k$ for the affine
Lie algebra (\cite{KL1}, \cite{KL2}, \cite{KL3}, \cite{KL4},
\cite{KL5}).  This work of Kazhdan-Lusztig in fact motivated the first
two authors to develop an analogous theory for vertex operator
algebras rather than for affine Lie algebras, as was explained in
detail in the introductory material in \cite{tensorAnnounce},
\cite{tensorK}, \cite{tensor1}, \cite{tensor2}, and \cite{HLaffine}.
However, this general theory, in its original form, did not apply to
Kazhdan-Lusztig's context, because the vertex-operator-algebra modules
considered in \cite{tensorAnnounce}, \cite{tensorK}, \cite{tensor1},
\cite{tensor2}, \cite{tensor3}, \cite{tensor4} are assumed to be the
direct sums of their weight spaces (with respect to $L(0)$), and the
non-semisimple modules considered by Kazhdan-Lusztig fail in general
to be the direct sums of their weight spaces.  Although their setup,
based on Lie theory, and ours, based on vertex operator algebra
theory, are very different (as was discussed in the introductory
material in our earlier papers), we expected to be able to recover
(and further extend) their results through our vertex operator
algebraic approach, which is very general, as we discussed above.
This motivated us, in the present work, to generalize the work of the
first two authors by considering modules with generalized weight
spaces, and especially, intertwining operators associated with such
generalized kinds of modules.  As we discuss below, this required us
to use logarithmic intertwining operators and logarithmic formal
calculus, and we have been able to construct braided tensor category
structure, and even vertex tensor category structure, on important
module categories that are not semisimple.  Using the present theory,
the third author (\cite{Z1}, \cite{Z2}) has indeed recovered the
braided tensor category structure of Kazhdan-Lusztig, and has also
extended it to vertex tensor category structure.  While in our theory,
logarithmic structure plays a fundamental role, in this
Kazhdan-Lusztig work, logarithmic structure does not show up
explicitly.

{}From the viewpoint of the general representation theory of vertex
operator algebras, it would be unnatural to study only semisimple
modules or only $L(0)$-semisimple modules; focusing only on such
modules would be analogous to focusing only on semisimple modules for
general (nonsemisimple) finite-dimensional Lie algebras.  And as we
have pointed out, working in this generality leads to logarithmic
structure; the general representation theory of vertex operator
algebras requires logarithmic structure.

Logarithmic structure in conformal field theory was in fact first
introduced by physicists to describe Wess-Zumino-Witten models on
supergroups (\cite{RoS}, \cite{SS}) and disorder phenomena \cite{Gu}.  A
lot of progress has been made on this subject.  We refer the
interested reader to the review articles \cite{Ga}, \cite{Fl2},
\cite{RT} and \cite{Fu}, and references therein.  One particularly
interesting class of logarithmic conformal field theories is the class
associated to the triplet $\mathcal{W}$-algebras 
$\mathcal{W}(1,p)$ introduced by Kausch \cite{K1}, of central charge
$1-6\frac{(p-1)^{2}}{p}$, $p=2,3,\dots$.  We will discuss these algebras, and
generalizations of them, including references, in Section 1.5
below.  The paper \cite{FHST} initiated a study of a possible
generalization of the Verlinde conjecture for rational conformal field
theories to these theories; see also \cite{FG}, \cite{FK},
\cite{GR2} and \cite{GT}.  The paper \cite{Fu}
assumed the existence of braided tensor category structures on the
categories of modules for the vertex operator algebras considered;
together with \cite{H13}, the present work gives a construction of
these structures. The paper \cite{CF} used the results in the present
work as announced in \cite{HLZ}.

Here is how such logarithmic structure also arises naturally in the
representation theory of vertex operator algebras: In the construction
of intertwining operator algebras, the first author proved (see
\cite{diff-eqn}) that if modules for the vertex operator algebra
satisfy a certain cofiniteness condition, then products of the usual
intertwining operators satisfy certain systems of differential
equations with regular singular points. In addition, it was proved in
\cite{diff-eqn} that if the vertex operator algebra satisfies certain
finite reductivity conditions, then the analytic extensions of
products of the usual intertwining operators have no logarithmic
terms.  In the case when the vertex operator algebra satisfies only
the cofiniteness condition but not the finite reductivity conditions,
the products of intertwining operators still satisfy systems of
differential equations with regular singular points.  But in this
case, the analytic extensions of such products of intertwining
operators might have logarithmic terms. This means that if we want to
generalize the results in \cite{tensorAnnounce},
\cite{tensorK}--\cite{tensor3}, \cite{tensor4} and \cite{diff-eqn} to
the case in which the finite reductivity properties are not always
satisfied, we have to consider intertwining operators involving
logarithmic terms.

Logarithmic structure also appears naturally in modular invariance
results for vertex operator algebras and in the genus-one parts of
conformal field theories.  For a vertex operator algebra $V$
satisfying certain finiteness and reductivity conditions, Zhu proved
in \cite{Zhu} a modular invariance result for $q$-traces of products
of vertex operators associated to $V$-modules.  Zhu's result was
generalized to the case involving twisted vertex operators by Dong, Li
and Mason in \cite{DLM} and to the case of $q$-traces of products of
one intertwining operator and arbitrarily many vertex operators by
Miyamoto in \cite{M1}.  In \cite{M2}, Miyamoto generalized Zhu's
modular invariance result to a modular invariance result involving the
logarithm of $q$ for vertex operator algebras not necessarily
satisfying the reductivity condition.  In \cite{Hmodular}, for vertex
operator algebras satisfying certain finiteness and reductivity
conditions, by overcoming the difficulties one encounters if one tries
to generalize Zhu's method, the first author was able to prove the
modular invariance for $q$-traces of products and iterates of more
than one intertwining operator, using certain differential equations
and duality properties for intertwining operators.  If the vertex
operator algebra satisfies only Zhu's cofiniteness condition but not
the reductivity condition, the $q$-traces of products and iterates of
intertwining operators still satisfy the same differential equations,
but now they involve logarithms of all the variables.  To generalize
the general Verlinde conjecture proved in \cite{HVerlindeconjecture}
and the modular tensor category structure on the category of
$V$-modules obtained in \cite{rigidity}, one will need such general
logarithmic modular invariance.  See \cite{FHST}, \cite{Fu},
\cite{GR2} and \cite{GT} for research in this direction.

In \cite{Mi}, Milas introduced and studied what he called
``logarithmic modules'' and ``logarithmic intertwining operators.''
See also \cite{Mi2}.
Roughly speaking, logarithmic modules are weak modules for a vertex
operator algebra that are direct sums of generalized eigenspaces for
the operator $L(0)$.  We will call such weak modules ``generalized
modules'' in this work.  Logarithmic intertwining operators are
operators that depend not only on powers of a (formal or complex)
variable $x$, but also on its logarithm $\log x$.

The special features of the logarithm function make the logarithmic
theory very subtle and interesting.  In order to develop our
logarithmic tensor category theory, we were required to considerably
develop:
\begin{enumerate}
\item
Formal calculus, beyond what had been developed in \cite{FLM2},
\cite{FHL}, \cite{tensor1}--\cite{tensor3}, \cite{tensor4} and
\cite{LL}, in particular.  (Formal calculus has been developed in a
great many works.)
\item
What we may call ``logarithmic formal calculus,'' which involves
arbitrary powers of formal variables and of their formal logarithms.
This logarithmic formal calculus has been extended and exploited by
Robinson \cite{Ro1}, \cite{Ro2}, \cite{Ro3}.
\item
Complex analysis involving series containing {\it arbitrary real}
powers of the variables.
\item
Complex analysis involving series containing nonnegative integral
powers of the logarithms of the variables, in the presence of arbitrary 
real powers of the variables.
\item
A blending of these themes in order to formulate and to prove
many interchange-of-limit results necessary for the construction of
the ingredients of the logarithmic tensor category theory and for the
proofs of the fundamental properties.
\end{enumerate}
Our methods intricately combine both algebra and analysis, and must do
so, since the statements of the results themselves are both algebraic
and analytic.  See Remark \ref{newinloggenerality} below for a
discussion of these methods and their roles in this work.

As we mentioned above, one important application of our generalization
is to the category ${\cal O}_\kappa$ of certain modules for an affine
Lie algebra studied by Kazhdan and Lusztig in their series of papers
\cite{KL1}--\cite{KL5}. It has been shown in \cite{Z1} and \cite{Z2}
by the third author that, among other things, all the conditions
needed to apply our theory to this module category are satisfied.  As
a result, it is proved in \cite{Z1} and \cite{Z2} that there is a
natural vertex tensor category structure on this module category,
giving in particular a new construction, in the context of general
vertex operator algebra theory, of the braided tensor category
structure on ${\cal O}_\kappa$. This construction does not use 
the Knizhnik-Zamolodchikov equations. The methods used in
\cite{KL1}--\cite{KL5} were very different; the Knizhnik-Zamolodchikov
equations play an essential role in their construction, while the
present theory is very general.

The triplet $\mathcal{W}$-algebras belong to a different class of
vertex operator algebras, satisfying certain finiteness, boundedness
and reality conditions.  In this case, it has been shown in \cite{H13}
by the first author that all the conditions needed to apply the theory
carried out in the present work to the category of grading-restricted
modules for the vertex operator algebra are also satisfied.  Thus, by
the results obtained in this work, there is a natural vertex tensor
category structure on this category.

In addition to these logarithmic issues, another way in which the
present work generalizes the earlier tensor category theory for module
categories for a vertex operator algebra is that we now allow the
algebras to be somewhat more general than vertex operator algebras, in
order, for example, to accommodate module categories for the vertex
algebras $V_L$ where $L$ is a nondegenerate even lattice that is not
necessarily positive definite (cf.\ \cite{B}, \cite{DL}); see
\cite{Z1}.

What we accomplish in this work, then, is the following: We generalize
essentially all the results in \cite{tensor1}, \cite{tensor2},
\cite{tensor3} and \cite{tensor4} {}from the category of modules for a
vertex operator algebra to categories of suitably generalized modules
for a conformal vertex algebra or a M\"obius vertex algebra equipped
with an additional suitable grading by an abelian group.  The algebras
that we consider include not only vertex operator algebras but also
such vertex algebras as $V_L$ where $L$ is a nondegenerate even
lattice, and the modules that we consider are not required to be the
direct sums of their weight spaces but instead are required only to be
the (direct) sums of their ``generalized weight spaces,'' in a
suitable sense.  In particular, in this work we carry out, in the
present greater generality, the construction theory for the
``$P(z)$-tensor product'' functor originally done in \cite{tensor1},
\cite{tensor2} and \cite{tensor3} and the associativity theory for
this functor---the construction of the natural associativity
isomorphisms between suitable ``triple tensor products'' and the proof
of their important properties, including the isomorphism
property---originally done in \cite{tensor4}. This leads, as in
\cite{tensorK}, \cite{tensor5}, to the proof of the coherence
properties for vertex tensor categories, and in the M\"obius case, the
coherence properties for M\"obius vertex tensor categories.

For simplicity, we present our theory only for a conformal vertex
algebra or a M\"obius vertex algebra and not for their superalgebraic
analogues, but in fact our theory generalizes routinely to a conformal
vertex superalgebra or a M\"obius vertex superalgebra equipped with an
additional suitable grading by an abelian group; here we are referring
only to the usual sign changes associated with the ``odd'' subspace of
a vertex superalgebra, and not to any superconformal structure.

The general structure of much of this work essentially follows that of
\cite{tensor1}, \cite{tensor2}, \cite{tensor3} and \cite{tensor4}.
However, the results here are much stronger and much more general than in
these earlier works, and in addition, many of the results here have no
counterparts in those works.  Moreover, many ideas, formulations and
proofs in this work are necessarily quite different {}from those in the
earlier papers, and we have chosen to give some proofs that are new
even in the finitely reductive case studied in the earlier papers.

Some of the new ingredients that we are introducing into the theory
here are: an analysis of logarithmic intertwining operators, including
``logarithmic formal calculus''; a notion of
``$P(z_1,z_2)$-intertwining map'' and a study of its properties; new
``compatibility conditions''; considerable generalizations of
virtually all of the technical results in \cite{tensor1},
\cite{tensor2}, \cite{tensor3} and \cite{tensor4}; and perhaps most
significantly, the analytic ideas and methods that are sketched in
Remark \ref{newinloggenerality} below.

The contents of the sections of this work are as follows: In the rest
of this Introduction we compare classical tensor product and tensor
category theory for Lie algebra modules with tensor product and tensor
category theory for vertex operator algebra modules.  One crucial
difference between the two theories is that in the vertex operator
algebra setting, the theory depends on an ``extra parameter'' $z$,
which must be understood as a (nonzero) complex variable rather than
as a formal variable (although one needs, and indeed we very heavily
use, an extensive ``formal calculus,'' or ``calculus of formal
variables,'' in order to develop the theory).  We also discuss recent
applications of the present theory and some related literature and state 
the main results of the present work.  In
Section 2 we recall and extend some basic concepts in the theory of
vertex (operator) algebras.  We use the treatments of \cite{FLM2},
\cite{FHL}, \cite{DL} and \cite{LL}; in particular, the
formal-calculus approach developed in these works is needed for the
present theory.  Readers can consult these works for further detail.
We also set up notation and terminology that will be used throughout
the present work, and we describe the main categories of (generalized)
modules that we will consider. In Section 3 we introduce the notion of
logarithmic intertwining operator as in \cite{Mi} and present a
detailed study of the basic properties of these operators.  At the
beginning of this section we introduce and prove results about
logarithmic formal calculus, including a general ``formal Taylor
theorem.''  In Sections 4 and 5 we present the notions of $P(z)$- and
$Q(z)$-intertwining maps, and based on this, the definitions and
constructions of $P(z)$- and $Q(z)$-tensor products, generalizing
considerations in \cite{tensor1}, \cite{tensor2} and \cite{tensor3}.
The constructions of the tensor product functors require certain
``compatibility conditions'' and ``local grading restriction
conditions.''  The proofs of some of the results in Section 5 are
postponed to Section 6.  In Section 7 the convergence condition for
products and iterates of intertwining maps introduced in
\cite{tensor4} is generalized to the present context.  More
importantly, in this section we start to develop the complex analysis
approach that we will heavily use in later sections.  The new notion
of $P(z_1,z_2)$-intertwining map, generalizing the corresponding
concept in \cite{tensor4}, is introduced and developed in Section 8.
This will play a crucial role in the construction of the natural
associativity isomorphisms.  In Section 9 we prove important
conditions that are satisfied by vectors in the dual space of the
vector-space tensor product of three modules that arise {}from
products and {}from iterates of intertwining maps.  This leads us to
study elements in this dual space satisfying suitable compatibility
and local grading restriction conditions.  In this section we
extensively use our complex analysis approach, including, in
particular, for proving that the order of many iterated summations can
be interchanged.  By relating the subspaces considered in Section 9,
we construct the associativity isomorphisms in Section 10.  In Section
11, we generalize certain sufficient conditions for the existence of
the associativity isomorphisms in \cite{tensor4}, and we prove the
relevant conditions using differential equations.  In Section 12, we
establish the coherence properties of our braided tensor category
structure.

\subsection{The Lie algebra case}\label{LA}

In this section and the next, we compare classical tensor product and
tensor category theory for Lie algebra modules with the present theory
for vertex operator algebra modules, and in fact it is heuristically
useful to start by considering tensor product theory for Lie algebra
modules in a somewhat unusual way in order to motivate our approach
for the case of vertex (operator) algebras.

In the theory of tensor products for modules for a Lie algebra, the
tensor product of two modules is defined, or rather, constructed, as
the vector-space tensor product of the two modules, equipped with a
Lie algebra module action given by the familiar diagonal action of the
Lie algebra.  In the vertex algebra case, however, the vector-space
tensor product of two modules for a vertex algebra is {\it not} the
correct underlying vector space for the tensor product of the
vertex-algebra modules.  In this section we therefore consider
another approach to the tensor category theory for modules for a Lie
algebra---an approach, based on ``intertwining maps,'' that will show
how the theory proceeds in the vertex algebra case.  Then, in the next
section, we shall lay out the corresponding ``road map'' for the
tensor category theory in the vertex algebra case, which we then carry
out in the body of this work.

We first recall the following elementary but crucial background about
tensor product vector spaces: Given vector spaces $W_1$ and $W_2$, the
corresponding tensor product structure consists of a vector space $W_1
\otimes W_2$ equipped with a bilinear map
$$W_1 \times W_2 \longrightarrow W_1 \otimes W_2,$$
denoted
$$(w_{(1)},w_{(2)}) \mapsto w_{(1)}\otimes w_{(2)}$$
for $w_{(1)}\in W_1$ and $w_{(2)}\in W_2$, such that for any vector
space $W_3$ and any bilinear map
$$B:W_1 \times W_2 \longrightarrow W_3,$$
there is a unique linear map
$$L:W_1 \otimes W_2 \longrightarrow W_3$$
such that
$$B(w_{(1)},w_{(2)}) = L(w_{(1)}\otimes w_{(2)})$$
for $w_{(i)}\in W_i$, $i=1, 2$.  This universal property characterizes the
tensor product structure $W_1 \otimes W_2$, equipped with its bilinear
map $\cdot \otimes \cdot$, up to unique isomorphism.  In addition, the
tensor product structure in fact exists.

As was illustrated in \cite{tensorK}, and as is well known, the notion
of tensor product of modules for a Lie algebra can be formulated in
terms of what can be called ``intertwining maps'': Let $W_1$, $W_2$,
$W_3$ be modules for a fixed Lie algebra $V$.  (We are calling our Lie
algebra $V$ because we shall be calling our vertex algebra $V$, and we
would like to emphasize the analogies between the two theories.)  An
{\it intertwining map of type ${W_3 \choose {W_1 W_2}}$} is a linear
map $I: W_1 \otimes W_2 \longrightarrow W_3$ (or equivalently, {}from
what we have just recalled, a bilinear map $W_1 \times W_2
\longrightarrow W_3$) such that
\begin{equation}\label{intwmap}
\pi_3 (v)I(w_{(1)}\otimes w_{(2)})=I(\pi_1 (v)w_{(1)}\otimes
w_{(2)})+I(w_{(1)}\otimes \pi_2 (v)w_{(2)})
\end{equation}
for $v\in V$ and $w_{(i)}\in W_i$, $i=1,2$, where $\pi_1, \pi_2,
\pi_3$ are the module actions of $V$ on $W_1$, $W_2$ and $W_3$,
respectively.  (Clearly, such an intertwining map is the same as a
module map {}from $W_1\otimes W_2$, equipped with the tensor product
module structure, to $W_3$, but we are now temporarily ``forgetting''
what the tensor product module is.)

A {\it tensor product of the $V$-modules $W_1$ and $W_2$} is then a
pair $(W_0,I_0)$, where $W_0$ is a $V$-module and $I_0$ is an
intertwining map of type ${W_0 \choose {W_1 W_2}}$ (which, again,
could be viewed as a suitable bilinear map $W_1 \times W_2
\longrightarrow W_0$), such that for any pair $(W,I)$ with $W$ a
$V$-module and $I$ an intertwining map of type ${W \choose {W_1
W_2}}$, there is a unique module homomorphism $\eta:
W_0\longrightarrow W$ such that $I=\eta \circ I_0$. This universal
property of course characterizes $(W_0, I_0)$ up to canonical
isomorphism.  Moreover, it is obvious that the tensor product in fact
exists, and may be constructed as the vector-space tensor product
$W_1\otimes W_2$ equipped with the diagonal action of the Lie algebra,
together with the identity map {}from $W_1\otimes W_2$ to itself (or
equivalently, the canonical bilinear map $W_1 \times W_2
\longrightarrow W_1 \otimes W_2$).  We shall denote the tensor product
$(W_0,I_0)$ of $W_1$ and $W_2$ by $(W_1 \boxtimes W_2, \boxtimes)$,
where it is understood that the image of $w_{(1)}\otimes w_{(2)}$
under our canonical intertwining map $\boxtimes$ is $w_{(1)}\boxtimes
w_{(2)}$.  Thus $W_1 \boxtimes W_2=W_1\otimes W_2$, and under our
identifications, $\boxtimes = 1_{W_1\otimes W_2}$.

\begin{rema}
{\rm This classical explicit construction of course shows that the
tensor product functor exists for the category of modules for a Lie
algebra.  For vertex algebras, it will be relatively straightforward
to {\it define} the appropriate tensor product functor(s) (see
\cite{tensorK}, \cite{tensor1}, \cite{tensor2}, \cite{tensor3}), but
it will be a nontrivial matter to {\it construct} this functor (or
more precisely, these functors) and thereby prove that the
(appropriate) tensor product of modules for a (suitable) vertex
algebra exists.  The reason why we have formulated the notion of
tensor product module for a Lie algebra in the way that we just did is
that this formulation motivates the correct notion of tensor product
functor(s) in the vertex algebra case.}
\end{rema}

\begin{rema}
{\rm Using this explicit construction of the tensor product functor
and our notation $w_{(1)}\boxtimes w_{(2)}$ for the tensor product of
elements, the standard natural associativity isomorphisms among tensor
products of triples of Lie algebra modules are expressed as follows:
Since $w_{(1)}\boxtimes w_{(2)} = w_{(1)}\otimes w_{(2)}$, we have
\begin{eqnarray}\
(w_{(1)}\boxtimes w_{(2)})\boxtimes w_{(3)}&=&
(w_{(1)}\otimes w_{(2)})\otimes w_{(3)},\nno\\
w_{(1)}\boxtimes(w_{(2)}\boxtimes w_{(3)})&=&
w_{(1)}\otimes(w_{(2)}\otimes w_{(3)})\nno
\end{eqnarray}
for $w_{(i)}\in W_i$, $i=1,2,3$, and so
the canonical identification between
$w_{(1)}\otimes(w_{(2)}\otimes w_{(3)})$ and $(w_{(1)}\otimes
w_{(2)})\otimes w_{(3)}$ gives the standard natural isomorphism
\begin{eqnarray}
(W_1\boxtimes W_2)\boxtimes W_3 &\to & W_1\boxtimes (W_2\boxtimes W_3)\nno\\
(w_{(1)}\boxtimes w_{(2)})\boxtimes w_{(3)} &\mapsto &
w_{(1)}\boxtimes(w_{(2)}\boxtimes w_{(3)}).\label{elemap}
\end{eqnarray}
This collection of natural associativity isomorphisms of course
satisfies the classical coherence conditions for associativity
isomorphisms among multiple nested tensor product modules---the
conditions that say that in nested tensor products involving any
number of tensor factors, the placement of parentheses (as in
(\ref{elemap}), the case of three tensor factors) is immaterial; we
shall discuss coherence conditions in detail later.  Now, as was
discovered in \cite{tensor4}, it turns out that maps analogous to
(\ref{elemap}) can also be constructed in the vertex algebra case,
giving natural associativity isomorphisms among triples of modules for
a (suitable) vertex operator algebra.  However, in the vertex algebra
case, the elements ``$w_{(1)}\boxtimes w_{(2)}$,'' which indeed exist
(under suitable conditions) and are constructed in the theory, lie in
a suitable ``completion'' of the tensor product module rather than in
the module itself.  Correspondingly, it is a nontrivial matter to
construct the triple-tensor-product elements on the two sides of
(\ref{elemap}); in fact, one needs to prove certain convergence, under
suitable additional conditions.  Even after the triple-tensor-product
elements are constructed (in suitable completions of the
triple-tensor-product modules), it is a delicate matter to construct
the appropriate natural associativity maps, analogous to
(\ref{elemap}), to prove that they are well defined, and to prove that
they are isomorphisms.  In the present work, we shall generalize these
matters (in a self-contained way) {}from the context of \cite{tensor4}
to a more general one.  In the rest of this section, for triples of
modules for a Lie algebra, we shall now describe a construction of the
natural associativity isomorphisms that will seem roundabout and
indirect, but this is the method of construction of these isomorphisms
that will give us the correct ``road map'' for the corresponding
construction (and theorems) in the vertex algebra case, as in
\cite{tensor1}, \cite{tensor2}, \cite{tensor3} and \cite{tensor4}.}
\end{rema}

A significant feature of the constructions in the earlier works (and
in the present work) is that the tensor product of modules $W_1$ and
$W_2$ for a vertex operator algebra $V$ is the contragredient module
of a certain $V$-module that is typically a {\it proper} subspace of
$(W_1\otimes W_2)^*$, the dual space of the vector-space tensor
product of $W_1$ and $W_2$.  In particular, our treatment, which
follows, of the Lie algebra case will use contragredient modules, and
we will therefore restrict our attention to {\it finite-dimensional}
modules for our Lie algebra.  It will be important that the
double-contragredient module of a Lie algebra module is naturally
isomorphic to the original module.  We shall sometimes denote the
contragredient module of a $V$-module $W$ by $W'$, so that $W''=W$.
(We recall that for a module $W$ for a Lie algebra $V$, the
corresponding contragredient module $W'$ consists of the dual vector
space $W^*$ equipped with the action of $V$ given by: $(v \cdot
w^*)(w) = - w^*(v \cdot w)$ for $v \in V$, $w^* \in W^*$, $w \in W$.)

Let us, then, now restrict our attention to finite-dimensional modules
for our Lie algebra $V$.  The dual space $(W_1\otimes W_2)^*$ carries
the structure of the classical contragredient module of the tensor
product module.  Given any intertwining map of type ${W_3 \choose {W_1
W_2}}$, using the natural linear isomorphism
\begin{equation}\label{corpd2}
\hom (W_1\otimes W_2, W_3)\tilde{\longrightarrow} \hom (W^*_3,
(W_1\otimes W_2)^*)
\end{equation}
we have a corresponding linear map in $\hom (W^*_3, (W_1\otimes
W_2)^*)$, and this must be a map of $V$-modules.  In the vertex
algebra case, given $V$-modules $W_1$ and $W_2$, it turns out that
with a suitable analogous setup, the union in the vector space
$(W_1\otimes W_2)^*$ of the ranges of all such $V$-module maps, as
$W_3$ and the intertwining map vary (and with $W^*_3$ replaced by the
contragredient module $W'_3$), is the correct candidate for the
contragredient module of the tensor product module $W_1\boxtimes W_2$.
Of course, in the Lie algebra situation, this union is $(W_1\otimes
W_2)^*$ itself (since we are allowed to take $W_3=W_1\otimes W_2$ and
the intertwining map to be the canonical map), but in the vertex
algebra case, this union is typically much smaller than $(W_1\otimes
W_2)^*$.  In the vertex algebra case, we will use the notation $W_1
\hboxtr \, W_2$ to designate this union, and if the tensor product
module $W_1\boxtimes W_2$ in fact exists, then
\begin{eqnarray}
W_1 \boxtimes W_2 &=& (W_1 \,\hboxtr \; W_2)',\label{LAhbox1}\\
W_1 \,\hboxtr \; W_2 &=& (W_1 \boxtimes W_2)'.\label{LAhbox2}
\end{eqnarray}
Thus in the Lie algebra case we will write
\begin{eqnarray}
W_1 \,\hboxtr \; W_2 = (W_1 \otimes W_2)^*,
\end{eqnarray}
and (\ref{LAhbox1}) and (\ref{LAhbox2}) hold.  (In the Lie algebra
case we prefer to write $(W_1 \otimes W_2)^*$ rather than $(W_1
\otimes W_2)'$, because in the vertex algebra case, $W_1 \otimes W_2$
is typically not a $V$-module, and so we will not be allowed to write
$(W_1 \otimes W_2)'$ in the vertex algebra case.)

The subspace $W_1 \,\hboxtr \; W_2$ of $(W_1\otimes W_2)^*$ was in
fact further described in the following terms in \cite{tensor1} and
\cite{tensor3}, in the vertex algebra case: For any map in $\hom
(W'_3, (W_1\otimes W_2)^*)$ corresponding to an intertwining map
according to (\ref{corpd2}), the image of any $w'_{(3)}\in W'_3$ under
this map satisfies certain subtle conditions, called the
``compatibility condition'' and the ``local grading restriction
condition''; these conditions are not ``visible'' in the Lie algebra
case.  These conditions in fact precisely describe the proper subspace
$W_1 \hboxtr \, W_2$ of $(W_1\otimes W_2)^*$.  We will discuss such
conditions further in Section 1.4 and in the body of this work.  As we
shall explain, this idea of describing elements in certain dual spaces
was also used in constructing the natural associativity isomorphisms
between triples of modules for a vertex operator algebra in
\cite{tensor4}.

In order to give the reader a guide to the vertex algebra case, we now
describe the analogue for the Lie algebra case of this construction of
the associativity isomorphisms.  To construct the associativity
isomorphism {}from $(W_1 \boxtimes W_2)\boxtimes W_3$ to $W_1\boxtimes
(W_2\boxtimes W_3)$, it is equivalent (by duality) to give a suitable
isomorphism {}from $W_1\hboxtr\, (W_2\boxtimes W_3)$ to $(W_1 \boxtimes
W_2)\, \hboxtr\, W_3$ (recall (\ref{LAhbox1}), (\ref{LAhbox2})).

Rather than directly constructing an isomorphism between these two
$V$-modules, it turns out that we want to embed both of them,
separately, into the single space $(W_1\otimes W_2 \otimes W_3)^*$.
Note that $(W_1\otimes W_2 \otimes W_3)^*$
is naturally a $V$-module, via the contragredient of the
diagonal action, that is,
\begin{eqnarray}\label{actiononW1W2W2*}
(\pi(v)\lambda)(w_{(1)}\otimes w_{(2)}\otimes w_{(3)})&=&
-\lambda(\pi_1(v)w_{(1)}\otimes w_{(2)}\otimes w_{(3)})\nno\\
& -&\lambda(w_{(1)}\otimes \pi_2(v)w_{(2)}\otimes w_{(3)})\nno\\
& -&\lambda(w_{(1)}\otimes w_{(2)}\otimes \pi_3(v)w_{(3)}),
\end{eqnarray}
for $v\in V$ and $w_{(i)}\in W_i$, $i=1,2,3$, where $\pi_1, \pi_2,
\pi_3$ are the module actions of $V$ on $W_1$, $W_2$ and $W_3$,
respectively.  A concept related to this is the notion of
{\it intertwining map {}from $W_1 \otimes W_2 \otimes W_3$ to a module
$W_4$}, a natural analogue of (\ref{intwmap}), defined to be a linear map
\begin{equation}\label{intwmap3}
F: W_1\otimes W_2 \otimes W_3 \longrightarrow W_4
\end{equation}
such that
\begin{eqnarray}\label{intwmapfor3} \pi_4(v)F(w_{(1)}\otimes w_{(2)}
\otimes w_{(3)})&=&F(\pi_1(v)w_{(1)} \otimes w_{(2)} \otimes
w_{(3)})\nno\\ &+&F(w_{(1)} \otimes \pi_2(v)w_{(2)} \otimes
w_{(3)})\nno\\ &+&F(w_{(1)} \otimes w_{(2)} \otimes
\pi_{(3)}(v)w_3),\label{11i}
\end{eqnarray}
with the obvious notation.  The
relation between (\ref{actiononW1W2W2*}) and (\ref{intwmapfor3})
comes directly {}from the natural linear isomorphism
\begin{equation}
\hom (W_1\otimes W_2\otimes W_3, W_4) \tilde{\longrightarrow} \hom
(W^*_4, (W_1\otimes W_2\otimes W_3)^*);
\end{equation}
given $F$, we have
\begin{eqnarray}
W^*_4 &\longrightarrow &(W_1\otimes W_2\otimes
W_3)^*\nno\\
\nu&\mapsto &\nu\circ F.
\end{eqnarray}
Under this natural linear isomorphism, the intertwining maps
correspond precisely to the $V$-module maps {}from $W^*_4$ to
$(W_1\otimes W_2\otimes W_3)^*$.  In the situation for vertex
algebras, as was the case for tensor products of two rather than three
modules, there are analogues of all of the notions and comments
discussed in this paragraph {\it except that we will not put
$V$-module structure onto the vector space} $W_1\otimes W_2\otimes
W_3$; as we have emphasized, we will instead base the theory on
intertwining maps.

Two important ways of constructing maps of the type (\ref{intwmap3})
are as follows: For modules $W_1$, $W_2$, $W_3$, $W_4$, $M_1$ and
intertwining maps $I_1$ and $I_2$ of types ${W_4 \choose {W_1 M_1}}$
and ${M_1 \choose {W_2 W_3}}$, respectively, by definition the
composition $I_1\circ (1_{W_1}\otimes I_2)$ is an intertwining map
{}from $W_1\otimes W_2 \otimes W_3$ to $W_4$.  Analogously, for
intertwining maps $I^1$, $I^2$ of types ${W_4 \choose {M_2 W_3}}$ and
${M_2 \choose {W_1 W_2}}$, respectively, with $M_2$ also a module, the
composition $I^1\circ (I^2\otimes 1_{W_3})$ is an intertwining map
{}from $W_1\otimes W_2 \otimes W_3$ to $W_4$. Hence we have two
$V$-module homomorphisms
\begin{eqnarray}
W^*_4 &\longrightarrow &(W_1\otimes W_2\otimes
W_3)^*\nno\\
\nu&\mapsto &\nu\circ F_1,
\label{injint1}
\end{eqnarray}
where $F_1$ is the intertwining map $I_1\circ (1_{W_1}\otimes I_2)$;
and
\begin{eqnarray}
W^*_4 &\longrightarrow & (W_1\otimes W_2\otimes
W_3)^*\nno\\
\nu&\mapsto &\nu\circ F_2,
\label{injint2}
\end{eqnarray}
where $F_2$ is the intertwining map $I^1\circ (I^2\circ 1_{W_3})$.

The special cases in which the modules $W_4$ are two iterated tensor
product modules and the ``intermediate'' modules $M_1$ and $M_2$ are
two tensor product modules are particularly interesting: When
$W_4=W_1\boxtimes (W_2\boxtimes W_3)$ and $M_1=W_2\boxtimes W_3$, and
$I_1$ and $I_2$ are the corresponding canonical intertwining maps,
(\ref{injint1}) gives the natural $V$-module homomorphism
\begin{eqnarray}
W_1\,\hboxtr \,(W_2 \boxtimes W_3)&\longrightarrow &(W_1\otimes
W_2\otimes W_3)^*\nno\\
\nu&\mapsto &(w_{(1)}\otimes w_{(2)}\otimes w_{(3)}\mapsto \nu
(w_{(1)}\boxtimes (w_{(2)}\boxtimes w_{(3)})));\nno\\
\label{inj1}
\end{eqnarray}
when $W_4=(W_1\boxtimes W_2)\boxtimes W_3$ and $M_2=W_1\boxtimes W_2$,
and $I^1$ and $I^2$ are the corresponding canonical intertwining maps,
(\ref{injint2}) gives the natural $V$-module homomorphism
\begin{eqnarray}
(W_1 \boxtimes W_2)\,\hboxtr \,W_3&\longrightarrow & (W_1\otimes
W_2\otimes W_3)^*\nno\\
\nu&\mapsto &(w_{(1)}\otimes w_{(2)}\otimes w_{(3)}\mapsto \nu
((w_{(1)}\boxtimes w_{(2)})\boxtimes w_{(3)})).\nno\\
\label{inj2}
\end{eqnarray}

Clearly, in our Lie algebra case, both of the maps (\ref{inj1}) and
(\ref{inj2}) are isomorphisms, since they both in fact amount to the
identity map on $(W_1\otimes W_2\otimes W_3)^*$.  However, in the
vertex algebra case the analogues of these two maps are only injective
homomorphisms, and typically not isomorphisms.  (Recall the analogous
situation, mentioned above, for double rather than triple tensor
products.)  These two maps enable us to identify both $W_1\,\hboxtr
\,(W_2 \boxtimes W_3)$ and $(W_1 \boxtimes W_2)\,\hboxtr \,W_3$ with
subspaces of $(W_1\otimes W_2\otimes W_3)^*$.  In the vertex algebra
case we will have certain ``compatibility conditions'' and ``local
grading restriction conditions'' on elements of $(W_1\otimes
W_2\otimes W_3)^*$ to describe each of the two subspaces.  In either
the Lie algebra or the vertex algebra case, the construction of our
desired natural associativity isomorphism between the two modules
$(W_1 \boxtimes W_2)\boxtimes W_3$ and $W_1\boxtimes(W_2 \boxtimes
W_3)$ follows {}from showing that the ranges of homomorphisms
(\ref{inj1}) and (\ref{inj2}) are equal to each other, which is of
course obvious in the Lie algebra case since both (\ref{inj1}) and
(\ref{inj2}) are isomorphisms to $(W_1\otimes W_2\otimes W_3)^*$.  It
turns out that, under this associativity isomorphism, (\ref{elemap})
holds in both the Lie algebra case and the vertex algebra case; in the
Lie algebra case, this is obvious because all the maps are the
``tautological'' ones.

Now we give the reader a preview of how, in the vertex algebra case,
these compatibility and local grading restriction conditions on
elements of $(W_1\otimes W_2\otimes W_3)^*$ will arise.  As we have
mentioned, in the Lie algebra case, an intertwining map {}from
$W_1\otimes W_2 \otimes W_3$ to $W_4$ corresponds to a module map {}from
$W^*_4$ to $(W_1\otimes W_2\otimes W_3)^*$.  As was discussed in
\cite{tensor4}, for the vertex operator algebra analogue, the image of
any $w'_{(4)}\in W'_4$ under such an analogous map satisfies certain
``compatibility'' and ``local grading restriction'' conditions, and so
these conditions must be satisfied by those elements of $(W_1\otimes
W_2\otimes W_3)^*$ lying in the ranges of the vertex-operator-algebra
analogues of either of the maps (\ref{inj1}) and (\ref{inj2}) (or the
maps (\ref{injint1}) and (\ref{injint2})).

Besides these two conditions, satisfied by the elements of the ranges
of the maps of both types (\ref{inj1}) and (\ref{inj2}), the elements
of the ranges of the analogues of the homomorphisms (\ref{inj1}) and
(\ref{inj2}) have their own separate properties.  First note that any
$\lambda\in (W_1\otimes W_2\otimes W_3)^*$ induces the two maps
\begin{eqnarray}\label{mu1}
\mu^{(1)}_\lambda: W_1 &\to & (W_2\otimes W_3)^*\nno\\
 w_{(1)}&\mapsto & \lambda(w_{(1)} \otimes \cdot \otimes \cdot)
\end{eqnarray}
and
\begin{eqnarray}\label{mu2}
\mu^{(2)}_\lambda: W_3 &\to & (W_1\otimes W_2)^*\nno\\
w_{(3)}& \mapsto &\lambda(\cdot \otimes \cdot \otimes w_{(3)}).
\end{eqnarray}
In the vertex operator algebra analogue \cite{tensor4}, if $\lambda$
lies in the range of (\ref{inj1}), then it must satisfy the condition
that the elements $\mu^{(1)}_\lambda(w_{(1)})$ all lie in a suitable
completion of the subspace $W_2\hboxtr\, W_3$ of $(W_2\otimes W_3)^*$,
and if $\lambda$ lies in the range of (\ref{inj2}), then it must
satisfy the condition that the elements $\mu^{(2)}_\lambda(w_{(3)})$
all lie in a suitable completion of the subspace $W_1\hboxtr\, W_2$ of
$(W_1\otimes W_2)^*$.  (Of course in the Lie algebra case, these
statements are tautological.)  In \cite{tensor4}, these important
conditions, that $\mu^{(1)}_\lambda(W_1)$ lies in a suitable
completion of $W_2\hboxtr\, W_3$ and that $\mu^{(2)}_\lambda(W_3)$
lies in a suitable completion of $W_1\hboxtr\, W_2$, are understood as
``local grading restriction conditions'' with respect to the two
different ways of composing intertwining maps.

In the construction of our desired natural associativity isomorphism,
since we want the ranges of (\ref{inj1}) and (\ref{inj2}) to be the
same submodule of $(W_1\otimes W_2\otimes W_3)^*$, the ranges of both
(\ref{inj1}) and (\ref{inj2}) should satisfy both of these conditions.
This amounts to a certain ``expansion condition'' in the vertex
algebra case.  When all these conditions are satisfied, it can in fact
be proved \cite{tensor4} that the associativity isomorphism does
indeed exist and that in addition, the ``associativity of intertwining
maps'' holds; that is, the ``product'' of two suitable intertwining
maps can be written, in a certain sense, as the ``iterate'' of two
suitable intertwining maps, and conversely.  This equality of products
with iterates, highly nontrivial in the vertex algebra case, amounts
in the Lie algebra case to the easy statement that in the notation
above, any intertwining map of the form $I_1\circ (1_{W_1}\otimes
I_2)$ can also be written as an intertwining map of the form $I^1\circ
(I^2\otimes 1_{W_3})$, for a suitable ``intermediate module'' $M_2$
and suitable intertwining maps $I^1$ and $I^2$, and conversely.  The
reason why this statement is easy in the Lie algebra case is that in
fact {\it any} intertwining map $F$ of the type (\ref{intwmap3}) can
be ``factored'' in either of these two ways; for example, to write $F$
in the form $I_1\circ (1_{W_1}\otimes I_2)$, take $M_1$ to be
$W_2\otimes W_3$, $I_2$ to be the canonical (identity) map and $I_1$
to be $F$ itself (with the appropriate identifications having been
made).

We are now ready to discuss the vertex algebra case.

\subsection{The vertex algebra case}

In this section, which should be carefully compared with the
previous one, we shall lay out our ``road map'' of the constructions
of the tensor product functors and the associativity isomorphisms for
a suitable class of vertex algebras, considerably generalizing, but
also following the ideas of, the corresponding theory developed in
\cite{tensor1}, \cite{tensor2}, \cite{tensor3} and \cite{tensor4} for
vertex operator algebras.  Without yet specifying the precise class of
vertex algebras that we shall be using in the body of this work,
except to say that our vertex algebras will be $\mathbb{Z}$-graded and
our modules will be $\mathbb{C}$-graded at first and then
$\mathbb{R}$-graded for the more substantial results, we now discuss
the vertex algebra case.  What follows applies to both the theory of
\cite{tensor1}, \cite{tensor2}, \cite{tensor3}, \cite{tensor4} and the
present new logarithmic theory.  In Remark \ref{newinloggenerality}
below, we comment on the substantial new features of the logarithmic
generality.

In the vertex algebra case, the concept of intertwining map involves the moduli
space of Riemann spheres with one negatively oriented puncture and two
positively oriented punctures and with local coordinates around each
puncture; the details of the geometric structures needed in this
theory are presented in \cite{H0} and \cite{H1}.  For each element of this moduli
space there is a notion of intertwining map adapted to the particular
element.  Let $z$ be a nonzero complex number and let $P(z)$ be the
Riemann sphere $\hat{\mathbb C}$ with one negatively oriented puncture
at $\infty$ and two positively oriented punctures at $z$ and $0$, with
local coordinates $1/w$, $w-z$ and $w$ at these three punctures,
respectively.

Let $V$ be a vertex algebra (on which appropriate assumptions,
including the existence of a suitable ${\mathbb Z}$-grading, will be made
later), and let $Y(\cdot,x)$ be the vertex operator map defining the
algebra structure (see Section 2 below for a brief summary of basic
notions and notation, including the formal delta function).  Let
$W_1$, $W_2$ and $W_3$ be modules for $V$, and let $Y_1(\cdot,x)$,
$Y_2(\cdot,x)$ and $Y_3(\cdot,x)$ be the corresponding vertex operator
maps.  (The cases in which some of the $W_i$ are $V$ itself, and some
of the $Y_i$ are, correspondingly, $Y$, are important, but the most
interesting cases are those where all three modules are different {}from
$V$.)  A ``$P(z)$-intertwining map of type ${W_3 \choose {W_1 W_2}}$''
is a linear map
\begin{equation}
I: W_1 \otimes W_2 \longrightarrow \overline{W}_3,
\end{equation}
where $\overline{W}_3$ is a certain completion of $W_3$, related to
its ${\mathbb C}$-grading, such that
\begin{eqnarray}\label{im-jacobi}
\lefteqn{x_0^{-1}\delta\left(\frac{ x_1-z}{x_0}\right)
Y_3(v, x_1)I(w_{(1)}\otimes w_{(2)})}\nno\\
&&=z^{-1}\delta\left(\frac{x_1-x_0}{z}\right)
I(Y_1(v, x_0)w_{(1)}\otimes w_{(2)})\nno\\
&&\hspace{2em}+x_0^{-1}\delta\left(\frac{z-x_1}{-x_0}\right)
I(w_{(1)}\otimes Y_2(v, x_1)w_{(2)})
\end{eqnarray}
for $v\in V$, $w_{(1)}\in W_1$, $w_{(2)}\in W_2$, where $x_0$,
$x_1$ and $x_2$ are commuting independent formal variables.  This
notion is motivated in detail in \cite{tensorK}, \cite{tensor1} and
\cite{tensor3}; we shall recall the motivation below.

\begin{rema}\label{formalandcomplexvariables}
{\rm In this theory, it is crucial to distinguish between formal
variables and complex variables. Thus we shall use the following
notational convention: {\it Throughout this work, unless we specify
otherwise, the symbols $x$, $x_0$, $x_1$, $x_2$, $\dots$ , $y$, $y_0$,
$y_1$, $y_2$, $\dots$ will denote commuting independent formal
variables, and by contrast, the symbols $z$, $z_0$, $z_1$, $z_2$,
$\dots$ will denote complex numbers in specified domains, not formal
variables.}}
\end{rema}

\begin{rema}\label{im-io}{\rm
Recall {}from \cite{FHL} the definition of the notion of intertwining
operator ${\cal Y}(\cdot, x)$ in the theory of vertex (operator)
algebras.  Given $(W_1,Y_1)$, $(W_2,Y_2)$ and $(W_3,Y_3)$ as above, an
intertwining operator of type ${W_3 \choose {W_1 W_2}}$ can be viewed
as a certain type of linear map ${\cal Y}(\cdot,x)\cdot$ {}from $W_1
\otimes W_2$ to the vector space of formal series in $x$ of the form
$\sum_{n \in {\mathbb C}} w(n) x^n$, where the coefficients $w(n)$ lie in
$W_3$, and where we are allowing arbitrary complex powers of $x$,
suitably ``truncated {}from below'' in this sum.  The main property of an
intertwining operator is the following ``Jacobi identity'':
\begin{eqnarray}\label{io-jacobi}
\lefteqn{\dps x^{-1}_0\delta \left( {x_1-x_2\over x_0}\right)
Y_3(v,x_1){\cal Y}(w_{(1)},x_2)w_{(2)}}\nno\\
&&\hspace{2em}- x^{-1}_0\delta \left( {x_2-x_1\over -x_0}\right)
{\cal Y}(w_{(1)},x_2)Y_2(v,x_1)w_{(2)}\nno \\
&&{\dps = x^{-1}_2\delta \left( {x_1-x_0\over x_2}\right)
{\cal Y}(Y_1(v,x_0)w_{(1)},x_2)
w_{(2)}}
\end{eqnarray}
for $v\in V$, $w_{(1)}\in W_1$ and $w_{(2)}\in W_2$.  (When all three
modules $W_i$ are $V$ itself and all four operators $Y_i$ and ${\cal
Y}$ are $Y$ itself, (\ref{io-jacobi}) becomes the usual Jacobi
identity in the definition of the notion of vertex algebra.  When
$W_1$ is $V$, $W_2=W_3$ and ${\cal Y}=Y_2=Y_3$, (\ref{io-jacobi})
becomes the usual Jacobi identity in the definition of the notion of
$V$-module.)  The point is that by ``substituting $z$ for $x_2$'' in
(\ref{io-jacobi}), we obtain (\ref{im-jacobi}), where we make the
identification
\begin{eqnarray}\label{intwmap=intwopatz}
I(w_{(1)} \otimes w_{(2)}) = {\cal Y}(w_{(1)},z)w_{(2)};
\end{eqnarray}
the resulting complex powers of the complex number $z$ are made
precise by a choice of branch of the $\log$ function.  The nonzero
complex number $z$ in the notion of $P(z)$-intertwining map thus
``comes {}from'' the substitution of $z$ for $x_2$ in the Jacobi
identity in the definition of the notion of intertwining operator.  In
fact, this correspondence (given a choice of branch of $\log$)
actually defines an isomorphism between the space of
$P(z)$-intertwining maps and the space of intertwining operators of
the same type (\cite{tensor1}, \cite{tensor3}); this will be
discussed.}
\end{rema}

There is a natural linear injection
\begin{equation}\label{homva}
\hom (W_1\otimes W_2, \overline{W}_3)\longrightarrow
\hom (W'_3, (W_1\otimes W_2)^*),
\end{equation}
where here and below we denote by $W'$ the (suitably defined)
contragredient module of a $V$-module $W$; we have $W''=W$.  Under
this injection, a map $I\in \hom (W_1\otimes W_2, \overline{W}_3)$
amounts to a map $I':
W'_3\longrightarrow (W_1\otimes W_2)^*$:
\begin{equation}\label{I'}
w'_{(3)} \mapsto \langle w'_{(3)}, I(\cdot\otimes \cdot)\rangle,
\end{equation}
where $\langle \cdot,\cdot \rangle$ denotes the natural pairing
between the contragredient of a module and its completion. If $I$ is a
$P(z)$-intertwining map, then as in the Lie algebra case (see above),
where such a map is a module map, the map (\ref{I'}) intertwines two
natural $V$-actions on $W'_3$ and $(W_1\otimes W_2)^*$. We will see
that in the present (vertex algebra) case, $(W_1\otimes W_2)^*$ is
typically not a $V$-module.  The images of all the elements $w'_{(3)}\in
W'_3$ under this map satisfy certain conditions, called the
``$P(z)$-compatibility condition'' and the ``$P(z)$-local grading
restriction condition,'' as formulated in \cite{tensor1} and
\cite{tensor3}; we shall be discussing these.

Given a category of $V$-modules and two modules $W_1$ and $W_2$ in
this category, as in the Lie algebra case, the ``$P(z)$-tensor product
of $W_1$ and $W_2$'' is then defined to be a pair $(W_0,I_0)$, where
$W_0$ is a module in the category and $I_0$ is a $P(z)$-intertwining
map of type ${W_0 \choose {W_1 W_2}}$, such that for any pair $(W,I)$
with $W$ a module in the category and $I$ a $P(z)$-intertwining map of
type ${W \choose {W_1 W_2}}$, there is a unique morphism $\eta:
W_0\longrightarrow W$ such that $I=\bar\eta \circ I_0$; here and
throughout this work we denote by $\bar\chi$ the linear map naturally
extending a suitable linear map $\chi$ {}from a graded space to its
appropriate completion. This universal property characterizes $(W_0,
I_0)$ up to canonical isomorphism, {\it if it exists}.  We will denote
the $P(z)$-tensor product of $W_1$ and $W_2$, if it exists, by
$(W_1\boxtimes_{P(z)} W_2, \boxtimes_{P(z)})$, and we will denote the
image of $w_{(1)}\otimes w_{(2)}$ under $\boxtimes_{P(z)}$ by
$w_{(1)}\boxtimes_{P(z)} w_{(2)}$, which is an element of
$\overline{W_1\boxtimes_{P(z)} W_2}$, not of $W_1\boxtimes_{P(z)}
W_2$.

{}From this definition and the natural map (\ref{homva}), we will see
that if the $P(z)$-tensor product of $W_1$ and $W_2$ exists, then its
contragredient module can be realized as the union of ranges of all
maps of the form (\ref{I'}) as $W'_3$ and $I$ vary.  Even if the
$P(z)$-tensor product of $W_1$ and $W_2$ does not exist, we denote
this union (which is always a subspace stable under a natural action
of $V$) by $W_1\hboxtr_{P(z)} W_2$.  If the tensor product does exist,
then
\begin{eqnarray}
W_1 \boxtimes_{P(z)} W_2 &=& (W_1 \hboxtr_{P(z)}
W_2)',\label{vertexhbox1}\\
W_1 \hboxtr_{P(z)} W_2 &=& (W_1 \boxtimes_{P(z)}
W_2)';\label{vertexhbox2}
\end{eqnarray}
examining (\ref{vertexhbox1}) will show the reader why the notation
$\hboxtr$ was chosen in the earlier papers ($\boxtimes = \hboxtr
\,'$!).  Several critical facts about $W_1\hboxtr_{P(z)} W_2$ were
proved in \cite{tensor1}, \cite{tensor2} and \cite{tensor3}, notably,
$W_1\hboxtr_{P(z)} W_2$ is equal to the subspace of $(W_1\otimes
W_2)^*$ consisting of all the elements satisfying the
$P(z)$-compatibility condition and the $P(z)$-local grading
restriction condition, and in particular, this subspace is $V$-stable;
and the condition that $W_1\hboxtr_{P(z)} W_2$ is a module is
equivalent to the existence of the $P(z)$-tensor product
$W_1\boxtimes_{P(z)} W_2$.  All these facts will be proved.

In order to construct vertex tensor category structure, we need to
construct appropriate natural associativity isomorphisms.  Assuming
the existence of the relevant tensor products, we in fact need to
construct an appropriate natural isomorphism {}from $(W_1
\boxtimes_{P(z_1-z_2)} W_2)\boxtimes_{P(z_2)} W_3$ to
$W_1\boxtimes_{P(z_1)} (W_2\boxtimes_{P(z_2)} W_3)$ for complex
numbers $z_1$, $z_2$ satisfying $|z_1|>|z_2|>|z_1-z_2|>0$.  Note that
we are using two distinct nonzero complex numbers, and that certain
inequalities hold.  This situation corresponds to the fact that a
Riemann sphere with one negatively oriented puncture and three
positively oriented punctures can be seen in two different ways as the
``product'' of two Riemann spheres each with one negatively oriented
puncture and two positively oriented punctures; the detailed geometric
motivation is presented in \cite{H0}, \cite{H1}, \cite{tensorK} and
\cite{tensor4}.

To construct this natural isomorphism, we first consider compositions
of certain intertwining maps.  As we have mentioned, a
$P(z)$-intertwining map $I$ of type ${W_3 \choose {W_1 W_2}}$ maps
into $\overline{W}_3$ rather than $W_3$.  Thus the existence of
compositions of suitable intertwining maps always entails certain
convergence.  In particular, the existence of the composition
$w_{(1)}\boxtimes_{P(z_1)} (w_{(2)}\boxtimes_{P(z_2)} w_{(3)})$ when
$|z_1|>|z_2|>0$ and the existence of the composition
$(w_{(1)}\boxtimes_{P(z_1-z_2)} w_{(2)})\boxtimes_{P(z_2)} w_{(3)}$
when $|z_2|>|z_1-z_2|>0$, for general elements $w_{(i)}$ of $W_i$,
$i=1,2,3$, requires the proof of certain convergence conditions.
These conditions will be discussed in detail.

Let us now assume these convergence conditions and let $z_1$, $z_2$
satisfy $|z_1|>|z_2|>|z_1-z_2|>0$. To construct the desired
associativity isomorphism {}from $(W_1 \boxtimes_{P(z_1-z_2)}
W_2)\boxtimes_{P(z_2)} W_3$ to $W_1\boxtimes_{P(z_1)}
(W_2\boxtimes_{P(z_2)} W_3)$, it is equivalent (by duality) to give a
suitable natural isomorphism {}from $W_1\hboxtr_{P(z_1)}
(W_2\boxtimes_{P(z_2)} W_3)$ to $(W_1 \boxtimes_{P(z_1-z_2)}
W_2)\hboxtr_{P(z_2)} W_3$. As we mentioned in the previous section,
instead of constructing this isomorphism directly, we shall embed both
of these spaces, separately, into the single space $(W_1\otimes W_2
\otimes W_3)^*$.

We will see that $(W_1\otimes W_2 \otimes W_3)^*$ carries a natural
$V$-action analogous to the contragredient of the diagonal action in
the Lie algebra case (recall the similar action of $V$ on $(W_1\otimes
W_2)^*$ mentioned above).  Also, for four $V$-modules $W_1$, $W_2$,
$W_3$ and $W_4$, we have a canonical notion of ``$P(z_1,
z_2)$-intertwining map {}from $W_1 \otimes W_2 \otimes W_3$ to
$\overline{W}_4$'' given by a vertex-algebraic analogue of
(\ref{11i}); for this notion, we need only that $z_1$ and $z_2$ are
nonzero and distinct. The relation between these two concepts comes
{}from the natural linear injection
\begin{eqnarray}
\hom (W_1\otimes W_2\otimes W_3, \overline{W}_4) &\longrightarrow &
\hom (W'_4, (W_1\otimes W_2\otimes W_3)^*)\nno\\ F&\mapsto & F',
\end{eqnarray}
where $F': W'_4\longrightarrow (W_1\otimes W_2\otimes W_3)^*$ is given
by
\begin{equation}\label{F'}
\nu\mapsto \nu\circ F,
\end{equation}
which is indeed well defined.  Under this natural map, the
$P(z_1,z_2)$-intertwining maps correspond precisely to the maps {}from
$W'_4$ to $(W_1\otimes W_2\otimes W_3)^*$ that intertwine the two
natural $V$-actions on $W'_4$ and $(W_1\otimes W_2\otimes W_3)^*$.

Now for modules $W_1$, $W_2$, $W_3$, $W_4$, $M_1$, and a
$P(z_1)$-intertwining map $I_1$ and a $P(z_2)$-intertwining map $I_2$
of types ${W_4 \choose {W_1 M_1}}$ and ${M_1 \choose {W_2 W_3}}$,
respectively, it turns out that the composition $I_1\circ
(1_{W_1}\otimes I_2)$ exists and is a $P(z_1, z_2)$-intertwining map
when $|z_1|>|z_2|>0$.  Analogously, for a $P(z_2)$-intertwining map
$I^1$ and a $P(z_1-z_2)$-intertwining map $I^2$ of types ${W_4 \choose
{M_2 W_3}}$ and ${M_2 \choose {W_1 W_2}}$, respectively, where $M_2$
is also a module, the composition $I^1\circ (I^2\otimes 1_{W_3})$ is a
$P(z_1, z_2)$-intertwining map when $|z_2|>|z_1-z_2|>0$.  Hence we
have two maps intertwining the $V$-actions:
\begin{eqnarray}
W'_4 &\longrightarrow &(W_1\otimes W_2\otimes
W_3)^*\nno\\
\nu&\mapsto &\nu\circ F_1,
\label{iiv1}
\end{eqnarray}
where $F_1$ is the intertwining map $I_1\circ (1_{W_1}\otimes I_2)$,
and
\begin{eqnarray}
W'_4 &\longrightarrow & (W_1\otimes W_2\otimes
W_3)^*\nno\\
\nu&\mapsto &\nu\circ F_2,
\label{iiv2}
\end{eqnarray}
where $F_2$ is the intertwining map $I^1\circ (I^2\circ 1_{W_3})$.

It is important to note that we can express these compositions
$I_1\circ (1_{W_1}\otimes I_2)$ and $I^1\circ (I^2\otimes 1_{W_3})$ in
terms of intertwining operators, as discussed in Remark \ref{im-io}.
Let ${\cal Y}_1$, ${\cal Y}_2$, ${\cal Y}^1$ and ${\cal Y}^2$ be the
intertwining operators corresponding to $I_1$, $I_2$, $I^1$ and $I^2$,
respectively.  Then the compositions $I_1\circ (1_{W_1}\otimes I_2)$
and $I^1\circ (I^2\otimes 1_{W_3})$ correspond to the ``product''
${\cal Y}_1(\cdot, x_1){\cal Y}_2(\cdot, x_2)\cdot$ and ``iterate''
${\cal Y}^1({\cal Y}^2(\cdot, x_0)\cdot, x_2)\cdot$ of intertwining
operators, respectively, and we make the ``substitutions'' (in the
sense of Remark \ref{im-io}) $x_1 \mapsto z_1$, $x_2 \mapsto z_2$ and
$x_0 \mapsto z_1-z_2$ in order to express the two compositions of
intertwining maps as the ``product'' ${\cal Y}_1(\cdot, z_1){\cal
Y}_2(\cdot, z_2)\cdot$ and ``iterate'' ${\cal Y}^1({\cal Y}^2(\cdot,
z_1-z_2)\cdot, z_2)\cdot$ of intertwining maps, respectively.  (These
products and iterates involve a branch of the $\log$ function and also
certain convergence.)

Just as in the Lie algebra case, the special cases in which the
modules $W_4$ are two iterated tensor product modules and the
``intermediate'' modules $M_1$ and $M_2$ are two tensor product
modules are particularly interesting: When $W_4=W_1\boxtimes_{P(z_1)}
(W_2\boxtimes_{P(z_2)} W_3)$ and $M_1=W_2\boxtimes_{P(z_2)} W_3$, and
$I_1$ and $I_2$ are the corresponding canonical intertwining maps,
(\ref{iiv1}) gives the natural $V$-homomorphism
\begin{eqnarray}
W_1\hboxtr_{P(z_1)}(W_2 \boxtimes_{P(z_2)} W_3)&\longrightarrow &
(W_1\otimes W_2\otimes W_3)^*\nno\\
\nu&\mapsto &(w_{(1)}\otimes w_{(2)}\otimes w_{(3)}\mapsto \nno\\
&&\nu(w_{(1)}\boxtimes_{P(z_1)} (w_{(2)}\boxtimes_{P(z_2)}
w_{(3)})))\nno;\\
\label{injva1}
\end{eqnarray}
when $W_4=(W_1\boxtimes_{P(z_1-z_2)} W_2)\boxtimes_{P(z_2)} W_3$ and
$M_2=W_1\boxtimes_{P(z_1-z_2)}W_2$, and $I^1$ and $I^2$ are the
corresponding canonical intertwining maps, (\ref{iiv2}) gives the
natural $V$-homomorphism
\begin{eqnarray}
(W_1 \boxtimes_{P(z_1-z_2)} W_2)\hboxtr_{P(z_2)} W_3&\longrightarrow &
(W_1\otimes W_2\otimes
W_3)^*\nno\\
\nu&\mapsto &(w_{(1)}\otimes w_{(2)}\otimes w_{(3)}\mapsto \nno\\
&&\nu((w_{(1)}\boxtimes_{P(z_1-z_2)} w_{(2)})\boxtimes_{P(z_2)}
w_{(3)})).\nno\\
\label{injva2}
\end{eqnarray}

It turns out that both of these maps are injections, as in
\cite{tensor4} (as we shall prove), so that we are embedding the
spaces $W_1\hboxtr_{P(z_1)}(W_2 \boxtimes_{P(z_2)} W_3)$ and $(W_1
\boxtimes_{P(z_1-z_2)} W_2)\hboxtr_{P(z_2)} W_3$ into the space
$(W_1\otimes W_2 \otimes W_3)^*$. Following the ideas in
\cite{tensor4}, we shall give a precise description of the ranges of
these two maps, and under suitable conditions, prove that the two
ranges are the same; this will establish the associativity
isomorphism.

More precisely, as in \cite{tensor4}, we prove that for any $P(z_1,
z_2)$-intertwining map $F$, the image of any $\nu\in W'_4$ under $F'$
(recall (\ref{F'})) satisfies certain conditions that we call the
``$P(z_1, z_2)$-compatibility condition'' and the ``$P(z_1,
z_2)$-local grading restriction condition.''  Hence, as special cases,
the elements of $(W_1\otimes W_2 \otimes W_3)^*$ in the ranges of
either of the maps (\ref{iiv1}) or (\ref{iiv2}), and in particular, of
(\ref{injva1}) or (\ref{injva2}), satisfy these conditions.

In addition, any $\lambda\in (W_1\otimes W_2\otimes W_3)^*$ induces
two maps $\mu^{(1)}_\lambda$ and $\mu^{(2)}_\lambda$ as in (\ref{mu1})
and (\ref{mu2}).  We will see that any element $\lambda$ of the range
of (\ref{iiv1}), and in particular, of (\ref{injva1}), must satisfy
the condition that the elements $\mu^{(1)}_\lambda(w_{(1)})$ all lie,
roughly speaking, in a suitable completion of the subspace
$W_2\hboxtr_{P(z_2)} W_3$ of $(W_2\otimes W_3)^*$, and any element
$\lambda$ of the range of (\ref{iiv2}), and in particular, of
(\ref{injva2}), must satisfy the condition that the elements
$\mu^{(2)}_\lambda(w_{(3)})$ all lie, again roughly speaking, in a
suitable completion of the subspace $W_1\hboxtr_{P(z_1-z_2)} W_2$ of
$(W_1\otimes W_2)^*$.  These conditions will be called the
``$P^{(1)}(z)$-local grading restriction condition'' and the
``$P^{(2)}(z)$-local grading restriction condition,'' respectively.

It turns out that the construction of the desired natural
associativity isomorphism follows {}from showing that the ranges of
both of (\ref{injva1}) and (\ref{injva2}) satisfy both of these
conditions.  This amounts to a certain ``expansion condition'' on our
module category.  When this expansion condition and a suitable
convergence condition are satisfied, we show that the desired
associativity isomorphisms do exist, and that in addition, the
associativity of intertwining maps holds.  That is, let $z_1$ and
$z_2$ be complex numbers satisfying the inequalities
$|z_1|>|z_2|>|z_1-z_2|>0$.  Then for any $P(z_1)$-intertwining map
$I_1$ and $P(z_2)$-intertwining map $I_2$ of types ${W_4\choose {W_1\,
M_1}}$ and ${M_1\choose {W_2\, W_3}}$, respectively, there is a
suitable module $M_2$, and a $P(z_2)$-intertwining map $I^1$ and a
$P(z_1-z_2)$-intertwining map $I^2$ of types ${W_4\choose{M_2\, W_3}}$
and ${M_2\choose{W_1\, W_2}}$, respectively, such that
\begin{equation}\label{iiii}
\langle w'_{(4)}, I_1(w_{(1)}\otimes I_2(w_{(2)}\otimes w_{(3)}))
\rangle =
\langle w'_{(4)}, I^1(I^2(w_{(1)}\otimes w_{(2)})\otimes w_{(3)})
\rangle
\end{equation}
for $w_{(1)}\in W_1, w_{(2)}\in W_2$, $w_{(3)} \in W_3$ and
$w'_{(4)}\in W'_4$; and conversely, given $I^1$ and $I^2$ as
indicated, there exist a suitable module $M_1$ and maps $I_1$ and
$I_2$ with the indicated properties.  In terms of intertwining
operators (recall the comments above), the equality (\ref{iiii}) reads
\begin{eqnarray}\label{yyyy}
\lefteqn{\langle w'_{(4)}, {\cal Y}_1(w_{(1)}, x_1){\cal
Y}_2(w_{(2)},x_2)
w_{(3)}\rangle|_{x_1=z_1,\; x_2=z_2}}\nno\\
&&=\langle w'_{(4)},{\cal Y}^1({\cal Y}^2(w_{(1)}, x_0)w_{(2)},
x_2)w_{(3)})\rangle |_{x_0=z_1-z_2,\; x_2=z_2},
\end{eqnarray}
where ${\cal Y}_1$, ${\cal Y}_2$, ${\cal Y}^1$ and ${\cal Y}^2$ are
the intertwining operators corresponding to $I_1$, $I_2$, $I^1$ and
$I^2$, respectively.  (As we have been mentioning, the substitution of
complex numbers for formal variables involves a branch of the $\log$
function and also certain convergence.)  In this sense, the
associativity asserts that the ``product'' of two suitable
intertwining maps can be written as the ``iterate'' of two suitable
intertwining maps, and conversely.

{}From this construction of the natural associativity isomorphisms we
will see, by analogy with (\ref{elemap}), that
$(w_{(1)}\boxtimes_{P(z_1-z_2)} w_{(2)})\boxtimes_{P(z_2)} w_{(3)}$ is
mapped naturally to $w_{(1)}\boxtimes_{P(z_1)}
(w_{(2)}\boxtimes_{P(z_2)} w_{(3)})$ under the natural extension of
the corresponding associativity isomorphism (these elements in general
lying in the algebraic completions of the corresponding tensor product
modules). In fact, this property
\begin{equation}
(w_{(1)}\boxtimes_{P(z_1-z_2)} w_{(2)})\boxtimes_{P(z_2)} w_{(3)}
\mapsto
w_{(1)}\boxtimes_{P(z_1)}
(w_{(2)}\boxtimes_{P(z_2)} w_{(3)})
\end{equation}
for $w_{(1)}\in W_{1}$, $w_{(2)}\in W_{2}$ and $w_{(3)}\in W_{3}$
characterizes
the associativity isomorphism
\begin{equation}
(W_{1}\boxtimes_{P(z_1-z_2)} W_{2})\boxtimes_{P(z_2)} W_{3}\to
W_{1}\boxtimes_{P(z_1)}
(W_{2}\boxtimes_{P(z_2)} W_{3})
\end{equation}
(cf. (\ref{elemap})).  The coherence property of the associativity
isomorphisms will follow {}from this fact.  We will of course have
mutually inverse associativity isomorphisms.

\begin{rema}{\rm
Note that equation (\ref{yyyy}) can be written as
\begin{equation}\label{yyyy2}
{\cal Y}_1(w_{(1)}, z_1){\cal Y}_2(w_{(2)},z_2)=
{\cal Y}^1({\cal Y}^2(w_{(1)}, z_1-z_2)w_{(2)},
z_2),
\end{equation}
with the appearance of the complex numbers being understood as
substitutions in the sense mentioned above, and with the ``generic''
vectors $w_{(3)}$ and $w'_{(4)}$ being implicit.  This (rigorous)
equation amounts to the ``operator product expansion'' in the physics
literature on conformal field theory; indeed, in our language, if we
expand the right-hand side of (\ref{yyyy2}) in powers of $z_1-z_2$, we
find that a product of intertwining maps is expressed as an expansion
in powers of $z_1-z_2$, with coefficients that are again intertwining
maps, of the form ${\cal Y}^1(w,z_2)$.  When all three modules are the
vertex algebra itself, and all the intertwining operators are the
canonical vertex operator $Y(\cdot,x)$ itself, this ``operator product
expansion'' follows easily {}from the Jacobi identity.  But for
intertwining operators in general, it is a deep matter to prove the
operator product expansion, that is, to prove the assertions involving
(\ref{iiii}) and (\ref{yyyy}) above.  This was proved in
\cite{tensor4} in the finitely reductive setting and is considerably
generalized in the present work to the logarithmic setting.}
\end{rema}

\begin{rema}{\rm
The constructions of the tensor product modules and of the
associativity isomorphisms previewed above for suitably general vertex
algebras follow those in \cite{tensor1}, \cite{tensor2},
\cite{tensor3} and \cite{tensor4}.  Alternative constructions are
certainly possible.  For example, an alternative construction of the
tensor product modules was given in \cite{Li}.  However, no matter
what construction is used for the tensor product modules of suitably
general vertex algebras, one cannot avoid constructing structures and
proving results equivalent to what is carried out in this work.  The
constructions in this work of the tensor product functors and of the
natural associativity isomorphisms are crucial in the deeper part of
the theory of vertex tensor categories. }
\end{rema}

\begin{rema}\label{newinloggenerality}{\rm We have outlined the construction 
of the tensor product functors and the associativity isomorphisms
without getting into the technical details. On the other hand, though
the general ideas of the constructions are the same for both the
semisimple theory developed in \cite{tensor1}, \cite{tensor2},
\cite{tensor3} and \cite{tensor4} and the nonsemisimple logarithmic
theory carried out in the present work, many of the proofs of the
results in the present work involve substantial new ideas and
techniques, making the nonsemisimple logarithmic theory vastly more
difficult technically than the semisimple theory.  First, we have had
to further develop formal calculus beyond what had been developed in
\cite{FLM2}, \cite{FHL}, \cite{tensor1}--\cite{tensor3},
\cite{tensor4}, \cite{LL} and many other works.  We have had to study
new kinds of combinations of formal delta function expressions in
several {\it formal and complex} variables.  Second, we have extended
formal calculus to include logarithms of formal variables.  In formal
calculus, logarithms of formal variables are in fact additional
independent formal variables.  We develop our ``logarithmic formal
calculus'' in a much more general setting than what we need for the
main results in this work.  In particular, we at first allow the
formal series in a formal variable and its logarithm to involve {\it
infinitely many arbitrary complex powers} of the logarithm.  This
study of logarithmic formal calculus has surprising connections with
various classes of combinatorial identities and has been extended
and exploited by Robinson \cite{Ro1}, \cite{Ro2}, \cite{Ro3}.  Third,
to construct the natural associativity isomorphisms and other data for
the tensor categories and to prove the coherence property, it is
necessary to use complex analysis.  We wanted to carry out our theory
under the most general natural sets of assumptions that would indeed
{\it yield} a theory.  This required us to work with series involving
{\it arbitrary real} powers of the complex variables, {\it with the
powers not even being lower bounded}.  We have in fact extended a
number of classical results in complex analysis to results that can be
applied to such series.  In particular, we have had to prove many
results that allow us to switch orders of infinite sums, by either
proving the multiconvergence of the corresponding multisums or by
using Taylor expansion for analytic functions.  Fourth, since our
theory also involves logarithms of complex variables, we have also had
to extend those same classical results in complex analysis to results
that can be applied still further to series involving logarithms of
the complex variables.  In particular, we prove that when the powers
of the logarithm of a complex variable are bounded above in a series
involving {\it arbitrary real powers} of the variable and nonnegative
integral powers of its logarithm, the convergence of suitable {\it
iterated} sums implies absolute convergence of the corresponding {\it
double} sums.  We also prove what we call the ``unique expansion
property'' for the set $\R\times \{0, \dots, N\}$ (see Proposition
\ref{real-exp-set}), which says that the coefficients of an absolutely
convergent series of the form just indicated are determined uniquely
by its sum.  One important difference from the logarithmic {\it
formal} calculus is that when we use complex analysis, it is necessary
for the powers of the logarithms to be bounded from above, essentially
because a complex variable $z$ can also be expressed as the sum of the
series $z=\sum_{n\in \N}\frac{(\log z)^{n}}{n!}$.  Fifth, we have had
to combine our results on formal calculus, on logarithmic formal
calculus, and on complex analysis for series with both arbitrary real
powers and also logarithms to prove our main results on the
construction of the tensor category structures.  In many proofs, we
encounter expressions involving both formal variables and complex
variables, and thus we have had to develop new and delicate methods
exploiting both the formal and complex analysis methods that we have
just mentioned.  The proofs, which are not short (and cannot be),
accomplish the necessary interchanges of order of summations.}
\end{rema}

\begin{rema}\label{hist-btc}{\rm
The operator product expansion and resulting braided tensor category
structure constructed by the theory in \cite{tensor1}, \cite{tensor2},
\cite{tensor3}, \cite{tensor4} were originally structures whose
existence was conjectured: It was in their important study of
conformal field theory that Moore and Seiberg \cite{MS1} \cite{MS}
first discovered a set of polynomial equations from a suitable axiom
system for a ``rational conformal field theory.'' Inspired by a
comment of Witten, they observed an analogy between the theory of
these polynomial equations and the theory of tensor categories. The
structures given by these Moore-Seiberg equations were called
``modular tensor categories'' by I. Frenkel.  However, in the work of
Moore and Seiberg, as they commented, neither tensor product structure
nor other related structures were either formulated or constructed
mathematically. Later, Turaev formulated a precise notion of modular
tensor category in \cite{T1} and \cite{T} and gave examples of such
tensor categories from representations of quantum groups at roots of
unity, based on results obtained by many people on quantum groups and
their representations, especially those in the pioneering work
\cite{RT1} and \cite{RT2} by Reshetikhin and Turaev on the
construction of knot and $3$-manifold invariants from representations
of quantum groups. On the other hand, on the ``rational conformal
field theory'' side, a modular tensor category structure in this sense
on certain module categories for affine Lie algebras, and much more
generally, on certain module categories for ``chiral algebras''
associated with rational conformal field theories, was then believed
to exist by both physicists and mathematicians, but such structure was
not in fact constructed at that time.  Moore and Seiberg observed the
analogy mentioned above based on the assumption of the existence of a
suitable tensor product functor (including a tensor product module)
and derived their polynomial equations based on the assumption of the
existence of a suitable operator product expansion for chiral vertex
operators, which is essentially equivalent to assuming the
associativity of intertwining maps, as we have expressed it above.  As
we have discussed, the desired tensor product modules and functors
were constructed under suitable conditions in the series of papers
\cite{tensor1}, \cite{tensor2} and \cite{tensor3}, and in
\cite{tensor4} the appropriate natural associativity isomorphisms
among tensor products of triples of modules were constructed, and it
was shown that this is equivalent to the desired associativity of
intertwining maps (and thus the existence of a suitable operator
product expansion).  In particular, this work \cite{tensor1},
\cite{tensor2}, \cite{tensor3} and \cite{tensor4} served to construct
the desired braided tensor category structure in the generality of
suitable vertex operator algebras, including those associated with
affine Lie algebras and the Virasoro algebra as a very special case;
see \cite{HLaffine} and \cite{H3}, respectively.  (For a discussion of
the remaining parts of the modular tensor category structure in this
generality, see below and \cite{Hfinitetensor}.)  The results in these
papers will be generalized in this work.  In the special case of
affine Lie algebras and also in the special case of Virasoro-algebraic
structures, using the work of Tsuchiya-Ueno-Yamada \cite{TUY} and
Beilinson-Feigin-Mazur \cite{BFM} combined with a formulation of
braided tensor category structure by Deligne \cite{De}, one can obtain
the braided tensor category structure discussed above (but not the
modular tensor category structure).  }
\end{rema}

\subsection{Some recent applications and related literature}\label{literature}

We begin with a discussion concerning the ``rational'' case, with
semisimple module categories.  We also refer the reader to the recent
review by Fuchs, Runkel and Schweigert \cite{FRS} on rational
conformal field theory, which also in fact briefly discusses
nonrational conformal field theories, including in particular
logarithmic conformal field theories.

After the important work \cite{MS1} and \cite{MS} of Moore and
Seiberg, it was widely believed that the category of modules for a
suitable vertex operator algebra must have a structure of braided
tensor category satisfying additional properties related to the
modular invariance of the vertex operator algebra. As is mentioned in
Remark \ref{hist-btc}, for a suitable vertex operator algebra, the
work \cite{tensor1}, \cite{tensor2}, \cite{tensor3} and \cite{tensor4}
constructed a structure of braided tensor category on the category of
modules for the vertex operator algebra; see also \cite{H3} and
\cite{HLaffine}.  On the other hand, the precise and conceptual
formulation of the notion of modular tensor category by Turaev
\cite{T1} led to a mathematical conjecture that the category of
modules for a suitable vertex operator algebra can be endowed in a
natural way with modular tensor category structure in this sense. It
was in 2005 that this conjecture was finally proved by the first
author in \cite{rigidity} (see also the announcement \cite{HPNAS} and
the exposition \cite{Hconference}).  The hardest part of the proof of
this conjecture was the proof of the rigidity property of the braided
tensor category constructed in \cite{tensor1}, \cite{tensor2},
\cite{tensor3} and \cite{tensor4}.

Even in the case of a vertex operator algebra associated to an affine
Lie algebra or the Virasoro algebra, there was no proof of rigidity
for the braided tensor category of modules in the literature, before
the proof discovered in \cite{rigidity}. The works of
Tsuchiya-Ueno-Yamada \cite{TUY} and Beilinson-Feigin-Mazur \cite{BFM}
can be used to construct a structure of braided tensor category on the
category of modules for such a vertex operator algebra, but neither
the rigidity property nor the other main axiom for modular tensor
category structure, called the
nondegeneracy property, of these braided tensor categories has ever
been proved using the results or methods in those works. Under the
assumption that the braided tensor category structure on the category
of integrable highest weight (standard) modules of a fixed positive
integral level for an affine Lie algebra was already known to have the
rigidity property, Finkelberg \cite{F1} \cite{F2} showed that this
braided tensor category structure could be recovered by transporting
to this category the corresponding rigid braided tensor category
structure previously constructed for negative levels by Kazhdan and
Lusztig \cite{KL1}--\cite{KL5}. But since the rigidity was an
assumption needed in the proof, the work \cite{F1}, \cite{F2} did not
actually serve to give a construction of the braided tensor category
structure at positive integral level.  The book \cite{BK} asserted
that one had a construction of the structure of modular tensor
category on the category of modules for a vertex operator algebra
associated to an affine Lie algebra at positive integral level, and
while a construction of the structure of braided tensor category was
indeed given, there was no proof of the rigidity property, so that
even in the cases of affine Lie algebras and the Virasoro algebra, the
construction of the corresponding modular tensor category structures
was still an unsolved open problem before 2005.

Under the assumption of the rigidity for positive integral level, the
work \cite{F1}, \cite{F2} of Finkelberg combined with the work
\cite{KL1}--\cite{KL5} of Kazhdan and Lusztig established the
important equivalence between the braided tensor category of a
semisimple subquotient of the category of modules for a quantum group
at a root of unity and the braided tensor category of integrable
highest weight modules of a positive integral level for an affine Lie
algebra. The proof of the rigidity of the braided tensor category of
integrable highest weight modules of a positive integral level for an
affine Lie algebra, as a special case in \cite{rigidity}, based on the
the braided tensor category structure constructed in \cite{HLaffine},
as a special case in \cite{tensor1}, \cite{tensor2}, \cite{tensor3}
and \cite{tensor4}, thus in fact provided the completion of the proof
of the equivalence theorem that was the goal in \cite{F1} and
\cite{F2} above. As we have mentioned, the only known proof of this
rigidity requires the work \cite{tensor1}, \cite{tensor2},
\cite{tensor3} and \cite{tensor4}, and in particular, in the affine
Lie algebra case, uses the work \cite{HLaffine}.

The proof of the rigidity in \cite{rigidity} is highly nontrivial.
The reason why the rigidity was so hard is that one needed to prove
the Verlinde conjecture for suitable vertex operator algebras in order
to prove the rigidity, and the Verlinde conjecture requires the
consideration of genus-one as opposed to genus-zero conformal field
theory. The nondegeneracy property of the modular tensor category also
follows from the truth of the Verlinde conjecture.  The Verlinde
conjecture was discovered by E. Verlinde \cite{V} in 1987, and as was
demonstrated by Moore and Seiberg \cite{MS1} \cite{MS} in 1988, the
validity of the conjecture follows from their axiom system for a
rational conformal field theory.  However, the construction of
rational conformal field theories is much harder than the construction
of modular tensor categories, and this in turn requires the proof of
the Verlinde conjecture without the assumption of the axioms for a
rational conformal field theory.  The Verlinde conjecture for suitable
vertex operator algebras was proved in 2004 by the first author in
\cite{HVerlindeconjecture} (without the assumption of the axioms for a
rational conformal field theory), and its proof in turn depended on
the aspects of the theory of intertwining operators (the genus-zero
theory) developed in \cite{diff-eqn} and on the aspects of the theory
of $q$-traces of products or iterates of intertwining operators and
their modular invariance (the genus-one theory) developed in
\cite{Hmodular}.  (These works in turn depended on \cite{tensor1},
\cite{tensor2}, \cite{tensor3} and \cite{tensor4}.)  The modular
invariance theorem proved in the pioneering work \cite{Zhu1},
\cite{Zhu} of Zhu actually turned out to be only a very special case
of the stronger necessary result proved in \cite{Hmodular}, and was
far from enough for the purpose of establishing either the required
rigidity property or the required nondegeneracy property of the
modular tensor category structure.  The paper \cite{Hmodular}
established the most general modular invariance result in the
semisimple case and also constructed all genus-one correlation
functions of the corresponding chiral rational conformal field
theories. After Zhu's modular invariance was proved in 1990, the
modular invariance for products or iterates of more than one
intertwining operator was an open problem for a long time.  In the
case of products or iterates of at most one intertwining operator and
any number of vertex operators for modules, a straightforward
generalization of Zhu's result using his same method gives the modular
invariance (see \cite{M1}).  But for products or iterates of more than
one intertwining operator, Zhu's method is not sufficient because the
commutator formula that he used to derive his recurrence formula in
his proof has no generalization for intertwining operators. This was
one of the main reasons that for about 15 years after 1990, there had
been not much progress toward the proof of the rigidity and
nondegeneracy properties. In \cite{Hmodular}, this difficulty was
overcome by means of a proof that $q$-traces of products or iterates
of intertwining operators satisfy modular invariant differential
equations with regular singular points; the need for a recurrence
formula was thus bypassed.

We have been discussing the case of rational conformal field theories.
The present work includes as a special case a complete treatment of
the work \cite{tensor1}, \cite{tensor2}, \cite{tensor3} and
\cite{tensor4}, with much stronger results added as well; this work is
required for the results that we have just discussed.  The main theme
of the present work being the logarithmic generalization of this
theory, allowing categories of modules that are not completely
reducible, we would now like to comment on some recent applications
and related literature in the (much greater) logarithmic generality,
and also, in this generality we are in addition able to replace vertex
operator algebras by much more general vertex algebras equipped with a
suitable additional grading by an abelian group.  (Allowing
logarithmic structures and allowing vertex algebras with a grading by
an abelian group are ``unrelated'' generalizations of the context of
\cite{tensor1}, \cite{tensor2}, \cite{tensor3} and \cite{tensor4}; in
the present work we are able to carry out both generalizations
simultaneously.)

The triplet $\mathcal{W}$-algebras $\mathcal{W}(1,p)$, mentioned above,
are a class of vertex operator algebras
of central charge $1-6\frac{(p-1)^{2}}{p}$ which in recent years have
attracted a lot of attention from physicists and mathematicians.
Introduced by Kausch \cite{K1}, they 
have been studied extensively by Flohr \cite{Fl1} \cite{Fl2},
Gaberdiel-Kausch \cite{GK1} \cite{GK3}, Kausch \cite{K2},
Fuchs-Hwang-Semikhatov-Tipunin \cite{FHST}, Abe \cite{A},
Feigin-Ga{\u\i}nutdinov-Semikhatov-Tipunin \cite{FGST1}
\cite{FGST3}, Carqueville-Flohr \cite{CF}, Flohr-Gaberdiel \cite{FG},
Fuchs \cite{Fu}, Adamovi\'{c}-Milas \cite{AM2} \cite{AM5}
\cite{AM7}, Flohr-Grabow-Koehn \cite{FGK}, Flohr-Knuth \cite{FK},
Gaberdiel-Runkel \cite{GR1} \cite{GR2}, Ga{\u\i}nutdinov-Tipunin
\cite{GT}, Pearce-Rasmussen-Ruelle \cite{PRR1} \cite{PRR2},
Nagatomo-Tsuchiya \cite{NT2} and Rasmussen \cite{Ra4}.  A triplet
$\mathcal{W}$-algebra $V=\coprod_{n\in \Z}V_{(n)}$ satisfies the
positive energy condition ($V_{(0)}=\C\mathbf{1}$ and $V_{(n)}=0$ for
$n<0$) and the $C_{2}$-cofiniteness condition (the quotient space
$V/C_{2}(V)$ is finite dimensional, where $C_{2}(V)$ is the subspace
of $V$ spanned by the elements of the form $u_{-2}v$ for $u, v\in V$).
The $C_{2}$-cofiniteness condition was proved by Abe \cite{A} in the
simplest $p=2$ case and by Carqueville-Flohr \cite{CF} and
Adamovi\'{c}-Milas \cite{AM2} in the general case.

In \cite{H13}, the first author proved that for a vertex operator
algebra $V$ satisfying the positive energy condition and the
$C_{2}$-cofiniteness condition, the category of grading-restricted
generalized $V$-modules satisfies the assumptions needed to invoke the
theory carried out in the present work. The present work, combined
with \cite{H13} (for proving the assumptions of the present work),
thus establishes the logarithmic operator product expansion and
constructs the logarithmic tensor category theory for any vertex
operator algebra satisfying the positive energy condition and the
$C_{2}$-cofiniteness condition.  For example, the logarithmic tensor
products used heavily in the papers \cite{M3}, \cite{M4} and \cite{M5}
of Miyamoto are in fact constructed in the present work together with
\cite{H13}.  In particular, for a triplet $\mathcal{W}$-algebra $V$
discussed above, the category of grading-restricted generalized
$V$-modules indeed has the natural braided tensor category structure
constructed in the present work.  Many of the assertions involving a
logarithmic operator product expansion and a logarithmic tensor
category theory in the works on triplet $\mathcal{W}$-algebras
mentioned above are mathematically formulated and established in the
present work together with the paper \cite{H13}, so that now, they do
not have to be taken as unproved assumptions in those works.

Based on the results of Feigin-Ga{\u\i}nutdinov-Semikhatov-Tipunin
\cite{FGST3} and of Fuchs-Hwang-Semikhatov-Tipunin \cite{FHST},
Feigin, Ga{\u\i}nutdinov, Semikhatov and Tipunin conjectured
\cite{FGST1} an equivalence between the braided finite tensor category
of grading-restricted generalized modules for a triplet
$\mathcal{W}$-algebra and the braided finite tensor category of
suitable modules for a restricted quantum group.  Their formulation of
the conjecture also includes the statement that the categories of
grading-restricted generalized modules for the triplet
$\mathcal{W}$-algebras considered in their paper are indeed braided
tensor categories.  Assuming the existence of the braided tensor
category structure on the triplet $\mathcal{W}$-algebra with $p=2$,
Feigin, Ga{\u\i}nutdinov, Semikhatov and Tipunin gave a proof of their
conjecture. However, in the case $p\ne 2$, Kondo and Saito \cite{KS}
showed that the tensor category of modules for the corresponding
restricted quantum group is not braided. Thus, the conjecture in the
case $p\ne 2$ cannot be true as it is stated,
although the equivalence between the abelian categories was
proved in \cite{NT2} for all $p$.  It is believed that the
correct formulation of the conjecture and the proof will be possible
only after the conformal-field-theoretic aspects of the
representations of triplet $\mathcal{W}$-algebras are studied
thoroughly.  As we mentioned above, the present work, the paper
\cite{H13} and the papers \cite{A}, \cite{CF} and \cite{AM2} provide a
proof of the assumption in their conjecture that the categories of
grading-restricted generalized modules for the triplet
$\mathcal{W}$-algebras are indeed braided tensor categories. We expect
that further studies of the tensor-categorical structures and
conformal-field-theoretic properties for triplet
$\mathcal{W}$-algebras will provide a correct formulation and proof of
suitable equivalence between categories of suitable modules for
triplet $\mathcal{W}$-algebras and for restricted quantum groups.

In \cite{H14}, the first author introduced a notion of generalized
twisted module associated to a general automorphism of a vertex
operator algebra, including an automorphism of infinite order.  The
first author in \cite{H14} also gave a construction of such
generalized twisted modules associated to the automorphisms obtained
by exponentiating weight $1$ elements of the vertex operator algebra.
If the automorphism of the vertex operator algebra does not act
semisimply, the twisted vertex operators for these generalized twisted
modules must involve the logarithm of a formal or complex variable,
and we need additional $\C/\Z$- or $\C$-gradings on these generalized
twisted modules.  As was noticed by Milas, the triplet
$\mathcal{W}$-algebras are fixed-point subalgebras of suitable vertex
operator algebras constructed from a one-dimensional lattice under an
automorphism obtained by exponentiating a weight $1$ element. In
particular, some logarithmic intertwining operators constructed in
\cite{AM4} are in fact twisted vertex operators. Thus the paper
\cite{H14} provided an orbifold approach to the representation theory
of triplet $\mathcal{W}$-algebras.  (This orbifold point of view is
one of the analogues of the orbifold point of view for vertex operator
algebras introduced in \cite{FLM2}.)  Since the automorphisms involved
indeed do not act on the vertex operator algebra semisimply, the
twisted vertex operators for the generalized twisted modules
associated to these automorphisms must involve the logarithm of the
variables, and we also need additional $\C/\Z$- or $\C$-gradings on
these generalized twisted modules.  Here $\C/\Z$ or $\C$ are instances
of the additional grading abelian group in the present work.  Thus we
need the general framework and results in the present work, including
both the logarithmic generality and also the additional abelian-group
gradings, for the study of these generalized twisted modules.

Many of the results on the representation theory of triplet
$\mathcal{W}$-algebras have also been generalized to the more general
case of $\mathcal{W}(p, q)$-algebras \cite{FGST2} of central charge
$1-6\frac{(p-q)^{2}}{pq}$, $q>p>0$ coprime (see for example
\cite{AM1}, \cite{S}, \cite{RP2}, \cite{Ra1}, \cite{Ra2},
\cite{Ra3}, \cite{Ra5}, \cite{GRW}, \cite{AM6} and \cite{Wo}), and to
$N=1$ triplet vertex operator superalgebras (see \cite{AM3},
\cite{AM5} and \cite{AM7}).  
Results have also been obtained for the vertex operator subalgebras of
the algebras $\mathcal{W}(p,q)$ generated by the Virasoro algebra (see
\cite{GK2}, \cite{EF}, \cite{PRZ}, \cite{ReS}, \cite{MR}, \cite{RP},
\cite{BFGT} and \cite{GV}).  The $C_{2}$-cofiniteness of the
$\mathcal{W}(2, q)$-algebras has been proved by Adamovi\'{c} and Milas
in \cite{AM6}. Thus using the results obtained in \cite{H14}, the
theory developed in the present work applies to these $\mathcal{W}(2,
q)$-algebras, yielding braided tensor categories.  The $N=1$ triplet
vertex operator superalgebras introduced by Adamovi\'{c} and Milas in
\cite{AM3} are also proved by these authors in \cite{AM5} to be
$C_{2}$-cofinite. As was mentioned above in Section 1.2, the theory
developed in this work also applies to vertex superalgebras.  The same
remarks apply to the results in \cite{H14}.  Thus the theory developed
in the present work applies to these $N=1$ triplet vertex operator
superalgebras, producing the corresponding braided tensor categories.

Finally, we would like to emphasize that it is interesting that the
methods developed and used in the present work, even in the special
case of categories of modules for an affine Lie algebra at negative
levels, are very different from those developed and used by Kazhdan
and Lusztig in \cite{KL1}--\cite{KL5}, and are much more general.  The
methods used in \cite{KL1}--\cite{KL5}, closely related to algebraic
geometry, depended heavily on the Knizhnik-Zamolodchikov equations.
In the present work, we use and develop the general theory of vertex
(operator) algebras (and generalizations), requiring both formal
calculus theory and complex analysis, and we do not use algebraic
geometry. Also, in the present work and in the work \cite{Z1} and
\cite{Z2}, which verified the assumptions needed for the application
of the present theory, although we need to show that products of
intertwining operators satisfy certain differential equations with
regular singular points, no explicit form of the equations, such as
the explicit form of the Knizhnik-Zamolodchikov equations, is
needed. In fact, because for a general vertex (operator) algebra
satisfying those assumptions in the present work or in \cite{H13} no
explicit form of the differential equations such as the form of the
Knizhnik-Zamolodchikov equations exists, it was crucial that in the
work \cite{tensor1}, \cite{tensor2}, \cite{tensor3}, \cite{tensor4},
\cite{diff-eqn}, the present work and \cite{H13}, we have developed
methods that are independent of the explicit form of the differential
equations.  Another interesting difference between the present general
theory and this work of Kazhdan and Lusztig is that logarithmic
structures (necessarily) pervade our theory, starting from the
vertex-algebraic foundations, while the logarithmic nature of
solutions of the Knizhnik-Zamolodchikov equations involved in
\cite{KL1}--\cite{KL5} did not have to be emphasized there.

\subsection{Main results of the present work}

In this section, we state the main results of the present work,
numbered as in the main text.  The reader is referred to the relevant
sections for definitions, notations and details.

Let $A$ be an abelian group and $\tilde{A}$ an abelian group
containing $A$ as a subgroup. Let $V$ be a strongly $A$-graded
M\"{o}bius or conformal vertex algebra, as defined in Section 2.  Let
$\mathcal{C}$ be a full subcategory of the category $\mathcal{M}_{sg}$
of strongly $\tilde{A}$-graded (ordinary) $V$-modules or the category
$\mathcal{G}\mathcal{M}_{sg}$ of strongly $\tilde{A}$-graded
generalized $V$-modules, closed under the contragredient functor and
under taking finite direct sums; see Section 2 and Assumptions
\ref{assum} and \ref{assum-c}.

In Section 4, the notions of $P(z)$- and $Q(z)$-tensor product functor
are defined in terms of $P(z)$- and $Q(z)$-intertwining maps and
$P(z)$- and $Q(z)$-products; intertwining maps are related to
logarithmic intertwining operators, defined and studied in Section
3. The symbols $P(z)$ and $Q(z)$ refer to the moduli space elements
described in Remarks \ref{P(z)geometry} and \ref{Q(z)geometry},
respectively.  In Section 5, we give a construction of the
$P(z)$-tensor product of two objects of $\mathcal{C}$, when this
structure exists.  For $W_{1}, W_{2}\in \ob \mathcal{C}$, define the
subset
\[
W_{1}\hboxtr_{P(z)}W_{2}\subset (W_{1}\otimes W_{2})^{*}
\]
of $(W_{1}\otimes W_{2})^{*}$ to be the union, or equivalently, the
sum, of the images
\[
I'(W')\subset (W_{1}\otimes W_{2})^{*}
\]
as $(W; I)$ ranges through all the $P(z)$-products of $W_{1}$ and
$W_{2}$ with $W\in \ob \mathcal{C}$, where $I'$ is a map corresponding
naturally to the $P(z)$-intertwining map $I$ and where $W'$ is the
contragredient (generalized) module of $W$.

The following two results give the construction of the $P(z)$-tensor
product:

\setcounter{section}{5}
\setcounter{rema}{36}
\begin{propo}
Let $W_{1}, W_{2}\in \ob \mathcal{C}$.  If $(W_1\hboxtr_{P(z)} W_2,
Y'_{P(z)})$ is an object of ${\cal C}$ (where $Y'_{P(z)}$ is the
natural action of $V$), denote by $(W_1\boxtimes_{P(z)} W_2,
Y_{P(z)})$ its contragredient (generalized) module:
\[
W_1\boxtimes_{P(z)} W_2 = (W_1\hboxtr_{P(z)} W_2)'.
\]
Then the $P(z)$-tensor product of $W_{1}$ and $W_{2}$ in ${\cal C}$
exists and is
\[
(W_1\boxtimes_{P(z)} W_2, Y_{P(z)}; i'),
\]
where $i$ is the natural inclusion {}from $W_1\hboxtr_{P(z)} W_2$ to
$(W_1\otimes W_2)^*$.  Conversely, let us assume that $\mathcal{C}$ is
closed under images.  If the $P(z)$-tensor product of $W_1$ and $W_2$
in ${\cal C}$ exists, then $(W_1\hboxtr_{P(z)} W_2, Y'_{P(z)})$ is an
object of ${\cal C}$.
\end{propo}

For 
\[
\lambda
\in (W_{1}\otimes W_{2})^{*},
\]
let $W_{\lambda}$ be the smallest doubly graded subspace of
$((W_1\otimes W_2)^*)_{[ {\mathbb C} ]}^{( \tilde A )}$ (the direct
sum of the homogeneous subspaces with respect to the gradings both by
conformal generalized weights and by $\tilde{A}$) containing $\lambda$
and stable under the component operators of the operators
$Y'_{P(z)}(v,x)$ for $v\in V$, $m\in {\mathbb Z}$, and under the
operators $L'_{P(z)}(-1)$, $L'_{P(z)}(0)$ and $L'_{P(z)}(1)$ (to
handle the M\"{o}bius but non-conformal case).  Let
\[
\comp_{P(z)}((W_1\otimes W_2)^*),
\]
\[
\lgr_{[\C]; P(z)}((W_1\otimes W_2)^*)
\]
and 
\[
\lgr_{(\C); P(z)}((W_1\otimes W_2)^*)
\]
be the spaces of elements of $(W_1\otimes W_2)^*$ satisfying the
$P(z)$-compatibility condition, the $P(z)$-local grading restriction
condition and the $L(0)$-semisimple $P(z)$-local grading restriction
condition, respectively, as defined in Section 4; the subscript $(\C)$
refers to the semisimplicity of the action of $L(0)$ in this case, so
that generalized weights are weights.

\setcounter{rema}{49}
\begin{theo}
Suppose that for every element 
\[
\lambda
\in \comp_{P(z)}((W_1\otimes W_2)^*)\cap \lgr_{[\C]; P(z)}((W_1\otimes W_2)^*)
\]
the space $W_{\lambda}$ (which is a (strongly-graded) generalized
module) is a generalized submodule of some object of ${\cal C}$
included in $(W_1\otimes W_2)^*$ (this holds vacuously if ${\cal C} =
{\cal GM}_{sg}$).  Then
\[
W_1\hboxtr_{P(z)}W_2=\comp_{P(z)}((W_1\otimes W_2)^*)
\cap \lgr_{[\C]; P(z)}((W_1\otimes W_2)^*).
\]
Suppose that ${\cal C}$ is a category of strongly-graded $V$-modules
(that is, ${\cal C}\subset {\cal M}_{sg}$) and that for every element
\[
\lambda
\in \comp_{P(z)}((W_1\otimes W_2)^*)\cap \lgr_{(\C); P(z)}((W_1\otimes W_2)^*)
\]
the space $W_{\lambda}$ (which is a (strongly-graded) $V$-module) is a
submodule of some object of ${\cal C}$ included in $(W_1\otimes
W_2)^*$ (which holds vacuously if ${\cal C} = {\cal M}_{sg}$).  Then
\[
W_1\hboxtr_{P(z)}W_2=\comp_{P(z)}((W_1\otimes W_2)^*)
\cap \lgr_{(\C); P(z)}((W_1\otimes W_2)^*).
\]
\end{theo}

The hard parts of the proof of Theorem
\ref{characterizationofbackslash} are given in Section 6.

We also give an analogous construction of $Q(z)$-tensor products in
these sections.

For the construction of the natural associativity isomorphism between
suitable pairs of triple tensor product functors, we assume that for
any object of ${\cal C}$, all the (generalized) weights are real
numbers and in addition there exists $K\in \Z_{+}$ such that
\[
(L(0)-L(0)_{s})^{K}=0
\]
on the module, $L(0)_{s}$ being the semisimple part of $L(0)$ (the
latter condition holding vacuously when $\mathcal{C}$ is in ${\cal
M}_{sg}$); see Assumption \ref{assum-exp-set}.

The main hard parts of the construction of the associativity
isomorphisms are presented in Section 9, after necessary preparation
in Section 8.  To discuss these results, we need the important
$P^{(1)}(z)$- and $P^{(2)}(z)$-local grading restriction conditions on
\[
\lambda\in (W_{1}\otimes W_{2}\otimes W_{3})^{*}
\]
(where $W_{1}$, $W_{2}$, and $W_{3}$ are objects of $\mathcal{C}$)
and their $L(0)$-semisimple versions. 
Here we state the (two-part) $P^{(2)}(z)$-local grading 
restriction condition, the other conditions being analogous: 

\begin{description}

\item{\bf The $P^{(2)}(z)$-local grading restriction condition}

(a) The {\em $P^{(2)}(z)$-grading condition}: For any $w_{(3)}\in
W_3$,  there exists  a formal series $\sum_{n\in \R}\lambda^{(2)}_{n}$ 
with
\[
\lambda^{(2)}_{n}\in \coprod_{\beta\in \tilde{A}}
((W_{1}\otimes W_{2})^{*})_{[n]}^{(\beta)}
\]
for $n\in \R$, an open neighborhood of $z'=0$, and $N \in \N$
such that for $w_{(1)}\in W_{1}$ and $w_{(2)}\in W_{2}$,
the series 
\[
\sum_{n\in \R}(e^{z'L'_{P(z)}(0)}\lambda^{(2)}_{n})(w_{(1)}\otimes w_{(2)})
\]
has the following properties:
\begin{enumerate}

\item[(i)] It can be written as the iterated series
\[
\sum_{n\in \R}e^{nz'}\left(\left(\sum_{i=0}^{N}\frac{(z')^{i}}{i!}
(L'_{P(z)}(0)-n)^{i}\lambda^{(2)}_{n}
\right)(w_{(1)}\otimes 
w_{(2)})\right).
\]

\item[(ii)] It is absolutely convergent for $z'\in \C$ in the neighborhood of
$z'=0$ above.

\item[(iii)] It is absolutely convergent to $\mu^{(2)}_{\lambda,
w_{(3)}}(w_{(1)}\otimes w_{(2)})$ when $z'=0$:
\[
\sum_{n\in \R}\lambda^{(2)}_{n}(w_{(1)}\otimes w_{(2)})
=\mu^{(2)}_{\lambda, w_{(3)}}(w_{(1)}\otimes w_{(2)})
=\lambda(w_{(1)}\otimes w_{(2)}\otimes w_{(3)})
\]
(the last equality being the definition of $\mu^{(2)}_{\lambda, w_{(3)}}$).

\end{enumerate}

(b) For any $w_{(3)}\in W_3$, let $W^{(2)}_{\lambda, w_{(3)}}$ be the
smallest doubly graded subspace
of $((W_1\otimes W_2)^{*})_{[\R]}^{(\tilde{A})}$ containing all the
terms $\lambda^{(2)}_{n}$ in the formal series in (a) and stable under
the component operators  of the operators
$Y'_{P(z)}(v, x)$ for $v\in V$, $m\in {\mathbb Z}$, and under the
operators $L'_{P(z)}(-1)$, $L'_{P(z)}(0)$ and $L'_{P(z)}(1)$.   Then
$W^{(2)}_{\lambda, w_{(3)}}$ has the properties
\begin{eqnarray*}
&\dim(W^{(2)}_{\lambda, w_{(3)}})^{(\beta)}_{[n]}<\infty,&\\
&(W^{(2)}_{\lambda, w_{(3)}})^{(\beta)}_{[n+k]}=0\;\;\mbox{ for
}\;k\in {\mathbb Z} \;\mbox{ sufficiently negative}&
\end{eqnarray*}
for any $n\in \R$ and $\beta\in \tilde A$, where the
subscripts denote the $\R$-grading by
$L'_{P(z)}(0)$-(generalized) eigenvalues and the superscripts denote
the $\tilde A$-grading.

\end{description}

The following result gives, among other things, the deep fact that
when $\lambda$ is obtained from a suitable product of intertwining
maps, the elements $\lambda^{(2)}_{n}$ for $n\in \R$ in the assumed
$P^{(2)}(z)$-local grading restriction condition for suitable $z\in
\C^{\times}$ satisfy the $P(z)$-compatibility condition:

\setcounter{section}{9}
\setcounter{rema}{16}
\begin{theo}
Assume that the convergence condition for intertwining maps in
$\mathcal{C}$ (see Section 7) holds and that
\[
|z_1|>|z_2|>|z_{1}-z_{2}|>0.
\]
Let $W_{1}$, $W_{2}$, $W_{3}$, $W_{4}$, $M_{1}$ and $M_{2}$ be objects
of $\mathcal{C}$ and let $I_{1}$, $I_{2}$, $I^1$ and $I^2$ be
$P(z_1)$-, $P(z_2)$-~, $P(z_2)$- and $P(z_{1}-z_{2})$-intertwining
maps of types ${W_4}\choose {W_1M_1}$, ${M_1}\choose {W_2W_3}$,
${W_4}\choose {M_2W_3}$ and ${M_2}\choose {W_1W_2}$, respectively.
Let $w'_{(4)}\in W'_4$.
\begin{enumerate}

\item
Suppose that $(I_1\circ (1_{W_1}\otimes I_2))'(w'_{(4)})$ satisfies
Part (a) of the $P^{(2)}(z_{1}-z_{2})$-local grading restriction condition,
that is, the $P^{(2)}(z_{1}-z_{2})$-grading condition (or the
$L(0)$-semisimple $P^{(2)}(z_{1}-z_{2})$-grading condition when $\mathcal{C}$
is in $\mathcal{M}_{sg}$).  For any $w_{(3)}\in W_{3}$, let
$\sum_{n\in \R}\lambda_{n}^{(2)}$ be a series weakly absolutely
convergent to
\[
\mu^{(2)}_{(I_1\circ (1_{W_1}\otimes I_2))'(w'_{(4)}), w_{(3)}} \in
(W_1 \otimes W_2)^*
\]
as indicated in the $P^{(2)}(z_{1}-z_{2})$-grading condition (or the
$L(0)$-semisimple $P^{(2)}(z_{1}-z_{2})$-grading condition), and suppose in
addition that the elements $\lambda_{n}^{(2)} \in (W_{1}\otimes
W_{2})^{*}$, $n\in \R$, satisfy the $P(z_{1}-z_{2})$-lower truncation
condition (Part (a) of the $P(z_{1}-z_{2})$-compatibility condition in
Section 5).  Then each $\lambda_{n}^{(2)}$ satisfies the (full)
$P(z_{1}-z_{2})$-compatibility condition.  Moreover, the corresponding space
\[
W^{(2)}_{(I_1\circ (1_{W_1}\otimes I_2))'(w'_{(4)}), w_{(3)}} \subset
(W_{1}\otimes W_{2})^{*},
\]
equipped with the vertex operator map given by $Y'_{P(z_{1}-z_{2})}$ and the
operators $L'_{P(z_{1}-z_{2})}(j)$ for $j=-1, 0, 1$, is a doubly-graded
generalized $V$-module, and when $\mathcal{C}$ is in
$\mathcal{M}_{sg}$, a doubly-graded $V$-module.  In particular, if
$(I_1\circ (1_{W_1}\otimes I_2))'(w'_{(4)})$ satisfies the full
$P^{(2)}(z_{1}-z_{2})$-local grading restriction condition (or the
$L(0)$-semisimple $P^{(2)}(z_{1}-z_{2})$-local grading restriction condition
when $\mathcal{C}$ is in $\mathcal{M}_{sg}$), then $W^{(2)}_{(I_1\circ
(1_{W_1}\otimes I_2))'(w'_{(4)}), w_{(3)}}$ is an object of
$\mathcal{GM}_{sg}$ (or $\mathcal{M}_{sg}$ when $\mathcal{C}$ is in
$\mathcal{M}_{sg}$); in this case, the assumption that each
$\lambda_{n}^{(2)}$ satisfies the $P(z_{1}-z_{2})$-lower truncation condition
is redundant.

\item
Analogously, suppose that $(I^1\circ (I^2\otimes 1_{W_3}))'(w'_{(4)})$
satisfies Part (a) of the $P^{(1)}(z_2)$-local grading restriction
condition, that is, the $P^{(1)}(z_2)$-grading condition (or the
$L(0)$-semisimple $P^{(1)}(z_2)$-grading condition when $\mathcal{C}$
is in $\mathcal{M}_{sg}$).  For any $w_{(1)}\in W_{1}$, let
$\sum_{n\in \R}\lambda_{n}^{(1)}$ be a series weakly absolutely
convergent to
\[
\mu^{(1)}_{(I^1\circ (I^2\otimes 1_{W_3}))'(w'_{(4)}), w_{(1)}} \in
(W_2 \otimes W_3)^*
\]
as indicated in the $P^{(1)}(z_2)$-grading condition (or the
$L(0)$-semisimple $P^{(1)}(z_2)$-grading condition), and suppose in
addition that the elements $\lambda_{n}^{(1)} \in (W_{2}\otimes
W_{3})^{*}$, $n\in \R$, satisfy the $P(z_{2})$-lower truncation
condition (Part (a) of the $P(z_{2})$-compatibility condition).  Then
each $\lambda_{n}^{(1)}$ satisfies the (full) $P(z_{2})$-compatibility
condition.  Moreover, the corresponding space
\[
W^{(1)}_{(I^1\circ (I^2\otimes 1_{W_3}))'(w'_{(4)}), w_{(1)}} \subset
(W_{2}\otimes W_{3})^{*},
\]
equipped with the vertex operator map given by $Y'_{P(z_{2})}$ and the
operators $L'_{P(z_{2})}(j)$ for $j=-1, 0, 1$, is a doubly-graded
generalized $V$-module, and when $\mathcal{C}$ is in
$\mathcal{M}_{sg}$, a doubly-graded $V$-module.  In particular, if
$(I^1\circ (I^2\otimes 1_{W_3}))'(w'_{(4)})$ satisfies the full
$P^{(1)}(z_2)$-local grading restriction condition (or the
$L(0)$-semisimple $P^{(1)}(z_2)$-local grading restriction condition
when $\mathcal{C}$ is in $\mathcal{M}_{sg}$), then $W^{(1)}_{(I^1\circ
(I^2\otimes 1_{W_3}))'(w'_{(4)}), w_{(1)}}$ is an object of
$\mathcal{GM}_{sg}$ (or $\mathcal{M}_{sg}$ when $\mathcal{C}$ is in
$\mathcal{M}_{sg}$); in this case, the assumption that each
$\lambda_{n}^{(1)}$ satisfies the $P(z_2)$-lower truncation condition
is redundant.
\end{enumerate}
\end{theo}

The following result, based heavily on the previous theorem,
establishes the associativity of intertwining maps:

\setcounter{rema}{22}
\begin{theo}
Assume that $\mathcal{C}$ is closed under images, that the convergence
condition for intertwining maps in $\mathcal{C}$ holds and that
\[
|z_1|>|z_2|>|z_{1}-z_{2}|>0.
\]
Let $W_{1}$, $W_{2}$, $W_{3}$, $W_{4}$, $M_{1}$ and $M_{2}$ be objects
of $\mathcal{C}$.  Assume also that $W_1\boxtimes_{P(z_{1}-z_{2})} W_2$ and
$W_2\boxtimes_{P(z_2)} W_3$ exist in $\mathcal{C}$.

\begin{enumerate}

\item  Let
$I_{1}$ and $I_{2}$ be $P(z_1)$- and $P(z_2)$-intertwining maps of 
types ${W_4}\choose {W_1M_1}$ and
${M_1}\choose {W_2W_3}$, respectively. 
Suppose that for each $w'_{(4)} \in W'_{4}$,
\[
\lambda=(I_1\circ (1_{W_1}\otimes I_2))'(w'_{(4)}) \in
(W_{1}\otimes W_{2} \otimes W_{3})^{*}
\]
satisfies the $P^{(2)}(z_{1}-z_{2})$-local grading restriction condition (or
the  $L(0)$-semisimple $P^{(2)}(z_{1}-z_{2})$-local grading restriction
condition when $\mathcal{C}$ is in $\mathcal{M}_{sg}$). For
$w'_{(4)}\in W'_{4}$ and $w_{(3)}\in W_{3}$, let $\sum_{n\in
\R}\lambda_{n}^{(2)}$ be the (unique) series weakly absolutely convergent to
$\mu^{(2)}_{\lambda, w_{(3)}}$ as indicated in the
$P^{(2)}(z_{1}-z_{2})$-grading condition (or the $L(0)$-semisimple
$P^{(2)}(z_{1}-z_{2})$-grading condition).  Suppose also that for each $n \in
\R$, $w'_{(4)} \in W'_4$ and $w_{(3)} \in W_3$, the generalized
$V$-submodule of $W^{(2)}_{\lambda, w_{(3)}}$ generated by
$\lambda_{n}^{(2)}$ is a generalized $V$-submodule of some object of
$\mathcal{C}$ included in $(W_1 \otimes W_2)^*$.  Then the product
\[
I_1\circ (1_{W_1}\otimes I_2)
\]
can be expressed as an
iterate, and in fact, there exists a unique
$P(z_2)$-intertwining map $I^{1}$ of type ${W_4\choose
W_1\boxtimes_{P(z_{1}-z_{2})} W_2\,\,W_3}$ such that
\[
\langle w'_{(4)},I_1(w_{(1)} \otimes I_2(w_{(2)} \otimes w_{(3)}))\rangle
=\langle w'_{(4)}, I^{1}((w_{(1)}\boxtimes_{P(z_{1}-z_{2})} w_{(2)})\otimes 
w_{(3)})\rangle
\]
for all $w_{(1)}\in W_1$, $w_{(2)}\in W_2$, $w_{(3)}\in W_3$ and
$w'_{(4)}\in W'_4$.

\item Analogously, let $I^1$ and $I^2$ be $P(z_2)$- and
$P(z_{1}-z_{2})$-intertwining maps of types ${W_4}\choose {M_2W_3}$ and
${M_2}\choose {W_1W_2}$, respectively. Suppose that for each
$w'_{(4)} \in W'_{4}$,
\[
\lambda=(I^1\circ (I^2 \otimes 1_{W_3}))'(w'_{(4)}) \in
(W_{1}\otimes W_{2} \otimes W_{3})^{*}
\]
satisfies the $P^{(1)}(z_2)$-local grading restriction condition (or
the $L(0)$-semisimple $P^{(1)}(z_2)$-local grading restriction
condition when $\mathcal{C}$ is in $\mathcal{M}_{sg}$). For
$w'_{(4)}\in W'_{4}$ and $w_{(1)}\in W_{1}$, let $\sum_{n\in
\R}\lambda_{n}^{(1)}$ be the (unique) series weakly absolutely
convergent to $\mu^{(1)}_{\lambda, w_{(1)}}$ as indicated in the
$P^{(1)}(z_2)$-grading condition (or the $L(0)$-semisimple
$P^{(1)}(z_2)$-grading condition).  Suppose also that for each
$n \in \R$, $w'_{(4)} \in W'_4$ and $w_{(1)} \in W_1$, the generalized
$V$-submodule of $W^{(1)}_{\lambda, w_{(1)}}$ generated by
$\lambda_{n}^{(1)}$ is a generalized $V$-submodule of some object of
$\mathcal{C}$ included in $(W_2 \otimes W_3)^*$.  Then 
the iterate
\[
I^1\circ (I^2 \otimes 1_{W_3})
\]
can be expressed as a product, and in fact, there exists a unique
$P(z_1)$-intertwining map $I_{1}$ of type
${W_4\choose W_1\,\,W_2\boxtimes_{P(z_2)} W_3}$ such that
\[
\langle w'_{(4)}, I^1(I^2(w_{(1)}\otimes w_{(2)})\otimes w_{(3)})\rangle
=\langle w'_{(4)}, I_{1}(w_{(1)}\otimes (w_{(2)}\boxtimes_{P(z_2)}w_{(3)}))
\rangle
\]
for all $w_{(1)}\in W_1$, $w_{(2)}\in W_2$, $w_{(3)}\in W_3$ and
$w'_{(4)}\in W'_4$. 

\end{enumerate}
\end{theo}

The hard part of the proof of this theorem is the proof of Lemma 
\ref{intertwine-tau}. This associativity of intertwining maps 
immediately gives the following important associativity of 
logarithmic intertwining operators, which is a strong version of 
logarithmic operator product expansion:

\begin{corol}
Assume that $\mathcal{C}$ is closed under images, that the convergence
condition for intertwining maps in $\mathcal{C}$ holds and that
\[
|z_1|>|z_2|>|z_{1}-z_{2}|>0.
\]
Let $W_{1}$, $W_{2}$, $W_{3}$, $W_{4}$, $M_{1}$ and $M_{2}$ be objects
of $\mathcal{C}$.  Assume also that $W_1\boxtimes_{P(z_{1}-z_{2})} W_2$ and
$W_2\boxtimes_{P(z_2)} W_3$ exist in $\mathcal{C}$.

\begin{enumerate}

\item Let $\Y_{1}$ and $\Y_{2}$ be logarithmic intertwining operators
(ordinary intertwining operators in the 
case that $\mathcal{C}$ is in $\mathcal{M}_{sg}$) 
of types ${W_4}\choose {W_1M_1}$ and ${M_1}\choose {W_2W_3}$,
respectively.  Suppose that for each $w'_{(4)} \in W'_{4}$, the
element $\lambda\in (W_{1}\otimes W_{2}\otimes W_{3})^{*}$ given by
\[
\lambda(w_{(1)}\otimes w_{(2)}\otimes w_{(3)})
=\langle w'_{(4)}, \Y_1(w_{(1)}, x_{1})\Y_2(w_{(2)}, x_{2})w_{(3)}\rangle
\lbar_{x_{1}=z_{1},\;x_{2}=z_{2}}
\]
for $w_{(1)}\in W_{1}$, $w_{(2)}\in
W_{2}$ and $w_{(3)}\in W_{3}$ satisfies the $P^{(2)}(z_{1}-z_{2})$-local
grading restriction condition (or the $L(0)$-semisimple
$P^{(2)}(z_{1}-z_{2})$-local grading restriction condition when $\mathcal{C}$
is in $\mathcal{M}_{sg}$). For $w'_{(4)}\in W'_{4}$ and $w_{(3)}\in
W_{3}$, let $\sum_{n\in \R}\lambda_{n}^{(2)}$ be the (unique) series
weakly absolutely convergent to $\mu^{(2)}_{\lambda, w_{(3)}}$ as
indicated in the $P^{(2)}(z_{1}-z_{2})$-grading condition (or the
$L(0)$-semisimple $P^{(2)}(z_{1}-z_{2})$-grading condition).  Suppose also
that for each $n \in \R$, $w'_{(4)} \in W'_4$ and $w_{(3)} \in W_3$,
the generalized $V$-submodule of $W^{(2)}_{\lambda, w_{(3)}}$
generated by $\lambda_{n}^{(2)}$ is a generalized $V$-submodule of
some object of $\mathcal{C}$ included in $(W_1 \otimes W_2)^*$.  Then
there exists a unique logarithmic intertwining operator (a unique
ordinary intertwining operator in the case that $\mathcal{C}$ is in
$\mathcal{M}_{sg}$) $\Y^{1}$ of type ${W_4\choose
W_1\boxtimes_{P(z_{1}-z_{2})} W_2\,\,W_3}$ such that
\begin{eqnarray*}
\lefteqn{\langle w'_{(4)},\Y_1(w_{(1)}, x_{1}) \Y_2(w_{(2)}, x_{2}) w_{(3)}\rangle
\lbar_{x_{1}=z_{1},\;x_{2}=z_{2}}}\nn
&&=\langle w'_{(4)}, \Y^{1}(\Y_{\boxtimes_{P(z_{1}-z_{2})}, 0}
(w_{(1)}, x_{0})w_{(2)}, x_{2})w_{(3)})\rangle
\lbar_{x_{0}=z_{1}-z_{2},\;x_{2}=z_{2}}
\end{eqnarray*}
for all $w_{(1)}\in
W_1$, $w_{(2)}\in W_2$, $w_{(3)}\in W_3$ and $w'_{(4)}\in W'_4$.  In
particular, the product of the logarithmic intertwining operators
(ordinary intertwining operators in the case that $\mathcal{C}$ is in
$\mathcal{M}_{sg}$) $\Y_{1}$ and $\Y_{2}$ evaluated at $z_{1}$ and
$z_{2}$, respectively, can be expressed as an iterate (with the
intermediate generalized $V$-module $W_1\boxtimes_{P(z_{1}-z_{2})} W_2$) of
logarithmic intertwining operators (ordinary intertwining operators in
the case that $\mathcal{C}$ is in $\mathcal{M}_{sg}$) evaluated at
$z_{2}$ and $z_{1}-z_{2}$.

\item Analogously, let $\Y^1$ and $\Y^2$ be logarithmic intertwining
operators (ordinary intertwining operators in the 
case that $\mathcal{C}$ is in $\mathcal{M}_{sg}$) 
of types ${W_4}\choose {M_2W_3}$ and ${M_2}\choose
{W_1W_2}$, respectively. Suppose that for each $w'_{(4)}\in W'_{4}$,
the element $\lambda\in (W_{1}\otimes W_{2}\otimes W_{3})^{*}$ given
by
\[
\lambda(w_{(1)}\otimes w_{(2)}\otimes w_{(3)})
=\langle w'_{(4)}, \Y^{1}(\Y^2(w_{(1)}, x_{0})w_{(2)}, x_{2})w_{(3)}\rangle
\lbar_{x_{0}=z_{1}-z_{2},\;x_{2}=z_{2}}
\]
satisfies the $P^{(1)}(z_2)$-local
grading restriction condition (or the $L(0)$-semisimple
$P^{(1)}(z_2)$-local grading restriction condition when $\mathcal{C}$
is in $\mathcal{M}_{sg}$). For $w'_{(4)}\in W'_{4}$ and $w_{(1)}\in
W_{1}$, let $\sum_{n\in \R}\lambda_{n}^{(1)}$ be the (unique) series
weakly absolutely convergent to $\mu^{(1)}_{\lambda, w_{(1)}}$ as
indicated in the $P^{(1)}(z_2)$-grading condition (or the
$L(0)$-semisimple $P^{(1)}(z_2)$-grading condition).  Suppose also
that for each $n \in \R$, $w'_{(4)} \in W'_4$ and $w_{(1)} \in W_1$,
the generalized $V$-submodule of $W^{(1)}_{\lambda, w_{(1)}}$
generated by $\lambda_{n}^{(1)}$ is a generalized $V$-submodule of
some object of $\mathcal{C}$ included in $(W_2 \otimes W_3)^*$.  Then
there exists a unique logarithmic intertwining operator (a unique
ordinary intertwining operator in the case that $\mathcal{C}$ is in
$\mathcal{M}_{sg}$) $\Y_{1}$ of type ${W_4\choose
W_1\,\,W_2\boxtimes_{P(z_2)} W_3}$ such that
\begin{eqnarray*}
\lefteqn{\langle w'_{(4)}, \Y^1(\Y^2(w_{(1)}, x_{0})w_{(2)}, x_{2})w_{(3)}\rangle
\lbar_{x_{0}=z_{1}-z_{2},\;x_{2}=z_{2}}}\nn
&&=\langle w'_{(4)}, \Y_{1}(w_{(1)}, x_{1})\Y_{\boxtimes_{P(z_2)}, 0}(w_{(2)}, x_{2})w_{(3)}
\rangle\lbar_{x_{1}=z_{1},\;x_{2}=z_{2}}
\end{eqnarray*}
for all
$w_{(1)}\in W_1$, $w_{(2)}\in W_2$, $w_{(3)}\in W_3$ and $w'_{(4)}\in
W'_4$.  In particular, the iterate of the logarithmic intertwining
operators (ordinary intertwining operators in the case that
$\mathcal{C}$ is in $\mathcal{M}_{sg}$) $\Y^{1}$ and $\Y^{2}$
evaluated at $z_{2}$ and $z_{1}-z_{2}$, respectively, can be expressed as a
product (with the intermediate generalized $V$-module
$W_2\boxtimes_{P(z_2)} W_3$) of logarithmic intertwining operators
(ordinary intertwining operators in the case that $\mathcal{C}$ is in
$\mathcal{M}_{sg}$) evaluated at $z_{1}$ and $z_{2}$.

\end{enumerate}
\end{corol}

In Section 10, we construct the associativity isomorphisms, under
certain assumptions: In addition to the assumptions above, we assume
that $\mathcal{C}$ is closed under images and that for some $z\in
\C^{\times}$ (and hence for every $z\in \C^{\times}$), $\mathcal{C}$
is closed under $P(z)$-tensor products; see Assumption
\ref{assum-assoc}.  Besides the convergence condition (Section 7), at
the end of Section 9 we introduce what we call the ``expansion
condition,'' which, roughly speaking, states that an element of
$(W_{1}\otimes W_{2}\otimes W_{3})^{*}$ obtained from a product or an
iterate of intertwining maps satisfies the $P^{(2)}(z)$- or
$P^{(1)}(z)$-local grading restriction condition, respectively, for
suitable $z\in \C^{\times}$, along with certain other ``minor''
conditions. Then we have:

\setcounter{section}{10}
\setcounter{rema}{2}
\begin{theo}
Assume that the convergence condition and the expansion condition for
intertwining maps in ${\cal C}$ both hold.  Let $z_1$, $z_2$ be
complex numbers satisfying
\[
|z_1|>|z_2|>|z_1-z_{2}|>0
\]
(so that in particular $z_1\neq 0$, $z_2\neq 0$ and $z_1\neq z_2$).
Then there exists a unique natural isomorphism
\[
\mathcal{A}_{P(z_{1}), P(z_{2})}^{P(z_{1}-z_{2}), P(z_{2})}:
\boxtimes_{P(z_{1})}\circ (1 \times \boxtimes_{P(z_2)}) \to
\boxtimes_{P(z_2)}\circ (\boxtimes_{P(z_1-z_2)}\times 1)
\]
such that for all $w_{(1)}\in W_1$, $w_{(2)}\in W_2$ and
$w_{(3)}\in W_3$, with $W_j$ objects of ${\cal C}$,
\[
\overline{\mathcal{A}_{P(z_{1}), P(z_{2})}^{P(z_{1}-z_{2}), P(z_{2})}}
(w_{(1)}\boxtimes_{P(z_1)}
(w_{(2)}\boxtimes_{P(z_2)} w_{(3)})) = (w_{(1)}\boxtimes_{P(z_1-z_2)}
w_{(2)})\boxtimes_{P(z_2)} w_{(3)},
\]
where for simplicity we use the same notation $\mathcal{A}_{P(z_{1}),
P(z_{2})}^{P(z_{1}-z_{2}), P(z_{2})}$ to denote the isomorphism of
generalized modules
\[
\mathcal{A}_{P(z_{1}), P(z_{2})}^{P(z_{1}-z_{2}), P(z_{2})}:
W_1\boxtimes_{P(z_{1})}
(W_2\boxtimes_{P(z_2)} W_3) \longrightarrow (W_1\boxtimes_{P(z_1-z_2)}
W_2)\boxtimes_{P(z_2)} W_3.
\]
\end{theo}

Here we are using the notation
\[
\overline{\eta}:\overline{W_1} \to \overline{W_2}
\]
to denote the natural extension of a map $\eta:W_1 \to W_2$ of
generalized modules to the (suitably defined) formal completions; such
natural extensions enter into many of the constructions in this work.

In Section 11, we give results which will allow us to verify the
convergence and expansion conditions. We need the ``convergence and
extension property'' for products or iterates and the ``convergence
and extension property without logarithms'' for products or
iterates. Here we only give the convergence and extension property for
products:

Given objects $W_1$, $W_2$, $W_3$, $W_4$ and $M_1$ of the category
${\cal C}$, let ${\cal Y}_1$ and ${\cal Y}_2$ be logarithmic
intertwining operators of types ${W_4}\choose {W_1M_1}$ and
${M_1}\choose {W_2W_3}$, respectively.

\begin{description}

\item[Convergence and extension property for products] For any
$\beta\in \tilde{A}$, there exists an integer $N_{\beta}$ depending
only on ${\cal Y}_1$, ${\cal Y}_2$ and $\beta$, and for any
weight-homogeneous elements $w_{(1)}\in (W_1)^{(\beta_{1})}$ and
$w_{(2)}\in (W_2)^{(\beta_{2})}$ ($\beta_{1}, \beta_{2} \in
\tilde{A}$) and any $w_{(3)}\in W_3$ and $w'_{(4)}\in W'_4$ such that
\[
\beta_{1}+\beta_{2}=-\beta,
\]
there exist $M\in{\mathbb N}$, $r_{k}, s_{k}\in {\mathbb R}$, $i_{k},
j_{k}\in \N$, $k=1,\dots,M$; $K \in \Z_{+}$ independent of $w_{(1)}$ and
$w_{(2)}$ such that each $i_{k} < K$; and analytic functions $f_{k}(z)$ on
$|z|<1$, $k=1, \dots, M$, satisfying
\[
\wt w_{(1)}+\wt w_{(2)}+s_{k}>N_{\beta}, \;\;\;k=1, \dots, M,
\]
such that
\[
\langle w'_{(4)}, {\cal Y}_1(w_{(1)}, x_2) {\cal Y}_2(w_{(2)},
x_2)w_{(3)}\rangle_{W_4} \lbar_{x_1= z_1, \;x_2=z_2}
\]
is absolutely convergent when $|z_1|>|z_2|>0$ and can be analytically
extended to the multivalued analytic function
\[
\sum_{k=1}^{M}z_2^{r_{k}}(z_1-z_2)^{s_{k}}(\log z_2)^{i_{k}}
(\log(z_1-z_2))^{j_{k}}f_{k}\left(\frac{z_1-z_2}{z_2}\right)
\]
(here $\log (z_{1}-z_{2})$ and $\log z_{2}$, and in particular, the
powers of the variables, mean the multivalued functions, not the
particular branch we have been using) in the region
$|z_2|>|z_1-z_2|>0$.

\end{description}

\setcounter{section}{11}
\setcounter{rema}{3}

\begin{theo}
Suppose that the following two conditions are satisfied:
\begin{enumerate}
\item Every finitely-generated lower bounded doubly-graded (as defined
in Section 11) generalized $V$-module is an object of ${\cal C}$ (or
every finitely-generated lower bounded doubly-graded $V$-module is an
object of ${\cal C}$, when ${\cal C}$ is in $\mathcal{M}_{sg}$).

\item The convergence and extension property for either products or
iterates holds in ${\cal C}$ (or the convergence and extension
property without logarithms for either products or iterates holds in
${\cal C}$, when ${\cal C}$ is in $\mathcal{M}_{sg}$).
\end{enumerate}
Then the convergence and expansion conditions for intertwining
maps in ${\cal C}$ both hold.
\end{theo}

In the following two results, we assume that that the grading abelian
groups $A$ and $\tilde{A}$ are trivial.  Set
\[
V_{+}=\coprod_{n>0}V_{(n)}.
\]
Let $W$ be a
generalized $V$-module and let 
\[
C_{1}(W)={\rm span}\{u_{-1}w\;|\; u\in V_{+},\;
w\in W\}.  
\]
If $W/C_{1}(W)$ is finite dimensional, we say that $W$ is
$C_{1}$-cofinite or satisfies the $C_{1}$-cofiniteness condition. If
for any $N\in {\mathbb R}$, $\coprod_{n<N}W_{[n]}$ is finite
dimensional, we say that $W$ is quasi-finite dimensional or satisfies
the quasi-finite-dimensionality condition.  The following result in
Section 11 allows us to verify the convergence and extension
properties and thus the convergence and expansion conditions:

\setcounter{rema}{5}

\begin{theo}
Let $W_{i}$ for $i=0, \dots, n+1$ be generalized $V$-modules
satisfying the $C_{1}$-cofiniteness condition and the
quasi-finite-dimensionality condition. Then for any $w'_{(0)}\in W'_{0}$, 
$w_{(1)}\in W_{1}, \dots,
w_{(n+1)}\in W_{n+1}$, there exist
\[
a_{k, \;l}(z_{1}, \dots, z_{n})\in {\mathbb C}[z_{1}^{\pm 1}, \dots,
z_{n}^{\pm 1}, (z_{1}-z_{2})^{-1}, (z_{1}-z_{3})^{-1}, \dots,
(z_{n-1}-z_{n})^{-1}],
\]
for $k=1, \dots, m$ and $l=1, \dots, n,$ such that the following 
holds: For any generalized
$V$-modules $\widetilde{W}_{1}, \dots, \widetilde{W}_{n-1}$, and any logarithmic 
intertwining operators 
\[
{\cal Y}_{1}, {\cal Y}_{2}, \dots, {\cal
Y}_{n-1}, {\cal Y}_{n}
\]
of types 
\[
{W_{0}\choose
W_{1}\widetilde{W}_{1}}, {\widetilde{W}_{1}\choose
W_{2}\widetilde{W}_{2}}, \dots, {\widetilde{W}_{n-2}\choose
W_{n-1}\widetilde{W}_{n-1}}, {\widetilde{W}_{n-1}\choose
W_{n}W_{n+1}},
\]
respectively, the series
\[
\langle w'_{(0)}, {\cal Y}_{1}(w_{(1)}, z_{1})\cdots {\cal Y}_{n}(w_{(n)},
z_{n})w_{(n+1)}\rangle
\]
satisfies the system of differential equations
\[
\frac{\partial^{m}\varphi}{\partial z_{l}^{m}}+ \sum_{k=1}^{m}
\iota_{|z_{1}|>\cdots >|z_{n}|>0}(a_{k,
\;l}(z_{1}, \dots, z_{n})) \frac{\partial^{m-k}\varphi}{\partial
z_{l}^{m-k}}=0,\;\;\;l=1, \dots, n
\]
in the region $|z_{1}|>\cdots >|z_{n}|>0$, where 
\[
\iota_{|z_{1}|>\cdots >|z_{n}|>0}(a_{k,
\;l}(z_{1}, \dots, z_{n}))
\]
for $k=1, \dots, m$ and $l=1, \dots, n$ are the (unique) Laurent
expansions of $a_{k, \;l}(z_{1}, \dots, z_{n})$ in the region
$|z_{1}|>\cdots >|z_{n}|>0$.  Moreover, for any set of possible
singular points of the system
\[
\frac{\partial^{m}\varphi}{\partial z_{l}^{m}}+ \sum_{k=1}^{m} a_{k,
\;l}(z_{1}, \dots, z_{n})\frac{\partial^{m-k}\varphi}{\partial
z_{l}^{m-k}}=0,\;\;\;l=1, \dots, n
\]
such that either $z_i = 0$ or $z_i = \infty$ for some $i$ or $z_i =
z_j$ for some $i \ne j$, the $a_{k, \;l}(z_{1}, \dots, z_{n})$ can be
chosen for $k=1, \dots, m$ and $l=1, \dots, n$ so that these singular
points are regular.
\end{theo}

Using this result, we prove the following:

\setcounter{rema}{7}

\begin{theo}
Suppose that all generalized $V$-modules in ${\cal C}$ satisfy the
$C_{1}$-cofiniteness condition and the quasi-finite-dimensionality
condition.  Then:

\begin{enumerate}

\item The convergence and extension properties for products
and iterates hold in ${\cal C}$. If $\mathcal{C}$ is in
$\mathcal{M}_{sg}$ and if every object of $\mathcal{C}$ is a direct sum
of irreducible objects of $\mathcal{C}$ and there are only finitely
many irreducible objects of $\mathcal{C}$ (up to equivalence), then
the convergence and extension properties without logarithms for
products and iterates hold in ${\cal C}$. 

\item For any $n\in{\mathbb Z}_+$, any objects $W_1, \dots, W_{n+1}$
and $\widetilde{W}_1, \dots, \widetilde{W}_{n-1}$ of $\mathcal{C}$,
any logarithmic intertwining operators
\[
{\cal Y}_{1}, {\cal Y}_{2}, \dots, {\cal
Y}_{n-1}, {\cal Y}_{n}
\]
of types 
\[
{W_{0}\choose
W_{1}\widetilde{W}_{1}}, {\widetilde{W}_{1}\choose
W_{2}\widetilde{W}_{2}}, \dots, {\widetilde{W}_{n-2}\choose
W_{n-1}\widetilde{W}_{n-1}}, {\widetilde{W}_{n-1}\choose
W_{n}W_{n+1}},
\]
respectively, and any $w_{(0)}'\in W_{0}'$, $w_{(1)}\in W_{1}, \dots,
w_{(n+1)}\in W_{n+1}$, the series
\begin{equation} 
\langle
w_{(0)}', {\cal Y}_{1}(w_{(1)}, z_{1})\cdots {\cal Y}_{n}(w_{(n)},
z_{n})w_{(n+1)}\rangle
\end{equation}
is absolutely convergent in the region $|z_{1}|>\cdots> |z_{n}|>0$ and
its sum can be analytically extended to a multivalued analytic
function on the region given by $z_{i}\ne 0$, $i=1, \dots, n$,
$z_{i}\ne z_{j}$, $i\ne j$, such that for any set of possible singular
points with either $z_{i}=0$, $z_{i}=\infty$ or $z_{i}= z_{j}$ for
$i\ne j$, this multivalued analytic function can be expanded near the
singularity as a series having the same form as the expansion near the
singular points of a solution of a system of differential equations
with regular singular points.

\end{enumerate}
\end{theo}

We now return to the assumptions before Theorem \ref{sys}, that is, we
do not assume that $A$ and $\tilde{A}$ are trivial.  To construct the
braided tensor category structure, we need more assumptions in
addition to those mentioned above, which are collected in Assumption
\ref{assum-assoc}.  We assume in addition that the M\"{o}bius or
conformal vertex algebra $V$, viewed as a $V$-module, is an object of
$\mathcal{C}$; and also that the product of three logarithmic
intertwining operators is absolutely convergent in a suitable region
and can be analytically extended to a multivalued analytic function,
admitting suitable expansions as series in powers of the variables and
their logarithms near its singularities (expansions that hold for
solutions of systems of differential equations with regular
singularities), on a suitable largest possible region containing the
original region for the convergence of the product. See Assumptions
\ref{assum-V} and \ref{assum-con} for the precise statements.  Under
these assumptions, we construct, in addition to the tensor product
bifunctor $\boxtimes=\boxtimes_{P(1)}$, a braiding isomorphism
$\mathcal{R}$, an associativity isomorphism $\mathcal{A}$ (for the
braided tensor category structure, different from the associativity
isomorphisms above), a left unit isomorphism $l$ and a right unit
isomorphism $r$.  The following main results of this work are given in
Section 12:

\setcounter{section}{12}
\setcounter{rema}{14}
\begin{theo}
Let $V$ be a M\"{o}bius or conformal vertex algebra and $\mathcal{C}$
a full subcategory of $\mathcal{M}_{sg}$ or $\mathcal{GM}_{sg}$
satisfying Assumptions \ref{assum-assoc}, \ref{assum-V} and
\ref{assum-con}.  Then the category $\mathcal{C}$, equipped with the
tensor product bifunctor $\boxtimes$, the unit object $V$, the
braiding isomorphism $\mathcal{R}$, the associativity isomorphism
$\mathcal{A}$, and the left and right unit isomorphisms $l$ and $r$,
is an additive braided monoidal category.
\end{theo}

\begin{corol}
If the category $\mathcal{C}$ is an abelian category, then
$\mathcal{C}$, equipped with the tensor product bifunctor $\boxtimes$,
the unit object $V$, the braiding isomorphism $\mathcal{R}$, the
associativity isomorphism $\mathcal{A}$, and the left and right unit
isomorphisms $l$ and $r$, is a braided tensor category.
\end{corol}

\newpage

\setcounter{section}{1}
\setcounter{equation}{0}
\setcounter{rema}{0}

\section{The setting: strongly graded conformal and M\"obius vertex
algebras and their generalized modules}

In this section we define and discuss the basic structures and
introduce some notation that will be used in this work. More
specifically, we first introduce the notions of ``conformal vertex
algebra'' and ``M\"obius vertex algebra.''  A conformal vertex algebra
is just a vertex algebra equipped with a conformal vector satisfying
the usual axioms; a M\"obius vertex algebra is a variant of a
``quasi-vertex operator algebra'' as in \cite{FHL}, with the
difference that the two grading restriction conditions in the
definition of vertex operator algebra are not required. We then define
the notion of module for each of these types of vertex algebra.
Relaxing the $L(0)$-semisimplicity in the definition of module we
obtain the notion of ``generalized module.''  Finally, we notice that
in order to have a contragredient functor on the module category under
consideration, we need to impose a stronger grading condition.  This
leads to the notions of ``strong gradedness'' of M\"obius vertex
algebras and their generalized modules. In this work we are mainly
interested in certain full subcategories of the category of strongly
graded generalized modules for certain strongly graded M\"obius vertex
algebras.

Throughout the work we shall assume some familiarity with the material
in \cite{B}, \cite{FLM2}, \cite{FHL}, \cite{DL} and \cite{LL}.

In particular, we recall the necessary basic material on ``formal
calculus,'' starting with the ``formal delta function.''  Formal
calculus will be needed throughout this work, and in fact, the theory
of formal calculus will be considerably developed, whenever new
formal-calculus ideas are needed for the formulations and for the
proofs of the results.

Throughout, we shall use the notation ${\mathbb N}$ for the nonnegative
integers and $\Z_{+}$ for the positive integers.

We shall continue to use the notational convention concerning formal
variables and complex variables given in Remark
\ref{formalandcomplexvariables}.  Recall {}from \cite{FLM2},
\cite{FHL} or \cite{LL} that the formal delta function is defined as
the formal series
\[
\delta(x)=\sum_{n\in {\mathbb Z}} x^n.
\]
in the formal variable $x$.  We will consistently use the {\em
binomial expansion convention}: For any complex number $\lambda$,
$(x+y)^\lambda$ is to be expanded as a formal series in nonnegative
integral powers of the second variable, i.e.,
\[
(x+y)^\lambda=\sum_{n\in {\mathbb N}} {\lambda \choose n} x^{\lambda -n}
y^n.
\]
Here $x$ or $y$ might be something other than a formal variable (or a
nonzero complex multiple of a formal variable); for instance, $x$ or
$y$ (but not both; this expansion is understood to be formal) might be
a nonzero complex number, or $x$ or $y$ might be some more complicated
object.  The use of the binomial expansion convention will be clear in
context.

Objects like $\delta(x)$ and $(x+y)^\lambda$ lie in spaces of formal
series.  Some of the spaces that we will use are, with $W$ a vector
space (over $\C$) and $x$ a formal variable:
\[
W[x]=\biggl\{\sum_{n\in \mathbb{N}}a_{n}x^{n}| a_{n}\in W,\;
\mbox{all but finitely many} \; a_{n} = 0\biggr\}
\]
(the space of formal polynomials with coefficients in $W$),
\[
W[x,x^{-1}]=\biggl\{\sum_{n\in \mathbb{Z}}a_{n}x^{n}| a_{n}\in W,\;
\mbox{all but finitely many} \; a_{n} = 0\biggr\}
\]
(the formal Laurent polynomials),
\[
W[[x]]=\biggl\{\sum_{n\in \mathbb{N}}a_{n}x^{n}| a_{n}\in W\;
\mbox{(with possibly infinitely many} \; a_{n} \; \mbox{not} \; 0)\biggr\}
\]
(the formal power series),
\[
W((x))=\biggl\{\sum_{n\in \mathbb{Z}}a_{n}x^{n}| a_{n}\in W,\;
a_{n} = 0 \; \mbox{for sufficiently small} \; n \biggr\}
\]
(the truncated formal Laurent series), and
\[
W[[x,x^{-1}]]=\biggl\{\sum_{n\in \mathbb{Z}}a_{n}x^{n}| a_{n}\in W\;
\mbox{(with possibly infinitely many} \; a_{n} \; \mbox{not} \; 0)\biggr\}
\]
(the formal Laurent series).  We will also need the space
\begin{equation}\label{formalserieswithcomplexpowers}
W\{ x\}=\biggl\{\sum_{n\in \mathbb{C}}a_{n}x^{n}| a_{n}\in W \;
\mbox{for} \; n\in {\mathbb C}\biggr\}
\end{equation}
as in \cite{FLM2}; here the powers of the formal variable are complex,
and the coefficients may all be nonzero.  We will also use analogues
of these spaces involving two or more formal variables.  Note that for
us, a ``formal power series'' involves only nonnegative integral
powers of the formal variable(s), and a ``formal Laurent series'' can
involve all the integral powers of the formal variable(s).

The following formal version of Taylor's theorem is easily verified by
direct expansion (see Proposition 8.3.1 of \cite{FLM2}):  For $f(x)
\in W\{ x\}$,
\begin{equation}\label{formalTaylortheorem}
e^{y \frac{d}{dx}}f(x) = f(x + y),
\end{equation}
where the exponential denotes the formal exponential series, and where
we are using the binomial expansion convention on the right-hand side.
It is important to note that this formula holds for arbitrary formal
series $f(x)$ with complex powers of $x$, where $f(x)$ need not be an
expansion in any sense of an analytic function (again, see Proposition
8.3.1 of \cite{FLM2}).

The formal delta function $\delta(x)$ has the following simple and
fundamental property: For any $f(x)\in W[x, x^{-1}]$,
\begin{equation}
f(x)\delta(x)=f(1)\delta(x).
\end{equation}
(Here we are taking the liberty of writing complex numbers to the
right of vectors in $W$.)  This is proved immediately by observing its
truth for $f(x)=x^n$ and then using linearity.  This property has many
important variants; in general, whenever an expression is multiplied
by the formal delta function, we may formally set the argument
appearing in the delta function equal to 1, provided that the relevant
algebraic expressions make sense.  For example, for any
$$X(x_{1},
x_{2})\in (\mbox{End }W)[[x_{1}, x_{1}^{-1}, x_{2}, x_{2}^{-1}]]$$
such that
\begin{equation}\label{limx1approachesx2}
\lim_{x_{1}\to x_{2}}X(x_{1}, x_{2})=X(x_{1},
x_{2})\lbar_{x_{1}=x_{2}}
\end{equation}
exists, we have
\begin{equation}\label{Xx1x2=Xx2x2}
X(x_{1}, x_{2})\delta\left(\frac{x_{1}}{x_{2}}\right)=X(x_{2}, x_{2})
\delta\left(\frac{x_{1}}{x_{2}}\right).
\end{equation}
The existence of the ``algebraic limit'' defined in
(\ref{limx1approachesx2}) means that for an arbitrary vector $w\in W$,
the coefficient of each power of $x_{2}$ in the formal expansion
$X(x_{1}, x_{2})w\lbar_{x_{1}=x_{2}}$ is a finite sum.  In general,
the existence of such ``algebraic limits,'' and also such products of
formal sums, always means that the coefficient of each monomial in the
relevant formal variables gives a finite sum.  Often, proving the
existence of the relevant algebraic limits (or products) is a much
more subtle matter than computing such limits (or products), just as
in analysis.  (In this work, we will typically use ``substitution
notation'' like $\lbar_{x_{1}=x_{2}}$ or $X(x_{2},x_{2})$ rather than
the formal limit notation on the left-hand side of
(\ref{limx1approachesx2}).)  Below, we will give a more sophisticated
analogue of the delta-function substitution principle
(\ref{Xx1x2=Xx2x2}), an analogue that we will need in this work.

This analogue, and in fact, many fundamental principles of vertex
operator algebra theory, are based on certain delta-function
expressions of the following type, involving three (commuting and
independent, as usual) formal variables:
\[
x_{0}^{-1}\delta\left(\frac{x_{1}-x_{2}}{x_{0}}\right)=\sum_{n\in
{\mathbb Z}}
\frac{(x_{1}-x_{2})^{n}}{x_{0}^{n+1}}=\sum_{m\in {\mathbb N},\; n\in
{\mathbb Z}}
(-1)^{m}{{n}\choose {m}} x_{0}^{-n-1}x_{1}^{n-m}x_{2}^{m};
\]
here the binomial expansion convention is of course being used.

The following important identities involving such three-variable
delta-function expressions are easily proved (see \cite{FLM2} or
\cite{LL}, where extensive motivation for these formulas is also
given):
\begin{equation}\label{2termdeltarelation}
x_{2}^{-1}\delta\left(\frac{x_{1}-x_{0}}{x_{2}}\right)=
x_{1}^{-1}\delta\left(\frac{x_{2}+x_{0}}{x_{1}}\right),
\end{equation}
\begin{equation}\label{3termdeltarelation}
x_{0}^{-1}\delta\left(\frac{x_{1}-x_{2}}{x_{0}}\right)-
x_{0}^{-1}\delta\left(\frac{x_{2}-x_{1}}{-x_{0}}\right)=
x_{2}^{-1}\delta\left(\frac{x_{1}-x_{0}}{x_{2}}\right).
\end{equation}
Note that the three terms in (\ref{3termdeltarelation}) involve
nonnegative integral powers of $x_2$, $x_1$ and $x_0$, respectively.
In particular, the two terms on the left-hand side of
(\ref{3termdeltarelation}) are unequal formal Laurent series in three
variables, even though they might appear equal at first glance.  We
shall use these two identities extensively.

\begin{rema}\label{deltafunctionsubstitutionremark}{\rm
Here is the useful analogue, mentioned above, of the delta-function
substitution principle (\ref{Xx1x2=Xx2x2}):  Let
\begin{equation}
f(x_1,x_2,y) \in (\mbox{End }W)[[x_{1}, x_{1}^{-1}, x_{2},
x_{2}^{-1},y,y^{-1}]]
\end{equation}
be such that
\begin{equation}
\lim_{x_{1}\to x_{2}}f(x_1,x_2,y) \; \mbox{exists}
\end{equation}
and such that for any $w \in W$,
\begin{equation}
f(x_1,x_2,y)w \in W[[x_1,x_{1}^{-1},x_2,x_{2}^{-1}]]((y)).
\end{equation}
Then
\begin{equation}\label{deltafunctionsubstitutionformula}
x_{1}^{-1}\delta\left(\frac{x_{2}-y}{x_{1}}\right)f(x_1,x_2,y)=
x_{1}^{-1}\delta\left(\frac{x_{2}-y}{x_{1}}\right)f(x_2-y,x_2,y).
\end{equation}
For this principle, see Remark 2.3.25 of \cite{LL}, where the proof is
also presented.}
\end{rema}

The following formal residue notation will be useful: For
\[
f(x) = \sum_{n\in \mathbb{C}}a_{n}x^{n} \in W\{x\}
\]
(note that the powers of $x$ need not be integral),
\[
{\res_x}f(x) = a_{-1}.
\]
For instance, for the expression in (\ref{2termdeltarelation}),
\begin{equation}
\res_{x_2}x_{2}^{-1}\delta\left(\frac{x_{1}-x_{0}}{x_{2}}\right)=1.
\end{equation}

For a vector space $W$, we will denote its vector space dual by $W^*$
($=\hom_{\mathbb C}(W,{\mathbb C})$), and we will use the notation
$\langle\cdot,\cdot\rangle_W$, or $\langle\cdot,\cdot\rangle$ if the
underlying space $W$ is clear, for the canonical pairing between $W^*$
and $W$.

We will use the following version of the notion of ``conformal vertex
algebra'': A conformal vertex algebra is a vertex algebra (in the
sense of Borcherds \cite{B}; see \cite{LL}) equipped with a ${\mathbb
Z}$-grading and with a conformal vector satisfying the usual
compatibility conditions.  Specifically:

\begin{defi}\label{cva}
{\rm A {\it conformal vertex algebra} is a ${\mathbb Z}$-graded vector
space
\begin{equation}\label{Vgrading}
V=\coprod_{n\in{\mathbb Z}} V_{(n)}
\end{equation}
(for $v\in V_{(n)}$, we say the {\it weight} of $v$ is $n$ and we write
$\mbox{wt}\, v=n$)
equipped with a linear map $V\otimes V\to V[[x,
x^{-1}]]$, or equivalently,
\begin{eqnarray}\label{YforV}
V&\to&({\rm End}\; V)[[x, x^{-1}]] \nno\\
v&\mapsto& Y(v,
x)={\displaystyle \sum_{n\in{\mathbb Z}}}v_{n}x^{-n-1} \;\;(
\mbox{where }v_{n}\in{\rm End} \;V),
\end{eqnarray}
$Y(v, x)$ denoting the {\it vertex operator associated with} $v$, and
equipped also with two distinguished vectors ${\bf 1}\in V_{(0)}$ (the
{\it vacuum vector}) and $\omega\in V_{(2)}$ (the {\it conformal
vector}), satisfying the following conditions for $u,v \in V$:
the {\it lower truncation condition}:
\begin{equation}\label{ltc}
u_{n}v=0\;\;\mbox{ for }n\mbox{ sufficiently large}
\end{equation}
(or equivalently, $Y(u, x)v\in V((x))$); the {\it vacuum property}:
\begin{equation}\label{1left}
Y({\bf 1}, x)=1_V;
\end{equation}
the {\it creation property}:
\begin{equation}\label{1right}
Y(v, x){\bf 1} \in V[[x]]\;\;\mbox{ and }\;\lim_{x\rightarrow
0}Y(v, x){\bf 1}=v
\end{equation}
(that is, $Y(v, x){\bf 1}$ involves only nonnegative integral powers
of $x$ and the constant term is $v$); the {\it Jacobi identity} (the
main axiom):
\begin{eqnarray}
&x_0^{-1}\delta
\bigg({\displaystyle\frac{x_1-x_2}{x_0}}\bigg)Y(u, x_1)Y(v,
x_2)-x_0^{-1} \delta
\bigg({\displaystyle\frac{x_2-x_1}{-x_0}}\bigg)Y(v, x_2)Y(u,
x_1)&\nno \\ &=x_2^{-1} \delta
\bigg({\displaystyle\frac{x_1-x_0}{x_2}}\bigg)Y(Y(u, x_0)v,
x_2)&\label{Jacobi}
\end{eqnarray}
(note that when each expression in (\ref{Jacobi}) is applied to any
element of $V$, the coefficient of each monomial in the formal
variables is a finite sum; on the right-hand side, the notation
$Y(\cdot, x_2)$ is understood to be extended in the obvious way to
$V[[x_0, x^{-1}_0]]$); the {\em Virasoro algebra relations}:
\begin{equation}\label{vir1}
[L(m), L(n)]=(m-n)L(m+n)+{\displaystyle\frac{1}{12}}
(m^3-m)\delta_{n+m,0}c
\end{equation}
for $m, n \in {\mathbb Z}$, where
\begin{equation}\label{vir2}
L(n)=\omega _{n+1}\;\; \mbox{ for } \;  n\in{\mathbb Z}, \;\;
\mbox{i.e.},\;\;Y(\omega, x)=\sum_{n\in{\mathbb Z}}L(n)x^{-n-2},
\end{equation}
\begin{equation}\label{vir3}
c\in {\mathbb C}
\end{equation}
(the {\it central charge} or {\it rank} of $V$);
\begin{equation}\label{L-1derivativeproperty}
{\displaystyle \frac{d}{dx}}Y(v,
x)=Y(L(-1)v, x)
\end{equation}
(the {\it  $L(-1)$-derivative property}); and
\begin{equation}\label{L0gradingproperty}
L(0)v=nv=(\wt v)v \;\; \mbox{ for }\; n\in {\mathbb Z}\; \mbox{ and }\;
v\in V_{(n)}.
\end{equation}
}
\end{defi}

This completes the definition of the notion of conformal vertex
algebra.  We will denote such a conformal vertex algebra by $(V,Y,{\bf
1},\omega)$ or simply by $V$.

The only difference between the definition of conformal vertex algebra
and the definition of {\it vertex operator algebra} (in the sense of
\cite{FLM2} and \cite{FHL}) is that a vertex operator algebra $V$
also satisfies the two {\it grading restriction conditions}
\begin{equation}\label{gr1}
V_{(n)}=0 \;\; \mbox{ for }n\mbox{ sufficiently negative},
\end{equation}
and
\begin{equation}\label{gr2}
{\rm dim} \; V_{(n)} < \infty \;\;\mbox{ for }n \in {\mathbb Z}.
\end{equation}
(As we mentioned above, a {\it vertex algebra} is the same thing as a
conformal vertex algebra but without the assumptions of a grading or a
conformal vector, or, of course, the $L(n)$'s.)

\begin{rema}\label{va>cva}{\rm
Of course, not every vertex algebra is conformal. For example, it is
well known \cite{B} that any commutative associative algebra $A$ with
unit $1$, together with a derivation $D:A\to A$ can be equipped with a
vertex algebra structure, by:
\[
Y(\cdot,x)\cdot : A\times A\to A[[x]],\;\;Y(a,x)b=(e^{xD}a)b,
\]
and ${\bf 1}=1$.  In particular, $u_n=0$ for any $u\in A$ and $n\geq
0$.  If $\omega$ is a conformal vector for such a vertex algebra, then
for any $u\in A$, $Du=u_{-2}{\bf 1}=L(-1)u$ {}from (\ref{1right}) and
(\ref{L-1derivativeproperty}), so $D=L(-1)=\omega_0$, which equals $0$
because $\omega=L(0)\omega/2=\omega_1\omega/2=0$.  Thus a vertex
algebra constructed {}from a commutative associative algebra with
nonzero derivation in this way cannot be conformal.  }
\end{rema}

\begin{rema}\label{motivate-Mobius}{\rm
The theory of vertex tensor categories inherently uses the whole
moduli space of spheres with two positively oriented punctures and one
negatively oriented puncture (and in fact, more generally, with
arbitrary numbers of positively oriented punctures and one negatively
oriented puncture) equipped with general (analytic) local coordinates
vanishing at the punctures.  Because of the analytic local
coordinates, our constructions require certain conditions on the
Virasoro algebra operators.  However, recalling the definition of the moduli
space elements $P(z)$ {}from Section 1.4, we point out that if we
restrict our attention to elements of the moduli space of only the
type $P(z)$, then the relevant operations of sewing and subsequently
decomposing Riemann spheres continue to yield spheres of the same
type, and rather than general conformal transformations around the
punctures, only M\"obius (projective) transformations around the
punctures are needed.  This makes it possible to develop the essential
structure of our tensor product theory by working entirely with
spheres of this special type; the general vertex tensor category
theory then follows {}from the structure thus developed.  This is why,
in the present work, we are focusing on the theory of $P(z)$-tensor
products.  Correspondingly, it turns out that it is very natural for
us to consider, along with the notion of conformal vertex algebra
(Definition \ref{cva}), a weaker notion of vertex algebra
involving only the three-dimensional subalgebra of the Virasoro
algebra corresponding to the group of M\"obius transformations.  That
is, instead of requiring an action of the whole Virasoro algebra, we
use only the action of the Lie algebra ${\mathfrak s}{\mathfrak l}(2)$
generated by $L(-1)$, $L(0)$ and $L(1)$.  Thus we get a notion
essentially identical to the notion of ``quasi-vertex operator
algebra'' in \cite{FHL}; the reason for focusing on this notion here
is the same as the reason why it was considered in \cite{FHL}.  Here
we designate this notion by the term ``M\"obius vertex algebra''; the
only difference between the definition of M\"obius vertex algebra and
the definition of quasi-vertex operator algebra \cite{FHL} is that a
quasi-vertex operator algebra $V$ also satisfies the two grading
restriction conditions (\ref{gr1}) and (\ref{gr2}).}
\end{rema}

Thus we formulate:

\begin{defi}\label{mobdef}{\rm
The notion of {\it M\"obius vertex algebra} is defined in the same way
as that of conformal vertex algebra except that in addition to the
data and axioms concerning $V$, $Y$ and ${\bf 1}$ (through
(\ref{Jacobi}) in Definition \ref{cva}), we assume (in place of
the existence of the conformal vector $\omega$ and the Virasoro
algebra conditions (\ref{vir1}), (\ref{vir2}) and (\ref{vir3})) the
following: We have a representation $\rho$ of ${\mathfrak s}{\mathfrak l}(2)$
on $V$ given by
\begin{equation}\label{Lrho}
L(j)=\rho(L_j),\;\;j=0,\pm 1,
\end{equation}
where $\{L_{-1}, L_0, L_1\}$ is a basis of ${\mathfrak s}{\mathfrak l}(2)$
with Lie brackets
\begin{equation}\label{L_*}
[L_0, L_{-1}]=L_{-1},\;\;[L_0,L_1]=-L_1,{\rm\;\;and\;\;}
[L_{-1},L_1]=-2L_0,
\end{equation}
and the following conditions hold for $v \in V$:
\begin{equation}\label{sl2-1}
[L(-1), Y(v,x)]=Y(L(-1)v,x),
\end{equation}
\begin{equation}\label{sl2-2}
[L(0), Y(v,x)]=Y(L(0)v,x)+xY(L(-1)v,x),
\end{equation}
\begin{equation}\label{sl2-3}
[L(1),Y(v,x)]=Y(L(1)v,x)+2xY(L(0)v,x)+x^2Y(L(-1)v,x),
\end{equation}
and also, (\ref{L-1derivativeproperty}) and
(\ref{L0gradingproperty}).  Of course, (\ref{sl2-1})--(\ref{sl2-3})
can be written as
\begin{eqnarray}\label{sl2-all}
[L(j),Y(v,x)]&=&\sum_{k=0}^{j+1}{j+1\choose k}x^kY(L(j-k)v,x)\nno\\
&=&\sum_{k=0}^{j+1}{j+1\choose k}x^{j+1-k}Y(L(k-1)v,x)
\end{eqnarray}
for $j=0,\pm 1$.}
\end{defi}

We will denote such a M\"obius vertex algebra by $(V,Y,{\bf 1}, \rho)$
or simply by $V$. Note that there is no notion of central charge (or
rank) for a M\"obius vertex algebra.  Also, a conformal vertex algebra
can certainly be viewed as a M\"obius vertex algebra in the obvious
way.  (Of course, a conformal vertex algebra could have other ${\mathfrak
s}{\mathfrak l}(2)$-structures making it a M\"obius vertex algebra in a
different way.)

\begin{rema}{\rm
By (\ref{Lrho}) and (\ref{L_*}) we have $[L(0), L(j)]=-jL(j)$ for
$j=0,\pm 1$. Hence
\begin{equation}\label{degL(j)}
L(j)V_{(n)}\subset V_{(n-j)},\;\; \mbox{ for }\;j=0,\pm1.
\end{equation}
Moreover, {}from (\ref{sl2-1}), (\ref{sl2-2}) and (\ref{sl2-3}) with
$v={\bf 1}$ we get, by (\ref{1left}) and (\ref{1right}),
\[
L(j){\bf 1}=0\;\; \mbox{ for } \; j=0,\pm 1.
\]
}
\end{rema}

\begin{rema}{\rm
Not every M\"obius vertex algebra is conformal. As an example, take
the commutative associative algebra ${\mathbb C}[t]$ with derivation
$D=-d/dt$, and form a vertex algebra as in Remark \ref{va>cva}. By
Remark \ref{va>cva}, this vertex algebra is not conformal. However,
define linear operators
\[
L(-1)=D,\quad L(0)=tD,\quad L(1)=t^2D
\]
on ${\mathbb C}[t]$. Then it is straightforward to verify that ${\mathbb
C}[t]$ becomes a M\"obius vertex algebra with these operators giving a
representation of ${\mathfrak s}{\mathfrak l}(2)$ having the desired
properties and with the ${\mathbb Z}$-grading (by nonpositive integers)
given by the eigenspace decomposition with respect to $L(0)$.}
\end{rema}

\begin{rema}{\rm
It is also easy to see that not every vertex algebra is M\"obius. For
example, take the two-dimensional commutative associative algebra
$A={\mathbb C}1\oplus{\mathbb C}a$ with $1$ as identity and $a^2=0$. The
linear operator $D$ defined by $D(1)=0$, $D(a)=a$ is a nonzero
derivation of $A$. Hence $A$ has a vertex algebra structure by Remark
\ref{va>cva}. Now if it is a module for ${\mathfrak s}{\mathfrak l}(2)$ as in
Definition \ref{mobdef}, since $A$ is two-dimensional and $L(0)1=0$,
$L(0)$ must act as $0$. But then $D=L(-1)=[L(0),L(-1)]=0$, a
contradiction.  }
\end{rema}

A module for a conformal vertex algebra $V$ is a module for $V$ viewed
as a vertex algebra such that the conformal element acts in the same
way as in the definition of vertex operator algebra. More precisely:

\begin{defi}\label{cvamodule}
{\rm  Given a conformal vertex algebra
$(V,Y,{\bf 1},\omega)$,  a {\it module} for $V$ is a ${\mathbb C}$-graded
vector space
\begin{equation}\label{Wgrading}
W=\coprod_{n\in{\mathbb C}} W_{(n)}
\end{equation}
(graded by {\it weights}) equipped with a linear map
$V\otimes W \rightarrow W[[x,x^{-1}]]$, or equivalently,
\begin{eqnarray}\label{YforW}
V &\rightarrow & (\mbox{End}\ W)[[x,x^{-1}]] \nno\\
v&\mapsto & Y(v,x) =\sum_{n\in {\mathbb Z}}v_nx^{-n-1}\;\;\;
(\mbox{where}\;\; v_n \in \mbox{End}\ W)
\end{eqnarray}
(note that the sum is over ${\mathbb Z}$, not ${\mathbb C}$), $Y(v,x)$
denoting the {\it vertex operator on $W$ associated with $v$}, such
that all the defining properties of a conformal vertex algebra that
make sense hold.  That is, the following conditions are satisfied: the
lower truncation condition: for $v \in V$ and $w \in W$,
\begin{equation}\label{ltc-w}
v_nw = 0 \;\;\mbox{ for }\;n\;\mbox{ sufficiently large}
\end{equation}
(or equivalently, $Y(v, x)w\in W((x))$); the vacuum property:
\begin{equation}\label{m-1left}
Y(\mbox{\bf 1},x) = 1_W;
\end{equation}
the Jacobi identity for vertex operators on $W$: for $u, v \in V$,
\begin{eqnarray}\label{m-Jacobi}
&{\dps x^{-1}_0\delta \bigg( {x_1-x_2\over x_0}\bigg)
Y(u,x_1)Y(v,x_2) - x^{-1}_0\delta \bigg( {x_2-x_1\over -x_0}\bigg)
Y(v,x_2)Y(u,x_1)}&\nno\\
&{\dps = x^{-1}_2\delta \bigg( {x_1-x_0\over x_2}\bigg)
Y(Y(u,x_0)v,x_2)}
\end{eqnarray}
(note that on the right-hand side, $Y(u,x_0)$ is the operator on $V$
associated with $u$); the Virasoro algebra relations on $W$ with
scalar $c$ equal to the central charge of $V$:
\begin{equation}\label{m-vir1}
[L(m), L(n)]=(m-n)L(m+n)+{\displaystyle\frac{1}{12}}
(m^3-m)\delta_{n+m,0}c
\end{equation}
for $m,n \in {\mathbb Z}$, where
\begin{equation}\label{m-vir2}
L(n)=\omega _{n+1}\;\; \mbox{ for }n\in{\mathbb Z}, \;\;{\rm
i.e.},\;\;Y(\omega, x)=\sum_{n\in{\mathbb Z}}L(n)x^{-n-2};
\end{equation}
\begin{equation}\label{L-1}
\displaystyle \frac{d}{dx}Y(v, x)=Y(L(-1)v, x)
\end{equation}
(the $L(-1)$-derivative property); and
\begin{equation}\label{wl0}
(L(0)-n)w=0\;\;\mbox{ for }\;n\in {\mathbb C}\;\mbox{ and }\;w\in W_{(n)}.
\end{equation}
}
\end{defi}

This completes the definition of the notion of module for a conformal
vertex algebra.

\begin{rema}\label{virrelationsformodule}{\rm
The Virasoro algebra relations (\ref{m-vir1}) for a module action follow
{}from the corresponding relations (\ref{vir1}) for $V$ together with the
Jacobi identities (\ref{Jacobi}) and (\ref{m-Jacobi}) and the
$L(-1)$-derivative properties (\ref{L-1derivativeproperty}) and
(\ref{L-1}), as we recall {}from (for example) \cite{FHL} or \cite{LL}.}
\end{rema}

We also have:

\begin{defi}\label{moduleMobius}{\rm
The notion of {\it module} for a M\"obius vertex algebra is defined in
the same way as that of module for a conformal vertex algebra except
that in addition to the data and axioms concerning $W$ and $Y$
(through (\ref{m-Jacobi}) in Definition \ref{cvamodule}), we assume
(in place of the Virasoro algebra conditions (\ref{m-vir1}) and
(\ref{m-vir2})) a representation $\rho$ of ${\mathfrak s}{\mathfrak l}(2)$ on
$W$ given by (\ref{Lrho}) and the conditions (\ref{sl2-1}),
(\ref{sl2-2}) and (\ref{sl2-3}), for operators acting on $W$, and
also, (\ref{L-1}) and (\ref{wl0}).  }
\end{defi}

In addition to modules, we have the following notion of {\em
generalized module} (or {\em logarithmic module}, as in, for example,
\cite{Mi}):

\begin{defi}\label{definitionofgeneralizedmodule}{\rm
A {\it generalized module} for a conformal (respectively, M\"obius)
vertex algebra is defined in the same way as a module for a conformal
(respectively, M\"obius) vertex algebra except that in the grading
(\ref{Wgrading}), each space $W_{(n)}$ is replaced by $W_{[n]}$, where
$W_{[n]}$ is the generalized $L(0)$-eigenspace corresponding to the
(generalized) eigenvalue $n\in {\mathbb C}$; that is, (\ref{Wgrading})
and (\ref{wl0}) in the definition are replaced by
\begin{equation}\label{Wgeneralizedgrading}
W=\coprod_{n\in{\mathbb C}} W_{[n]}
\end{equation}
and
\begin{equation}\label{gerwt}
\mbox{for }\;n\in {\mathbb C}\;\mbox{ and }\;w\in W_{[n]},\;
(L(0)-n)^mw=0 \;\mbox{ for }\;m\in {\mathbb N}\;\mbox{ sufficiently
large},
\end{equation}
respectively.  For $w \in W_{[n]}$, we still write $\wt w = n$ for the
(generalized) weight of $w$. }
\end{defi}

We will denote such a module or generalized module just defined by
$(W,Y)$, or sometimes by $(W,Y_W)$ or simply by $W$. We will use the
notation
\begin{equation}\label{pi_n}
\pi_n: W\to W_{[n]}
\end{equation}
for the projection {}from $W$ to its subspace of (generalized) weight
$n$, and for its natural extensions to spaces of formal series with
coefficients in $W$. In either the conformal or M\"obius case, a
module is of course a generalized module.

\begin{rema}\label{generalizedeigenspacedecomp}
{\rm For any vector space $U$ on which an operator, say, $L(0)$, acts
in such a way that
\begin{equation}\label{U=directsum}
U=\coprod_{n\in{\mathbb C}} U_{[n]}
\end{equation}
where for $n\in {\mathbb C}$,
\[
U_{[n]} = \{ u \in U | (L(0)-n)^m u=0 \;\mbox{ for }\;m\in {\mathbb
N}\;\mbox{ sufficiently large} \},
\]
we shall typically use the same projection notation
\begin{equation}
\pi_n: U\to U_{[n]}
\end{equation}
as in (\ref{pi_n}).  If instead of (\ref{U=directsum}) we have only
\[
U=\sum_{n\in{\mathbb C}} U_{[n]},
\]
then in fact this sum is indeed direct, and for any $L(0)$-stable
subspace $T$ of $U$, we have
\[
T=\coprod_{n\in{\mathbb C}} T_{[n]}
\]
(as with ordinary rather than generalized eigenspaces).}
\end{rema}

\begin{rema}\label{modulesaremodules}{\rm
A module for a conformal vertex algebra $V$ is obviously again a
module for $V$ viewed as a M\"obius vertex algebra, and conversely, a
module for $V$ viewed as a M\"obius vertex algebra is a module for $V$
viewed as a conformal vertex algebra, by Remark
\ref{virrelationsformodule}.  Similarly, the generalized modules for a
conformal vertex algebra $V$ are exactly the generalized modules for
$V$ viewed as a M\"obius vertex algebra.}
\end{rema}

\begin{rema}{\rm
A conformal or M\"obius vertex algebra is a module for itself (and 
in particular, a generalized module for itself).}
\end{rema}

\begin{rema}{\rm
In either the conformal or M\"obius vertex algebra case, we have the
obvious notions of $V$-{\em module homomorphism}, {\em submodule},
{\em quotient module}, and so on; in particular, homomorphisms are
understood to be grading-preserving.  We sometimes write the vector
space of (generalized-) module maps (homomorphisms) $W_1 \to W_2$ for
(generalized) $V$-modules $W_1$ and $W_2$ as $\hom_{V}(W_1,W_2)$.  }
\end{rema}

\begin{rema}{\rm
We have chosen the name ``generalized module'' here because the vector
space underlying the module is graded by generalized
eigenvalues. (This notion is different {}from the notion of
``generalized module'' used in \cite{tensor1}. A generalized module
for a vertex operator algebra $V$ as defined in, for example,
Definition 2.11 of \cite{tensor1} is precisely a module for $V$ viewed
as a conformal vertex algebra.)  }
\end{rema}

We will use the following notion of (formal algebraic) completion of a
generalized module:

\begin{defi}\label{Wbardef}{\rm
Let $W=\coprod_{n\in{\mathbb C}}W_{[n]}$ be a generalized module for a
M\"obius (or conformal) vertex algebra. We denote by $\overline{W}$
the (formal) completion of $W$ with respect to the ${\mathbb
C}$-grading, that is,
\begin{equation}\label{Wbar}
\overline{W}=\prod _{n\in {\mathbb C}} W_{[n]}.
\end{equation}
We will use the same notation $\overline{U}$ for any ${\mathbb
C}$-graded subspace $U$ of $W$.  We will continue to use the notation
$\pi_n$ for the projection {}from $\overline{W}$ to $W_{[n]}$:
\[
\pi_n: \overline{W} \to W_{[n]}.
\]
We will also continue to use the notation
$\langle\cdot,\cdot\rangle_W$, or $\langle\cdot,\cdot\rangle$ if the
underlying space is clear, for the canonical pairing between the
subspace $\coprod _{n\in {\mathbb C}} (W_{[n]})^*$ of $W^*$, and
$\overline{W}$.  We are of course viewing $(W_{[n]})^*$ as embedded in
$W^*$ in the natural way, that is, for $w^*\in (W_{[n]})^*$,
\begin{equation}\label{Wnstar}
\langle w^*, w\rangle_W=\langle w^*, w_n\rangle_{W_{[n]}}
\end{equation}
for any $w=\sum_{m\in {\mathbb C}} w_m$ (finite sum) in $W$, where
$w_m\in W_{[m]}$.}
\end{defi}

The following weight formula holds for generalized modules,
generalizing the corresponding formula in the module case
(cf.\ \cite{Mi}):

\begin{propo}\label{gweight}
Let $W$ be a generalized module for a M\"obius (or conformal) vertex
algebra $V$. Let both $v\in V$ and $w\in W$ be homogeneous. Then
\begin{eqnarray}
&\wt(v_n w)=\wt v +\wt w-n-1 \;\; \mbox{ for any }\; n\in {\mathbb Z},
&\label{set:wtvn}\\
&\wt(L(j)w)=\wt w-j \;\;\mbox{ for }\; j=0,\pm 1.&\label{set:wtsl2}
\end{eqnarray}
\end{propo}
\pf Applying the $L(-1)$-derivative property (\ref{L-1}) to formula
(\ref{sl2-2}), with the operators acting on $W$, and extracting the
coefficient of $x^{-n-1}$, we obtain:
\begin{equation}\label{[L(0),v_n]}
[L(0), v_n]=(L(0)\,v)_n+(-n-1)v_n.
\end{equation}
This can be written as
\[
(L(0)-(\wt v-n-1))v_n=v_n L(0),
\]
and so we have
\[
(L(0)-(\wt v+m-n-1))v_n=v_n (L(0)-m)
\]
for any $m\in {\mathbb C}$. Applying this repeatedly we get
\[
(L(0)-(\wt v+m-n-1))^t v_n=v_n (L(0)-m)^t
\]
for any $t\in {\mathbb N}$, $m\in {\mathbb C}$, and (\ref{set:wtvn})
follows.

For (\ref{set:wtsl2}), since as operators acting on $W$ we have
\begin{equation}\label{set:0j}
[L(0),L(j)]=-jL(j)
\end{equation}
for $j=0,\pm 1$, we get $(L(0)+j)L(j)=L(j)L(0)$ so that
\[
(L(0)-m+j)L(j)=L(j)(L(0)-m)
\]
for any $m\in {\mathbb C}$. Thus
\[
(L(0)-m+j)^tL(j)=L(j)(L(0)-m)^t
\]
for any $t\in {\mathbb N}$, $m\in {\mathbb C}$, and (\ref{set:wtsl2})
follows. \epf

\begin{rema}\label{congruent}{\rm
{}From Proposition \ref{gweight} we see that a generalized $V$-module
$W$ decomposes into submodules corresponding to the congruence classes
of its weights modulo $\Z$: For $\mu \in \C/\Z$, let
\begin{equation}
W_{[\mu]} = \coprod_{\bar n = \mu} W_{[n]},
\end{equation}
where $\bar n$ denotes the equivalence class of $n \in \C$ in
$\C/\Z$.  Then
\begin{equation}
W = \coprod_{\mu \in \C/\Z} W_{[\mu]}
\end{equation}
and each $W_{[\mu]}$ is a $V$-submodule of $W$.  Thus if a generalized
module $W$ is indecomposable (in particular, if it is irreducible),
then all complex numbers $n$ for which $W_{[n]}\neq 0$ are congruent
modulo ${\mathbb Z}$ to each other. }
\end{rema}

\begin{rema}\label{set:L(0)s}{\rm
Let $W$ be a generalized module for a M\"obius (or conformal) vertex
algebra $V$.  We consider the ``semisimple part'' $L(0)_s\in
\mbox{End}\ W$ of the operator $L(0)$:
\[
L(0)_sw=nw \;\; \mbox{ for } \; w\in W_{[n]},\; n \in {\mathbb C}.
\]
Then on $W$ we have
\begin{eqnarray}
&{}[L(0)_s, v_n]=[L(0), v_n]\;\;\mbox{ for all }\;v\in V\;\mbox{ and }
\; n\in {\mathbb Z};&\label{L0s,vn}\\
&{}[L(0)_s, L(j)]=[L(0), L(j)]\;\;\mbox{ for }\; j=0,\pm 1.
&\label{L0s,Lj}
\end{eqnarray}
Indeed, for homogeneous elements $v\in V$ and $w\in W$,
(\ref{set:wtvn}) and (\ref{[L(0),v_n]}) imply that
\begin{eqnarray*}
[L(0)_s, v_n]w&=&L(0)_s(v_nw)-v_n(L(0)_sw)\\
&=&(\wt v +\wt w-n-1)v_nw-(\wt w)v_nw\\
&=&(\wt v)v_nw+(-n-1)v_nw\\
&=&(L(0)v)_nw+(-n-1)v_nw\\
&=&[L(0), v_n]w.
\end{eqnarray*}
Similarly, for any homogeneous element $w\in W$ and $j=0,\pm 1$,
(\ref{set:wtsl2}) and (\ref{set:0j}) imply that
\begin{eqnarray*}
[L(0)_s, L(j)]w&=&L(0)_s(L(j)w)-L(j)(L(0)_sw)\\
&=&(\wt w-j)L(j)w-(\wt w)L(j)w\\
&=&-jL(j)w\\
&=&[L(0), L(j)]w.
\end{eqnarray*}
Thus the ``locally nilpotent part'' $L(0)-L(0)_s$ of $L(0)$ commutes
with the action of $V$ and of ${\mathfrak s}{\mathfrak l}(2)$ on $W$.  In
other words, $L(0)-L(0)_s$ is a $V$-homomorphism {}from $W$ to itself.}
\end{rema}

Now suppose that $L(1)$ acts locally nilpotently on a M\"obius (or
conformal) vertex algebra $V$, that is, for
any $v\in V$, there is $m\in {\mathbb N}$ such that $L(1)^mv=0$.  Then
generalizing formula (3.20) in \cite{tensor1} (the case of ordinary
modules for a vertex operator algebra), we define the {\it opposite
vertex operator} on a generalized $V$-module $(W,Y_W)$ associated to
$v\in V$ by
\begin{equation}\label{yo}
Y^o_W(v,x)=Y_W(e^{xL(1)}(-x^{-2})^{L(0)}v,x^{-1}),
\end{equation}
that is, for $k\in {\mathbb Z}$ and $v\in V_{(k)}$,
\begin{eqnarray}\label{yo1}
Y^o_W(v,x)&=&\sum_{n \in {\mathbb Z}} v^o_n x^{-n-1}\nno\\
&=&\sum_{n\in {\mathbb Z}}\bigg((-1)^k\sum_{m\in {\mathbb N}}
\frac{1}{m!}(L(1)^mv)_{-n-m-2+2k}\bigg)x^{-n-1},
\end{eqnarray}
as in \cite{tensor1}. (In the present work, we are replacing the 
symbol $*$ used in \cite{tensor1} for opposite vertex operators
by the symbol $o$; see also Section 5.1 below.)
Here we are defining the component operators
\begin{equation}\label{v^o_n}
v^o_n=(-1)^k\sum_{m\in {\mathbb N}}\frac{1}{m!}(L(1)^mv)_{-n-m-2+2k}
\end{equation}
for $v\in V_{(k)}$ and $n,k\in {\mathbb Z}$. Note that the $L(1)$-local
nilpotence ensures well-definedness here.  Clearly, $v \mapsto
Y^o_W(v,x)$ is a linear map $V \to (\mbox{End} \ W)[[x,x^{-1}]]$ such
that $V \otimes W \to W((x^{-1}))$ ($v \otimes w \mapsto
Y^o_W(v,x)w)$.

By (\ref{v^o_n}), (\ref{degL(j)}) and (\ref{set:wtvn}), we see that
for $n,k\in {\mathbb Z}$ and $v\in V_{(k)}$, the operator $v^o_n$ is of
generalized weight $n+1-k\;(=n+1-\wt v)$, in the sense that
\begin{equation}\label{v^o-deg}
v^o_nW_{[m]}\subset W_{[m+n+1-k]}\;\;\mbox{ for any }\;m\;\in{\mathbb C}.
\end{equation}

As mentioned in \cite{tensor1} (see (3.23) in \cite{tensor1}), the
proof of the Jacobi identity in Theorem 5.2.1 of \cite{FHL} proves the
following {\it opposite Jacobi identity} for $Y^o_{W}$ in the case
where $V$ is a vertex operator algebra and $W$ is a $V$-module:
\begin{eqnarray}\label{op-jac-id}
\lefteqn{\dps x_{0}^{-1}\delta\bigg(\frac{x_{1}-x_{2}}{x_{0}}\bigg)
Y_{W}^o(v, x_{2})Y^o_{W}(u, x_{1})}\nno\\
&&\hspace{2em}-x_{0}^{-1}\delta\bigg(\frac{x_{2}-x_{1}}{-x_{0}}\bigg)
Y_{W}^o(u, x_{1})Y^o_{W}(v, x_{2})\nno\\
&&{\dps =x_{2}^{-1}\delta\bigg(\frac{x_{1}-x_{0}}{x_{2}}\bigg)
Y_{W}^o(Y(u, x_{0})v, x_{2})}
\end{eqnarray}
for $u,v \in V$, and taking ${\rm Res}_{x_{0}}$ gives us the {\it
opposite commutator formula}.  Similarly, the proof of the
$L(-1)$-derivative property in Theorem 5.2.1 of \cite{FHL} proves the
following {\it $L(-1)$-derivative property} for $Y^o_{W}$ in the same
case:
\begin{equation}\label{yo-l-1}
\frac{d}{dx} Y^o_{W}(v,x) = Y^o_{W}(L(-1)v,x).
\end{equation}
The same proofs carry over and prove the opposite Jacobi identity and
the $L(-1)$-derivative property
for $Y^o_{W}$  in the present case, where $V$ is a M\"obius (or
conformal) vertex algebra with $L(1)$ acting locally nilpotently and
where $W$ is a generalized $V$-module.  In the case in which $V$ is a
conformal vertex algebra, we have
\begin{equation}\label{Yoppositeomega}
Y^o_{W}(\omega,x) = Y_{W}(x^{-4}\omega,x^{-1})
=\sum_{n\in {\mathbb Z}}L(n)x^{n-2}
\end{equation}
since $L(1)\omega = 0$.

For opposite vertex operators, we have the following analogues of
(\ref{sl2-1})--(\ref{sl2-all}) in the M\"obius case:

\begin{lemma}\label{sl2opposite}
For $v \in V$,
\begin{equation}\label{sl2opp-1}
[Y^o_W(v,x),L(1)]=Y^o_W(L(-1)v,x),
\end{equation}
\begin{equation}\label{sl2opp-2}
[Y^o_W(v,x),L(0)]=Y^o_W(L(0)v,x)+xY^o_W(L(-1)v,x),
\end{equation}
\begin{equation}\label{sl2opp-3}
[Y^o_W(v,x),L(-1)]=Y^o_W(L(1)v,x)+2xY^o_W(L(0)v,x)+x^2Y^o_W(L(-1)v,x).
\end{equation}
Equivalently,
\begin{eqnarray}\label{sl2opp-all}
[Y^o_W(v,x),L(-j)]&=&\sum_{k=0}^{j+1}{j+1\choose k}x^k
Y^o_W(L(j-k)v,x)\nno\\
&=&\sum_{k=0}^{j+1}{j+1\choose k}x^{j+1-k}
Y^o_W(L(k-1)v,x)
\end{eqnarray}
for $j=0,\pm 1$.
\end{lemma}
\pf
For $j=0, \pm 1$, by definition and (\ref{sl2-all}) we have
\begin{eqnarray}\label{sl2opp-all-1}
[Y^o_W(v,x),L(j)]&=&-[L(j), Y_W(e^{xL(1)}(-x^{-2})^{L(0)}v,x^{-1})]\nn
&=&-\sum_{k=0}^{j+1}{j+1\choose k}x^{-k}Y_{W}(L(j-k)
e^{xL(1)}(-x^{-2})^{L(0)}v,x^{-1}).
\end{eqnarray}
By (5.2.14) in \cite{FHL} and the fact that
\begin{equation}\label{xL(0)L(j)}
x^{L(0)}L(j)x^{-L(0)}=x^{-j}L(j)
\end{equation}
(easily proved by applying to a homogeneous vector),
\begin{eqnarray}\label{sl2opp-all-2}
\lefteqn{L(-1)e^{xL(1)}(-x^{-2})^{L(0)}}\nn
&&=
e^{xL(1)}L(-1)(-x^{-2})^{L(0)}-2xe^{xL(1)}L(0)(-x^{-2})^{L(0)}
+x^{2}e^{xL(1)}L(1)(-x^{-2})^{L(0)}\nn
&&=-x^{2}e^{xL(1)}(-x^{-2})^{L(0)}L(-1)-2xe^{xL(1)}(-x^{-2})^{L(0)}L(0)
-e^{xL(1)}(-x^{-2})^{L(0)}L(1)\nn
&&=-e^{xL(1)}(-x^{-2})^{L(0)}(x^{2}L(-1)+2xL(0)+L(1)).
\end{eqnarray}
We also have
\begin{eqnarray}\label{sl2opp-all-4}
L(1)e^{xL(1)}(-x^{-2})^{L(0)}&=&
e^{xL(1)}L(1)(-x^{-2})^{L(0)}\nn
&=&-x^{-2}e^{xL(1)}(-x^{-2})^{L(0)}L(1).
\end{eqnarray}
By (\ref{sl2opp-all-2}), (\ref{sl2opp-all-4}),
$L(0)=\frac{1}{2}[L(1), L(-1)]$ and $[L(1), L(0)]=L(1)$, we have
\begin{eqnarray}\label{sl2opp-all-3}
\lefteqn{L(0)e^{xL(1)}(-x^{-2})^{L(0)}}\nn
&&=\frac{1}{2}L(1)L(-1)e^{xL(1)}(-x^{-2})^{L(0)}
-\frac{1}{2}L(-1)L(1)e^{xL(1)}(-x^{-2})^{L(0)}\nn
&&=-\frac{1}{2}L(1)e^{xL(1)}(-x^{-2})^{L(0)}(x^{2}L(-1)+2xL(0)+L(1))\nn
&&\quad +\frac{1}{2}x^{-2}L(-1)e^{xL(1)}(-x^{-2})^{L(0)}L(1)\nn
&&=\frac{1}{2}x^{-2}e^{xL(1)}(-x^{-2})^{L(0)}L(1)(x^{2}L(-1)+2xL(0)+L(1))\nn
&&\quad
-\frac{1}{2}x^{-2}e^{xL(1)}(-x^{-2})^{L(0)}(x^{2}L(-1)+2xL(0)+L(1))L(1)\nn
&&=e^{xL(1)}(-x^{-2})^{L(0)}L(0)+x^{-1}e^{xL(1)}(-x^{-2})^{L(0)}L(1)\nn
&&=e^{xL(1)}(-x^{-2})^{L(0)}(L(0)+x^{-1}L(1)).
\end{eqnarray}
Thus we obtain
\begin{eqnarray*}
\lefteqn{[Y^o_W(v,x),L(1)]}\nn
&&=-\sum_{k=0}^{2}{2\choose k}x^{-k}Y_{W}(L(1-k)
e^{xL(1)}(-x^{-2})^{L(0)}v,x^{-1})\nn
&&=-Y_{W}(L(1)e^{xL(1)}(-x^{-2})^{L(0)}v,x^{-1})
-2x^{-1}Y_{W}(L(0)
e^{xL(1)}(-x^{-2})^{L(0)}v,x^{-1})\nn
&&\quad -x^{-2}Y_{W}(L(-1)
e^{xL(1)}(-x^{-2})^{L(0)}v,x^{-1})\nn
&&=x^{-2}Y_{W}(e^{xL(1)}(-x^{-2})^{L(0)}L(1)v,x^{-1})
\nn
&&\quad -2x^{-1}Y_{W}(
e^{xL(1)}(-x^{-2})^{L(0)}(L(0)+x^{-1}L(1))v,x^{-1})
\nn
&&\quad +x^{-2}Y_{W}(
e^{xL(1)}(-x^{-2})^{L(0)}(x^{2}L(-1)+2xL(0)+L(1))v,x^{-1})\nn
&&=Y_{W}(
e^{xL(1)}(-x^{-2})^{L(0)}L(-1)v,x^{-1})\nn
&&=Y_{W}^{o}(L(-1)v,x),
\end{eqnarray*}
\begin{eqnarray*}
\lefteqn{[Y^o_W(v,x),L(0)]}\nn
&&=-\sum_{k=0}^{1}{1\choose k}x^{-k}Y_{W}(L(-k)
e^{xL(1)}(-x^{-2})^{L(0)}v,x^{-1})\nn
&&=-Y_{W}(L(0)
e^{xL(1)}(-x^{-2})^{L(0)}v,x^{-1})
-x^{-1}Y_{W}(L(-1)
e^{xL(1)}(-x^{-2})^{L(0)}v,x^{-1})\nn
&&=-Y_{W}(e^{xL(1)}(-x^{-2})^{L(0)}(L(0)+x^{-1}L(1))v,x^{-1})\nn
&&\quad
+x^{-1}Y_{W}(
e^{xL(1)}(-x^{-2})^{L(0)}(x^{2}L(-1)+2xL(0)+L(1))v,x^{-1})\nn
&&=Y_{W}(
e^{xL(1)}(-x^{-2})^{L(0)}(xL(-1)+L(0))v,x^{-1})\nn
&&=Y_{W}^{o}(L(0)v,x)+xY_{W}^{o}(L(-1)v,x)
\end{eqnarray*}
and
\begin{eqnarray*}
[Y^o_W(v,x),L(-1)]&=&-Y_{W}(L(-1)
e^{xL(1)}(-x^{-2})^{L(0)}v,x^{-1})\nn
&=&Y_{W}(
e^{xL(1)}(-x^{-2})^{L(0)}(x^{2}L(-1)+2xL(0)+L(1))v,x^{-1})\nn
&=&Y_{W}^{o}(L(1)v,x)+2xY_{W}^{o}(L(0)v,x)+x^{2}Y_{W}^{o}(L(-1)v,x),
\end{eqnarray*}
proving the lemma.
\epfv

As in Section 5.2 of \cite{FHL}, we can define a $V$-action on $W^*$
as follows:
\begin{equation}\label{y'}
\langle Y'(v,x)w',w\rangle = \langle w', Y^o_W(v,x)w\rangle
\end{equation}
for $v\in V$, $w'\in W^*$ and $w\in W$; the correspondence $v\mapsto
Y'(v,x)$ is a linear map {}from $V$ to $({\rm End}\,W^*)[[x,x^{-1}]]$.
Writing
\[
Y'(v,x)=\sum_{n\in {\mathbb Z}} v_n x^{-n-1}
\]
($v_n\in {\rm End}\,W^*)$, we have
\begin{equation}\label{v'vo}
\langle v_n w', w\rangle = \langle w', v^o_n w\rangle
\end{equation}
for $v\in V$, $w'\in W^*$ and $w\in W$.  (Actually, in \cite{FHL} this
$V$-action was defined on a space smaller than $W^*$, but this
definition holds without change on all of $W^*$.)  In the case in
which $V$ is a conformal vertex algebra we define the operators
$L'(n)\;(n\in {\mathbb Z})$ by
\[
Y'(\omega,x)=\sum_{n\in {\mathbb Z}}L'(n)x^{-n-2};
\]
then, by extracting the coefficient of $x^{-n-2}$ in (\ref{y'}) with
$v=\omega$ and using the fact that $L(1)\omega=0$ we have
\begin{equation}\label{L'(n)}
\langle L'(n)w',w\rangle=\langle w',L(-n)w\rangle\;\;\mbox{ for }\;n\in
{\mathbb Z}
\end{equation}
(see (\ref{Yoppositeomega})), as in Section 5.2 of \cite{FHL}.
In the case where $V$ is only a M\"obius vertex algebra, we define operators
$L'(-1)$, $L'(0)$ and $L'(1)$ on $W^*$ by formula (\ref{L'(n)}) for $n=0,
\pm 1$. It follows {}from (\ref{set:wtsl2}) that
\begin{equation}\label{L'(n)2}
L'(j)(W_{[m]})^*\subset (W_{[m-j]})^*
\end{equation}
for $m\in {\mathbb C}$ and $j=0,\pm 1$.  By combining (\ref{v'vo}) with
(\ref{v^o-deg}) we get
\begin{equation}\label{stable0}
v_n(W_{[m]})^*\subset (W_{[m+k-n-1]})^*
\end{equation}
for any $n,k\in {\mathbb Z}$, $v\in V_{(k)}$ and $m\in {\mathbb C}$.

We have just seen that the $L(1)$-local nilpotence
condition enables us to define a natural vertex operator action on the
vector space dual of a generalized module for a M\"obius (or
conformal) vertex algebra. This condition is satisfied by all
vertex operator algebras, due to (\ref{degL(j)}) and the grading
restriction condition (\ref{gr1}). However, the functor $W\mapsto W^*$
is certainly
not involutive, and $W^*$ is not in general a generalized module. In this work
we will need certain module categories equipped with an involutive
``contragredient functor'' $W\mapsto W'$ which generalizes the
contragredient functor for the category of modules for vertex operator
algebras. For this purpose, we introduce the following:

\begin{defi}\label{def:dgv}
{\rm Let $A$ be an abelian group.  A M\"obius (or conformal) vertex
algebra
\[
V=\coprod_{n\in {\mathbb Z}} V_{(n)}
\]
is said to be {\em strongly graded with respect to $A$} (or {\em
strongly $A$-graded}, or just {\em strongly graded} if the abelian
group $A$ is understood) if $V$ is equipped with a second gradation,
by $A$,
\[
V=\coprod _{\alpha \in A} V^{(\alpha)},
\]
such that the following conditions are satisfied: the two gradations
are compatible, that is,
\[
V^{(\alpha)}=\coprod_{n\in {\mathbb Z}} V^{(\alpha)}_{(n)} \;\;
(\mbox{where}\;V^{(\alpha)}_{(n)}=V_{(n)}\cap V^{(\alpha)})\;
\mbox{ for any }\;\alpha \in A;
\]
for any $\alpha,\beta\in A$ and $n\in {\mathbb Z}$,
\begin{eqnarray}
&V^{(\alpha)}_{(n)}=0\;\;\mbox{ for }\;n\;\mbox{ sufficiently
negative};&\label{dua:ltc}\\
&\dim V^{(\alpha)}_{(n)} <\infty;&\label{dua:fin}\\
&{\bf 1}\in V^{(0)}_{(0)};&\\
&v_l V^{(\beta)} \subset V^{(\alpha+\beta)}\;\;
\mbox{ for any }\;v\in V^{(\alpha)},\;l\in {\mathbb Z};&\label{v_l-A}
\end{eqnarray}
and
\begin{equation}\label{L(n)-A}
L(j)V^{(\alpha)} \subset V^{(\alpha)}\;\;\mbox{ for }\;j=0,\pm 1.
\end{equation}
If $V$ is in fact a conformal vertex algebra, we in addition require
that
\begin{equation}\label{omega0}
\omega\in V^{(0)}_{(2)},
\end{equation}
so that for all $j\in {\mathbb Z}$, (\ref{L(n)-A}) follows {}from
(\ref{v_l-A}).  }
\end{defi}

\begin{rema}\label{rm1}{\rm
Note that the notion of conformal vertex algebra strongly
graded with respect to the trivial group is exactly the notion of
vertex operator algebra. Also note that (\ref{degL(j)}),
(\ref{dua:ltc}) and (\ref{L(n)-A}) imply the local nilpotence of
$L(1)$ acting on $V$, and hence we have the construction and
properties of opposite
vertex operators on a generalized module for a strongly graded
M\"obius (or conformal) vertex algebra. }
\end{rema}

For (generalized) modules for a strongly graded algebra we will also
have a second grading by an abelian group, and it is natural to allow
this group to be larger than the second grading group $A$ for the
algebra.  (Note that this already occurs for the {\em first} grading
group, which is ${\mathbb Z}$ for algebras and ${\mathbb C}$ for
(generalized) modules.)  We now define the notions of strongly graded
module and generalized module, and also, at the end of this
definition, the notions of lower bounded such structures.

\begin{defi}\label{def:dgw}{\rm
Let $A$ be an abelian group and $V$ a strongly $A$-graded M\"obius (or
conformal)
vertex algebra. Let $\tilde A$ be an abelian group containing $A$ as a
subgroup. A $V$-module (respectively, generalized $V$-module)
\[
W=\coprod_{n\in{\mathbb C}} W_{(n)} \;\;\;(\mbox{respectively, }\;
W=\coprod_{n\in{\mathbb C}} W_{[n]})
\]
is said to be {\em strongly graded with respect to $\tilde A$} (or
{\em strongly $\tilde A$-graded}, or just {\em strongly graded}) if
the abelian group $\tilde A$ is understood) if $W$ is equipped with a
second gradation, by $\tilde A$,
\begin{equation}\label{2ndgrd}
W=\coprod _{\beta \in \tilde A} W^{(\beta)},
\end{equation}
such that the following conditions are satisfied: the two gradations
are compatible, that is, for any $\beta \in \tilde A$,
\[
W^{(\beta)}=\coprod_{n\in {\mathbb C}} W^{(\beta)}_{(n)} \;\;(\mbox{where }\;
W^{(\beta)}_{(n)}=W_{(n)}\cap W^{(\beta)})
\]
\[
(\mbox{respectively, }\;
W^{(\beta)}=\coprod_{n\in {\mathbb C}} W^{(\beta)}_{[n]} \;\;(\mbox{where }\;
W^{(\beta)}_{[n]}=W_{[n]}\cap W^{(\beta)}));
\]
for any $\alpha\in A$, $\beta\in \tilde A$ and $n\in {\mathbb C}$,
\begin{eqnarray}
&W^{(\beta)}_{(n+k)}=0 \;\; (\mbox{respectively, } \; W^{(\beta)}_{[n+k]}=0)
\;\;\mbox{ for }\;k\in {\mathbb Z}\;\mbox{ sufficiently
negative};&\label{set:dmltc}\\
&\dim W^{(\beta)}_{(n)} <\infty \;\; (\mbox{respectively, } \;
\dim W^{(\beta)}_{[n]} <\infty);&\label{set:dmfin}\\
&v_l W^{(\beta)} \subset W^{(\alpha+\beta)}\;\;\mbox{ for any }\;v\in
V^{(\alpha)},\;l\in {\mathbb Z};&\label{m-v_l-A}
\end{eqnarray}
and
\begin{equation}\label{m-L(n)-A}
L(j)W^{(\beta)} \subset W^{(\beta)}\;\;\mbox{ for }\;j=0,\pm 1.
\end{equation}
(Note that if $V$ is in fact a conformal vertex algebra, then for all
$j\in {\mathbb Z}$, (\ref{m-L(n)-A}) follows {}from (\ref{omega0}) and
(\ref{m-v_l-A}).)  A strongly $\tilde{A}$-graded (generalized) $V$-module $W$
is said to be {\it lower bounded} if instead of (\ref{set:dmltc}),
it satisfies the stronger condition that for any $\beta \in \tilde{A}$,
\begin{equation}\label{set:dmltc-1}
W^{(\beta)}_{(n)}=0 \;\; (\mbox{respectively, } \; W^{(\beta)}_{[n]}=0)
\;\;\mbox{ for }\;\Re{(n)}\;\mbox{ sufficiently
negative}
\end{equation}
($n \in \C$).}
\end{defi}

\begin{rema}\label{v-str-module}{\rm
A strongly  $A$-graded conformal or M\"{o}bius vertex algebra is a 
strongly  $A$-graded module for itself (and in particular, a 
strongly  $A$-graded generalized module for itself), and is in fact
lower bounded.}
\end{rema}

\begin{rema}\label{moduleswiththetrivialgroup}{\rm
Let $V$ be a vertex operator algebra, viewed (equivalently) as a
conformal vertex algebra strongly graded with respect to the trivial
group (recall Remark \ref{rm1}).  Then the $V$-modules that are
strongly graded with respect to the trivial group (in the sense of
Definition \ref{def:dgw}) are exactly the ($\C$-graded) modules for
$V$ as a vertex operator algebra, with the grading restrictions as
follows: For $n \in \C$,
\begin{equation}\label{Wn+k=0}
W_{(n+k)}=0 \;\; \mbox { for }\;k\in {\mathbb Z}\;\mbox{ sufficiently
negative}
\end{equation}
and
\begin{equation}\label{dimWnfinite}
\dim W_{(n)} <\infty,
\end{equation}
and the lower bounded such structures have (\ref{Wn+k=0}) replaced by:
\begin{equation}\label{ReWn=0}
W_{(n)}=0 \;\; \mbox { for }\;\Re{(n)}\;\mbox{ sufficiently
negative}.
\end{equation}
Also, the generalized $V$-modules that are strongly graded with
respect to the trivial group are exactly the generalized $V$-modules
(in the sense of Definition \ref{definitionofgeneralizedmodule}) such
that for $n \in \C$,
\begin{equation}\label{W[n+k]=0}
W_{[n+k]}=0 \;\; \mbox { for }\;k\in {\mathbb Z}\;\mbox{ sufficiently
negative}
\end{equation}
and
\begin{equation}
\dim W_{[n]} <\infty,
\end{equation}
and the lower bounded ones have (\ref{Wn+k=0}) replaced by:
\begin{equation}\label{ReW[n]=0}
W_{[n]}=0 \;\; \mbox { for }\;\Re{(n)}\;\mbox{ sufficiently
negative}.
\end{equation}}
\end{rema}

\begin{rema}\label{homsaregradingpreserving}{\rm
In the strongly graded case, algebra and module homomorphisms are of
course understood to preserve the grading by $A$ or $\tilde A$.}
\end{rema}

\begin{exam}{\rm
An important source of examples of strongly graded conformal vertex
algebras and modules comes {}from the vertex algebras and modules
associated with even lattices.  Let $L$ be an even lattice, i.e., a
finite-rank free abelian group equipped with a nondegenerate symmetric
bilinear form $\langle\cdot,\cdot\rangle$, not necessarily positive
definite, such that $\langle \alpha,\alpha\rangle\in 2{\mathbb Z}$ for
all $\alpha\in L$.  Then there is a natural structure of conformal vertex algebra on
a certain vector space $V_L$; see \cite{B} and Chapter 8 of
\cite{FLM2}.  If the form $\langle\cdot,\cdot\rangle$ on $L$ is also
positive definite, then $V_L$ is a vertex operator algebra (that is,
the grading restrictions hold).  If $L$ is not necessarily positive
definite, then $V_L$ is equipped with a natural second grading given
by $L$ itself, making $V_L$ a strongly $L$-graded conformal vertex
algebra in the sense of Definition \ref{def:dgv}.  Any (rational)
sublattice $M$ of the ``dual lattice'' $L^{\circ}$ of $L$ containing
$L$ gives rise to a lower bounded strongly $M$-graded module for the strongly
$L$-graded conformal vertex algebra (see Chapter 8 of \cite{FLM2};
cf.\ \cite{LL}).}
\end{exam}

In the next two remarks, we mention certain important properties of
compositions of two or more vertex operators, properties that will
also be important in the further generality of logarithmic
intertwining operators in the future.

\begin{rema}{\rm
As mentioned in Remark \ref{rm1}, strong gradedness for a M\"obius (or
conformal) vertex algebra $V$ implies the
local nilpotence of $L(1)$ acting on $V$. In fact,
strong gradedness implies much more that will be important for us:
{}From (\ref{dua:ltc}),
(\ref{dua:fin}), (\ref{v_l-A}) and (\ref{L(n)-A}) (and (\ref{omega0})
in the conformal vertex algebra case), it is clear that strong gradedness
for $V$ implies the following {\em local grading restriction
condition on} $V$ (see \cite{H-codes}):
\begin{enumerate}
\item[(i)] for any $m>0$ and $v_{(1)}, \dots, v_{(m)}\in V$, there
exists $r\in {\mathbb Z}$ such that the coefficient of each monomial in
$x_1, \dots, x_{m-1}$ in the formal series
\[
Y(v_{(1)}, x_1)\cdots Y(v_{(m-1)}, x_{m-1})v_{(m)}
\]
lies in $\coprod_{n>r}V_{(n)}$;
\item[(ii)] in the conformal vertex algebra case: for any element of
the conformal vertex algebra $V$ homogeneous with respect to the
weight grading, the Virasoro-algebra submodule $M=\coprod_{n\in {\mathbb
Z}}M_{(n)}$ (where $M_{(n)}=M\cap V_{(n)}$) of $V$ generated by this
element satisfies the following grading restriction conditions:
$M_{(n)}=0$ when $n$ is sufficiently negative and $\dim
M_{(n)}<\infty$ for $n\in {\mathbb Z}$
\end{enumerate}
or
\begin{enumerate}
\item[(ii$'$)] in the M\"obius vertex algebra case: for any element of
the M\"obius vertex algebra $V$ homogeneous with respect to the weight
grading, the ${\mathfrak s}{\mathfrak l}(2)$-submodule $M=\coprod_{n\in {\mathbb
Z}}M_{(n)}$ (where $M_{(n)}=M\cap V_{(n)}$) of $V$ generated by this
element satisfies the following grading restriction conditions:
$M_{(n)}=0$ when $n$ is sufficiently negative and $\dim
M_{(n)}<\infty$ for $n\in {\mathbb Z}$.
\end{enumerate}
As was pointed out in \cite{H-codes}, Condition (i) above was first stated
in \cite{DL} (see formula (9.39), Proposition 9.17 and Theorem 12.33
in \cite{DL}) for generalized vertex algebras and abelian intertwining
algebras (certain generalizations of vertex algebras); it guarantees
the convergence, rationality and commutativity properties of the
matrix coefficients of products of more than two vertex operators.
Conditions (i) and (ii) (or (ii$'$)) together ensure that all the
essential results involving the Virasoro operators and the geometry of
vertex operator algebras in \cite{H0} and \cite{H1} still hold for these algebras.}
\end{rema}

\begin{rema}{\rm
Similarly, {}from (\ref{set:dmltc}), (\ref{set:dmfin}),
(\ref{m-v_l-A}) and (\ref{m-L(n)-A}) (and (\ref{omega0}) in the
conformal vertex algebra case), it is clear that strong gradedness for
(generalized) modules implies the following {\it local grading
restriction condition on a (generalized) module} $W$ for a strongly
graded M\"obius (or conformal) vertex algebra $V$:
\begin{enumerate}
\item[(i)] for any $m>0$, $v_{(1)}, \dots, v_{(m-1)}\in V$, $n\in{\mathbb
C}$ and $w\in W_{[n]}$, there exists $r\in {\mathbb Z}$ such that the
coefficient of each monomial in $x_1, \dots, x_{m-1}$ in the formal
series
\[
Y(v_{(1)}, x_1)\cdots Y(v_{(m-1)}, x_{m-1})w
\]
lies in $\coprod_{k>r}W_{[n+k]}$;
\item[(ii)] in the conformal vertex algebra case: for any $w\in
W_{[n]}$ ($n\in{\mathbb C}$), the Virasoro-algebra submodule
$M=\coprod_{k\in {\mathbb Z}}M_{[n+k]}$ (where $M_{[n+k]}=M\cap
W_{[n+k]}$) of $W$ generated by $w$ satisfies the following grading
restriction conditions: $M_{[n+k]}=0$ when $k$ is sufficiently
negative and $\dim M_{[n+k]}<\infty$ for $k\in {\mathbb Z}$
\end{enumerate}
or
\begin{enumerate}
\item[(ii$'$)] in the M\"obius vertex algebra case: for any $w\in
W_{[n]}$ ($n\in {\mathbb C}$), the ${\mathfrak s}{\mathfrak l}(2)$-submodule
$M=\coprod_{k\in {\mathbb Z}}M_{[n+k]}$ (where $M_{[n+k]}=M\cap
W_{[n+k]}$) of $W$ generated by $w$ satisfies the following grading
restriction conditions: $M_{[n+k]}=0$ when $k$ is sufficiently
negative and $\dim M_{[n+k]}<\infty$ for $k\in {\mathbb Z}$.
\end{enumerate}
Note that in the case of ordinary (as opposed to generalized) modules,
all the generalized weight spaces such as $W_{[n]}$ mentioned here are
ordinary weight spaces $W_{(n)}$.  Analogous statments of course hold
for lower bounded (generalized) modules.}
\end{rema}

With the strong gradedness condition on a (generalized) module, we can
now define the corresponding notion of contragredient module. First we
give:

\begin{defi}\label{defofWprime}{\rm
Let $W=\coprod_{\beta \in \tilde A,\;n\in{\mathbb C}}W^{(\beta)}_{[n]}$
be a strongly $\tilde A$-graded generalized module for a strongly
$A$-graded M\"obius (or conformal) vertex algebra.  For each $\beta\in
\tilde A$ and $n \in {\mathbb C}$, let us identify
$(W^{(\beta)}_{[n]})^*$ with the subspace of $W^*$ consisting of the
linear functionals on $W$ vanishing on each $W^{(\gamma)}_{[m]}$ with
$\gamma \neq \beta$ or $m \neq n$ (cf.\ (\ref{Wnstar})).  We define
$W'$ to be the $(\tilde A\times{\mathbb C})$-graded vector subspace of
$W^*$ given by
\begin{equation}
W'=\coprod_{\beta \in \tilde A,\; n\in{\mathbb C}}
(W')^{(\beta)}_{[n]}, \;\;\mbox{ where }\;
(W')^{(\beta)}_{[n]}=(W^{(-\beta)}_{[n]})^*;
\end{equation}
we also use the notations
\begin{equation}\label{W'beta}
(W')^{(\beta)}=\coprod_{n\in{\mathbb C}}(W^{(-\beta)}_{[n]})^*
\subset(W^{(-\beta)})^* \subset W^*
\end{equation}
(where $(W^{(\beta)})^*$ consists of the linear functionals on $W$
vanishing on all $W^{(\gamma)}$ with $\gamma \neq \beta$) and
\begin{equation}
(W')_{[n]}=\coprod_{\beta\in\tilde A}(W^{(-\beta)}_{[n]})^*
\subset(W_{[n]})^* \subset W^*
\end{equation}
for the homogeneous subspaces of $W'$ with respect to the $\tilde A$-
and ${\mathbb C}$-grading, respectively (The reason for the minus
signs here will become clear below.)  We will still use the notation
$\langle\cdot,\cdot\rangle_W$, or $\langle\cdot,\cdot\rangle$ when the
underlying space is clear, for the canonical pairing between $W'$ and
\[
\overline{W}\subset \prod_{\beta \in \tilde A,\; n\in{\mathbb
C}}W^{(\beta)}_{[n]}
\]
(recall (\ref{Wbar})).  }
\end{defi}

\begin{rema}{\rm
In the case of ordinary rather than generalized modules, Definition
\ref{defofWprime} still applies, and all of the generalized weight
subspaces $W_{[n]}$ of $W$ are ordinary weight spaces $W_{(n)}$.  In
this case, we can write $(W')_{(n)}$ rather than $(W')_{[n]}$ for the
corresponding subspace of $W'$.}
\end{rema}

Let $W$ be a strongly graded (generalized) module for a strongly
graded M\"obius (or conformal) vertex algebra $V$. Recall that we have
the action (\ref{y'}) of $V$ on $W^*$ and that (\ref{stable0}) holds.
Furthermore, (\ref{v^o_n}), (\ref{v'vo}) and (\ref{m-v_l-A}) imply
for any $n,k\in {\mathbb Z}$, $\alpha\in A$, $\beta\in \tilde A$, $v\in
V^{(\alpha)}_{(k)}$ and $m\in {\mathbb C}$,
\begin{equation}\label{shift}
v_n((W')^{(\beta)}_{[m]})=v_n((W^{(-\beta)}_{[m]})^*)\subset
(W^{(-\alpha-\beta)}_{[m+k-n-1]})^* =(W')^{(\alpha+\beta)}_{[m+k-n-1]}.
\end{equation}
Thus $v_n$ preserves $W'$ for $v \in V$, $n \in {\mathbb Z}$.  Similarly
(in the M\"obius case), (\ref{L'(n)}), (\ref{L'(n)2}) and
(\ref{m-L(n)-A}) imply that $W'$ is stable under the operators
$L'(-1)$, $L'(0)$ and $L'(1)$, and in fact
\[
L'(j)(W')^{(\beta)}_{[n]}\subset (W')^{(\beta)}_{[n-j]}
\]
for any $j=0,\pm 1$, $\beta\in \tilde A$ and $n\in {\mathbb C}$.  In
the case or ordinary rather than generalized modules, the symbols
$(W')^{(\beta)}_{[n]}$, etc., can be replaced by
$(W')^{(\beta)}_{(n)}$, etc.

For any fixed $\beta\in \tilde A$ and $n\in {\mathbb C}$, by
(\ref{gerwt}) and the finite-dimensionality (\ref{set:dmfin}) of
$W^{(-\beta)}_{[n]}$, there exists $N\in {\mathbb N}$
such that $(L(0)-n)^NW^{(-\beta)}_{[n]}=0$. But then for any $w'\in
(W')^{(\beta)}_{[n]}$,
\begin{equation}\label{L(0)N}
\langle (L'(0)-n)^Nw',w\rangle=\langle w',(L(0)-n)^Nw\rangle=0
\end{equation}
for all $w\in W$. Thus $(L'(0)-n)^Nw'=0$. So (\ref{gerwt}) holds with
$W$ replaced by $W'$.  In the case of ordinary modules, we of course
take $N=1$.

By (\ref{set:dmltc}) and (\ref{shift}) we have the lower truncation
condition for the action $Y'$ of $V$ on $W'$:
\begin{equation}\label{truncationforY'}
\mbox{For any }\;v\in V\;\mbox{ and }\;w'\in W',\; v_nw'=0\;\;
\mbox{ for }\;n\;\mbox{ sufficiently large}.
\end{equation}
As a consequence, the Jacobi identity can now be formulated on
$W'$. In fact, by the above, and using the same proofs as those of
Theorems 5.2.1 and 5.3.1 in \cite{FHL}, together with Lemma
\ref{sl2opposite}, we obtain:

\begin{theo}\label{set:W'}
Let $\tilde A$ be an abelian group containing $A$ as a subgroup and
$V$ a strongly $A$-graded M\"obius (or conformal) vertex algebra. Let
$(W,Y)$ be a strongly $\tilde A$-graded $V$-module (respectively,
generalized $V$-module). Then the pair $(W',Y')$ carries a strongly
$\tilde A$-graded $V$-module (respectively, generalized $V$-module)
structure, and
\[
(W'',Y'')=(W,Y).
\]
If $W$ is lower bounded, then so is $W'$.  \epf
\end{theo}

\begin{defi}{\rm
The pair $(W',Y')$ in Theorem \ref{set:W'} will be called the {\em
contragredient module} of $(W,Y)$.}
\end{defi}

Let $W_1$ and $W_2$ be strongly $\tilde A$-graded (generalized)
$V$-modules and let $f:W_1\to W_2$ be a module homomorphism (which is
of course understood to preserve both the ${\mathbb C}$-grading and the
$\tilde A$-grading, and to preserve the action of ${\mathfrak s}{\mathfrak
l}(2)$ in the M\"obius case). Then by (\ref{v'vo}) and (\ref{L'(n)}),
the linear map
\[
f':W'_2\to W'_1
\]
given by
\begin{equation}\label{fprime}
\langle f'(w'_{(2)}), w_{(1)}\rangle=\langle
w'_{(2)},f(w_{(1)})\rangle
\end{equation}
for any $w_{(1)}\in W_1$ and $w'_{(2)}\in W'_2$ is well defined and is
clearly a module homomorphism {}from $W'_2$ to $W'_1$.

\begin{nota}\label{MGM}
{\rm In this work we will be especially interested in the case where
$V$ is strongly $A$-graded, and we will be focusing on the category of all
strongly $\tilde{A}$-graded (ordinary) $V$-modules, for which we will use the
notation
\[
{\cal M}_{sg},
\]
or the category of all strongly $\tilde{A}$-graded generalized $V$-modules,
which we will call
\[
{\cal GM}_{sg}.
\]  
{}From the above we see that in
the strongly graded case we have contravariant functors
\[
(\cdot)': (W,Y)\mapsto (W',Y'),
\]
the {\it contragredient functors}, {}from ${\cal M}_{sg}$ to itself and
{}from ${\cal GM}_{sg}$ to itself, and also from the full
subcategories of lower bounded such structures to themselves.
We also know that $V$ itself is a (lower bounded)
object of ${\cal M}_{sg}$ (and thus of ${\cal GM}_{sg}$ as well);
recall Remark \ref{v-str-module}. Our main objects of study will be
certain full subcategories ${\cal C}$ of ${\cal M}_{sg}$ or ${\cal
GM}_{sg}$ that are closed under the contragredient functor and such
that $V\in \ob {\cal C}$.}
\end{nota}

\begin{rema}{\rm
In order to formulate certain results in this work, even in the case
when our M\"obius or conformal vertex algebra $V$ is strongly graded
we will in fact sometimes use the category whose objects are {\it all}
the modules for $V$ and whose morphisms are all the $V$-module
homomorphisms, and also the category of {\it all} the generalized
modules for $V$.  (If $V$ is conformal, then the category of all the
$V$-modules is the same whether $V$ is viewed as either conformal or
M\"obius, by Remark \ref{modulesaremodules}, and similarly for the
category of all the generalized $V$-modules.)  Note that in view of
Remark \ref{homsaregradingpreserving}, the categories ${\cal M}_{sg}$
and ${\cal GM}_{sg}$ are not full subcategories of these categories of
{\it all} modules and generalized modules.}
\end{rema}

We now recall {}from \cite{FLM2}, \cite{FHL}, \cite{DL} and \cite{LL}
the well-known principles that vertex operator algebras (which are
exactly conformal vertex algebras strongly graded with respect to the
trivial group; recall Remark \ref{rm1}) and their modules have
important ``rationality,'' ``commutativity'' and ``associativity''
properties, and that these properties can in fact be used as axioms
replacing the Jacobi identity in the definition of the notion of
vertex operator algebra.  (These principles in fact generalize to all
vertex algebras, as in \cite{LL}.)

In the propositions below,
\[
{\C}[x_{1}, x_{2}]_{S}
\]
is the ring of formal rational functions obtained by inverting (localizing with
respect to) the products of (zero or more) elements of the set $S$ of
nonzero homogeneous linear polynomials in $x_{1}$ and $x_{2}$. Also,
$\iota_{12}$ (which might also be written as $\iota_{x_{1}x_{2}}$) is
the operation of expanding an element of ${\C}[x_{1}, x_{2}]_{S}$,
that is, a polynomial in $x_{1}$ and $x_{2}$ divided by a product of
homogeneous linear polynomials in $x_{1}$ and $x_{2}$, as a formal
series containing at most finitely many negative powers of $x_{2}$
(using binomial expansions for negative powers of linear polynomials
involving both $x_{1}$ and $x_{2}$); similarly for $\iota_{21}$ and so
on. (The distinction between rational functions and formal Laurent
series is crucial.)

Let $V$ be a vertex operator algebra.  For $W$ a ($\C$-graded)
$V$-module (including possibly $V$ itself), the space $W'$ is just the
``restricted dual space''
\begin{equation}
W'=\coprod_{n\in{\mathbb C}} W_{(n)}^{*}.
\end{equation}

\begin{propo}\label{rationalityandcommutativity}
We have:
\begin{enumerate}
\item[(a)] (rationality of products) For $v$, $v_{1}$, $v_{2}\in V$
and $v'\in V'$, the formal series
\begin{equation}\label{v'Yv1v2v}
\left\langle v', Y(v_{1}, x_{1})Y(v_{2}, x_{2})v\right\rangle,
\end{equation}
which involves only finitely
many negative powers of $x_{2}$ and only finitely many positive powers
of $x_{1}$, lies in the image of the map $\iota_{12}$:
\begin{equation}\label{Yv1v2}
\left\langle v', Y(v_{1}, x_{1})Y(v_{2}, x_{2})v\right\rangle
=\iota_{12}f(x_{1}, x_{2}),
\end{equation}
where the (uniquely determined) element $f\in {\C}[x_{1},
x_{2}]_{S}$ is of the form
\begin{equation}
f(x_{1}, x_{2})={\displaystyle \frac{g(x_{1},
x_{2})}{x_{1}^{r}x_{2}^{s}(x_{1}-x_{2})^{t}}}
\end{equation}
for some $g\in {\C}[x_{1}, x_{2}]$ and $r, s, t\in {\Z}$.
\item[(b)] (commutativity) We also have 
\begin{equation}\label{Yv2v1}
\left\langle v', Y(v_{2}, x_{2})Y(v_{1}, x_{1})v\right\rangle
=\iota_{21}f(x_{1}, x_{2}).
\end{equation}
\end{enumerate}
\end{propo}

\begin{propo}\label{rationalityofiterates}
We have:
\begin{enumerate}
\item[(a)] (rationality of iterates) For $v$, $v_{1}$, $v_{2}\in V$
and $v'\in V'$, the formal series
\begin{equation}\label{v'YYv1v2v}
\left\langle v', Y(Y(v_{1}, x_{0})v_{2}, x_{2})v\right\rangle,
\end{equation}
which involves only finitely many negative powers of $x_{0}$ and only
finitely many positive powers of $x_{2}$, lies in the image of the map
$\iota_{20}$:
\begin{equation}
\left\langle v', Y(Y(v_{1}, x_{0})v_{2},
x_{2})v\right\rangle=\iota_{20}h(x_{0}, x_{2}),
\end{equation}
where the (uniquely determined) element $h\in {\C}[x_{0},
x_{2}]_{S}$ is of the form
\begin{equation}
h(x_{0}, x_{2})={\displaystyle \frac{k(x_{0},
x_{2})}{x_{0}^{r}x_{2}^{s}(x_{0}+x_{2})^{t}}}
\end{equation}
for some $k\in {\C}[x_{0}, x_{2}]$ and $r, s, t\in {\Z}$.
\item[(b)] The formal series $\left\langle v', Y(v_{1},
x_{0}+x_{2})Y(v_{2}, x_{2})v\right\rangle,$ which involves only
finitely many negative powers of $x_{2}$ and only finitely many
positive powers of $x_{0}$, lies in the image of $\iota_{02}$, and in
fact
\begin{equation}
\left\langle v', Y(v_{1}, x_{0}+x_{2})Y(v_{2},
x_{2})v\right\rangle=\iota_{02}h(x_{0}, x_{2}).
\end{equation}
\end{enumerate}
\end{propo}

\begin{propo}\label{associativity} 
(associativity) We have the following equality of formal rational
functions:
\begin{equation}
\iota_{12}^{-1}\left\langle v', Y(v_{1}, x_{1})Y(v_{2}, x_{2})v\right\rangle
=(\iota_{20}^{-1}\left\langle v', Y(Y(v_{1}, x_{0})v_{2},
x_{2})v\right\rangle)\lbar_{x_{0}=x_{1}-x_{2}},
\end{equation}
that is,
\[
f(x_1,x_2)=h(x_1-x_2,x_2).
\]
\end{propo}

\begin{propo}\label{commandassocequivtoJacobi} 
In the presence of the other axioms for the notion of vertex operator
algebra, the Jacobi identity follows {from} the rationality of
products and iterates, and commutativity and associativity.  In
particular, in the definition of vertex operator algebra, the Jacobi
identity may be replaced by these properties.
\end{propo}

The rationality, commutativity and associativity properties
immediately imply the following result, in which the formal variables
$x_1$ and $x_2$ are specialized to nonzero complex numbers in suitable
domains:

\begin{corol}\label{dualitywithcovergence}
The formal series obtained by specializing $x_1$ and $x_2$ to
(nonzero) complex numbers $z_1$ and $z_2$, respectively, in
(\ref{v'Yv1v2v}) converges to a rational function of $z_1$ and $z_2$
in the domain
\begin{equation}
|z_1| > |z_2| > 0       
\end{equation}
and the analogous formal series obtained by specializing $x_1$ and
$x_2$ to $z_1$ and $z_2$, respectively, in (\ref{Yv2v1}) converges to
the same rational function of $z_1$ and $z_2$ in the (disjoint) domain
\begin{equation}
|z_2| > |z_1| > 0. 
\end{equation}
Moreover, the formal series obtained by specializing $x_0$ and $x_2$
to $z_1-z_2$ and $z_2$, respectively, in (\ref{v'YYv1v2v}) converges
to this same rational function of $z_1$ and $z_2$ in the domain
\begin{equation}
|z_2| > |z_1-z_2| > 0. 
\end{equation}
In particular, in the common domain
\begin{equation}
|z_1| > |z_2| > |z_1-z_2| > 0,
\end{equation}
we have the equality
\begin{equation}\label{associativitywithz1,z2}
\left\langle v', Y(v_{1}, z_{1})Y(v_{2}, z_{2})v\right\rangle=
\left\langle v', Y(Y(v_{1}, z_1-z_2)v_{2},z_{2})v\right\rangle
\end{equation}
of rational functions of $z_1$ and $z_2$.
\end{corol}

\begin{rema}{\rm
These last five results also hold for modules for a vertex operator
algebra $V$; in the statements, one replaces the vectors $v$ and $v'$
by elements $w$ and $w'$ of a $V$-module $W$ and its restricted dual
$W'$, respectively, and Proposition \ref{commandassocequivtoJacobi}
becomes: Given a vertex operator algebra $V$, in the presence of the
other axioms for the notion of $V$-module, the Jacobi identity follows
{from} the rationality of products and iterates, and commutativity and
associativity.  In particular, in the definition of $V$-module, the
Jacobi identity may be replaced by these properties.  }
\end{rema}

For either vertex operator algebras or modules, it is sometimes
convenient to express the equalities of rational functions in
Corollary \ref{dualitywithcovergence} informally as follows:
\begin{equation}\label{commutativityasoperatorvaluedratfns}
Y(v_{1}, z_{1})Y(v_{2}, z_{2}) \sim Y(v_{2}, z_{2})Y(v_{1}, z_{1})
\end{equation}
and
\begin{equation}\label{associativityasoperatorvaluedratfns}
Y(v_{1}, z_{1})Y(v_{2}, z_{2}) \sim Y(Y(v_{1}, z_1-z_2)v_{2},z_{2}),
\end{equation}
meaning that these expressions, defined in the domains indicated in
Corollary \ref{dualitywithcovergence} when the ``matrix coefficients''
of these expressions are taken as in this corollary, agree as
operator-valued rational functions, up to analytic continuation.

\begin{rema}\label{OPE}{\rm
Formulas (\ref{commutativityasoperatorvaluedratfns}) and
(\ref{associativityasoperatorvaluedratfns}) (or more precisely,
(\ref{associativitywithz1,z2})), express the meromorphic, or
single-valued, version of ``duality,'' in the language of conformal
field theory.  Formulas (\ref{associativityasoperatorvaluedratfns})
(and (\ref{associativitywithz1,z2})) express the existence and
associativity of the single-valued, or meromorphic, operator product
expansion.  This is the statement that the product of two (vertex)
operators can be expanded as a (suitable, convergent) infinite sum of
vertex operators, and that this sum can be expressed in the form of an
iterate of vertex operators, parametrized by the complex numbers
$z_1-z_2$ and $z_2$, in the format indicated; the infinite sum comes
{}from expanding $Y(v_{1}, z_1-z_2)v_{2},z_{2})$ in powers of $z_1-z_2$.
A central goal of this work is to generalize
(\ref{commutativityasoperatorvaluedratfns}) and
(\ref{associativityasoperatorvaluedratfns}), or more precisely,
(\ref{associativitywithz1,z2}), to logarithmic intertwining operators
in place of the operators $Y(\cdot,z)$.  This will give the existence
and also the associativity of the general, nonmeromorphic operator
product expansion.  This was done in the non-logarithmic setting in
\cite{tensor1}--\cite{tensor3} and \cite{tensor4}.  In the next
section, we shall develop the concept of logarithmic intertwining
operator.
}
\end{rema}


\bigskip

\noindent {\small \sc Department of Mathematics, Rutgers University,
Piscataway, NJ 08854 (permanent address)}

\noindent {\it and}

\noindent {\small \sc Beijing International Center for Mathematical Research,
Peking University, Beijing, China}

\noindent {\em E-mail address}: yzhuang@math.rutgers.edu

\vspace{1em}

\noindent {\small \sc Department of Mathematics, Rutgers University,
Piscataway, NJ 08854}

\noindent {\em E-mail address}: lepowsky@math.rutgers.edu

\vspace{1em}

\noindent {\small \sc Department of Mathematics, Rutgers University,
Piscataway, NJ 08854}

\noindent {\em E-mail address}: linzhang@math.rutgers.edu

\end{document}